\documentclass[10pt]{article}

\usepackage{url}
\usepackage{mathtools}

\usepackage{amssymb}
\usepackage{amsthm}
\usepackage{empheq}
\usepackage{latexsym}
\usepackage[inline]{enumitem}
\usepackage{eurosym}
\usepackage{dsfont}
\usepackage{appendix}
\usepackage{color} 
\usepackage[unicode]{hyperref}
\usepackage{frcursive}
\usepackage[utf8]{inputenc}
\usepackage[T1]{fontenc}
\usepackage{geometry}
\usepackage{multirow}
\usepackage{todonotes}
\usepackage{lmodern}
\usepackage{anyfontsize}
\usepackage{pgfplots}
\usepackage{stmaryrd}
\usepackage{cleveref}
\usepackage{natbib}
\usepackage{upgreek}
\usepackage{comment}
\usepackage{subfig}
\usepackage{scalerel}
\usepackage[normalem]{ulem} %strikethorugh
\usepackage{mathrsfs,bm,xfakebold,textgreek,upgreek}
\usepackage{cancel}
\usepackage{extarrows}

\RequirePackage[subfigure]{tocloft}
\advance\cftbeforesecskip 0.5em\relax
\usepackage[nottoc]{tocbibind}

\input Acorn.fd

\pgfplotsset{compat=1.16}

\setcitestyle{numbers,open={[},close={]}}

\definecolor{red}{rgb}{0.7,0.15,0.15}
\definecolor{green}{rgb}{0,0.5,0}
\definecolor{blue}{rgb}{0,0,0.7}
\hypersetup{colorlinks, linkcolor={red},citecolor={green}, urlcolor={blue}}
			
\makeatletter \@addtoreset{equation}{section}

\newtheorem{theorem}{Theorem}[section]
\newtheorem{assumption}{Assumption}

\newtheorem{corollary}[theorem]{Corollary}
\newtheorem{example}[theorem]{Example}

\newtheorem{lemma}[theorem]{Lemma}
\newtheorem{proposition}[theorem]{Proposition}

\newtheorem{definition}[theorem]{Definition}
\newtheorem{remark}[theorem]{Remark}
\newtheorem*{notation*}{Notation}

\DeclareUnicodeCharacter{014D}{\=o}
\def \C{\mathbb{C}}

\def \E{\mathbb{E}}
\def \F{\mathbb{F}}
\def \G{\mathbb{G}}

\def \L{\mathbb{L}}

\def \N{\mathbb{N}}
\def \P{\mathbb{P}}
\def \Q{\mathbb{Q}}
\def \R{\mathbb{R}}
\def \S{\mathbb{S}}

\def \X{\mathbb{X}}

\def\Bc{{\cal B}}
\def\Cc{{\cal C}}

\def\Ec{{\cal E}}
\def\Fc{{\cal F}}
\def\Gc{{\cal G}}
\def\Hc{{\cal H}}
\def\Ic{{\cal I}}
\def\Jc{{\cal J}}

\def\Lc{{\cal L}}
\def\Mc{{\cal M}}
\def\Nc{{\cal N}}

\def\Pc{{\cal P}}

\def\Rc{{\cal R}}

\def\Tc{{\cal T}}
\def\Uc{{\cal U}}
\def\Vc{{\cal V}}
\def\Wc{{\cal W}}
\def\Xc{{\cal X}}

\def\balpha{\boldsymbol{\alpha}}

\def\bx{\setBold[0.2] x \unsetBold}

\def\bmu{\setBold[0.2] \mu \unsetBold}

\def\bvarphi{\boldsymbol{\varphi}}

\def\eps{\varepsilon}
\def\d{{\mathrm{d}}}
\newcommand{\diff}{\,\mathrm{d}}
\newcommand{\1}[1]{\mathds{1}_{\{#1\}}}
\def\e{{\mathrm{e}}}

\def\=c{ \overset{c}{=}}
\def\Gammac{\Gamma_{\! c}}

\def\as{\text{--a.s.}}
\def\ae{\text{--a.e.}}
\renewcommand{\t}{\top}
 
\DeclareMathOperator*{\argmin}{arg\,min}

\newcommand{\fder}[1]{\frac{\delta #1}{\delta\mu}}

\newcommand{\smallerfont}[1]{\text{\fontsize{5}{5}\selectfont$#1$}} 
\newcommand{\smallfont}[1]{\text{\fontsize{4}{4}\selectfont$#1$}}

\newcommand{\xdim}{m}
\newcommand{\bmdim}{d}

\newcommand{\muin}{\nu_{0}}
\newcommand{\mufin}{\nu_{T}}

\renewcommand{\inf}{\mathop{\mathrm{inf}\vphantom{\mathrm{sup}}}}

\newcommand{\mybar}[1]{\makebox[0pt]{$\phantom{#1}\overline{\phantom{#1}}$}#1}

\definecolor{reviewteal}{rgb}{0.0,0.45,0.45}

\allowdisplaybreaks

\begin{document}

\title{Optimal transportation under McKean--Vlasov dynamics}%: duality, particle limits and Schrödinger potentials}

\author{Camilo {\sc Hern\'andez} \footnote{University of Southern California, ISE department, USA. camilohe@usc.edu.} \and Ludovic {\sc Tangpi} \footnote{Princeton University, ORFE department, USA. ludovic.tangpi@princeton.edu.}}

\date{\today}

\maketitle

\begin{abstract}
We study a broad class of stochastic transport problems under controlled McKean--Vlasov dynamics, allowing control of both drift and volatility, with weak terminal constraints. 
In contrast to the classical Schrödinger problem, the reference dynamics depend on the evolving law of the process, leading to nonlinear optimization problems on path space, while the terminal constraint is formulated through weak optimal transport, permitting a flexible treatment of distributional targets. 
	This framework encompasses several problems of current interest, including the mean field Schrödinger problem, martingale optimal transport, and Schrödinger--Bass bridges. 
	We establish three main results. 
	First, we prove the existence of optimal controls under general coercivity, continuity, and convexity assumptions. 
	Second, we derive a convex dual representation that extends weak Kantorovich duality and stochastic optimal transport duality to controlled McKean--Vlasov dynamics. 
	Finally, for the mean field Schrödinger problem, we prove convergence, at the level of value functions, of the associated finite-particle Schrödinger bridge problems to the mean field limit, thereby identifying the nonlinear problem as the exact infinite-particle limit of finite-particle Schrödinger bridges. 
	As a byproduct, we establish the existence of an optimal mean field Schrödinger potential, providing a variational foundation for duality in nonlinear stochastic transport.
\end{abstract}
\small{\textbf{Keywords}: Mean-field Schrödinger problem; McKean–Vlasov control; Weak optimal transport; Schrödinger bridge; propagation of chaos; interacting particle systems; convex duality.\smallskip

\textbf{2020 Mathematics Subject Classification.} Primary 49Q22; Secondary 49N80, 93E20, 60H10, 60K35.
}

%\tableofcontents

% {\color{blue}Comments
% \begin{itemize}
% %    \item abstract does not mention vol control.
%  %   \item definition of $\Vc_c$ in the intro, are we minimizing over $\P$ or $u$? or we don't want to say yet?
%     \item At least in the introduction, we should be consistent between $X \sim \mu$, $Law(X) = \mu$ and $\P\circ X^{-1} = \mu$.
% %    \item Why the first footnote? I don't understand. \camilo{deleted}
% \end{itemize}}

\section{Introduction}\label{sec:intro}

A broad class of stochastic transportation problems seeks to interpolate between prescribed probability distributions through the law of a stochastic process, while minimizing a cost associated with its dynamics; see, for instance, \citet*{tan2013optimal} and \citet*{liu2019compactness}.
 	A fundamental example is the Schrödinger problem, which asks for the most likely evolution of a random system conditionally on the observation of its initial and final laws. 
	In its modern formulation, given a reference path measure $\Rc$ on the path space $C([0,T];\R^\xdim)$ and two probability measures $\muin,\mufin$, one minimizes the relative entropy $\Hc(\P\,\vert\,\Rc)$ among path measures $\P$ satisfying $\P\circ X_0^{-1}=\muin$ and $\P\circ X_T^{-1}=\mufin$. 
	Originating in Schrödinger’s work \cite{schrodinger1931uber,schrodinger1932sur} and developed through the theory of reciprocal processes and entropy minimization, this problem has become a central object at the intersection of probability, stochastic analysis, and optimal transport; see, for instance, \citet*{follmer1988random}, \citet*{follmer1997entropy}, and \citet*{leonard2014survey,leonard2016lazy}.\medskip

In this work, we study stochastic transportation problems in which the dynamics of the transported process, as well as the associated transportation cost, may depend on its evolving path measure. 
	This law dependence introduces a genuinely nonlinear feature: the transport mechanism is itself determined by the law of the process being optimized. 
	A prominent instance of this structure is the mean field Schrödinger problem introduced by \citet*{backhoff2020mean} and further analyzed by \citet*{hernandez2023propagation}. 
	There, the fixed reference diffusion of the classical Schrödinger problem is replaced by McKean-Vlasov dynamics whose drift depends on the marginal flow of the candidate path measure, so that the reference measure against which relative entropy is computed varies with the path measure being optimized. 
	Mean-field Schrödinger bridge formulations have also been used in computational studies of crowd navigation and opinion dynamics; see \citet*{liu2022deep} and \citet*{inoue2026nonlocal}.
	We consider such law-dependent transportation problems subject to terminal constraints formulated via weak optimal transport, thereby accommodating distributional targets beyond equality of marginals. 
	%{\color{cyan}
	An additional difficulty, which will play an important role in our duality analysis, is that the natural McKean-Vlasov formulation lacks the convexity in the terminal distribution that is present in classical stochastic transport theory.
	%}
	\medskip
	
We now describe this problem more precisely, while leaving its rigorous formulation to the main body of the paper. 
	We say informally that $\P$ is an admissible law, on some abstract probability space, whenever there exists a control pair $(\alpha,\sigma\sigma^\top)$ such that the controlled process $X$ satisfies
\begin{align}\label{eq.dyn.intro}
    X_t = X_0 + \int_0^t \big(\alpha_s+b(s,X,\mu_s)\big)\,\d s
    + \int_0^t \sigma_s\,\d W_s,\; 
    \mu_t=\P\circ X_{\cdot\wedge t}^{-1},\;
    t\in [0, T],
\end{align}
where $W$ denotes a $\P$–Brownian motion. 
	Rather than requiring the terminal law to coincide with a prescribed target $\nu_T$, we impose the more general constraint $\P\circ X_T^{-1}\stackrel{c}{=}\nu_T$, formulated through weak optimal transport. More precisely, we write $\mu\stackrel{c}{=}\nu$ whenever $\Wc_c(\mu,\nu)=0$, where
\[
 \Wc_c(\mu,\nu)\coloneqq\inf_{\pi\in\Pi(\mu,\nu)}
 \int_{\R^\xdim}c(x,\pi_x)\,\mu(\d x),
 \]
$\Pi(\mu,\nu)$ denotes the set of couplings of $\mu$ and $\nu$, $(\pi_x)_x$ denotes the disintegration of $\pi$ with respect to its first marginal, and $c:\R^\xdim\times\Pc(\R^\xdim)\longrightarrow\R\cup{\infty}$ is a given weak transport cost.\medskip

Combining the controlled dynamics \eqref{eq.dyn.intro} with this terminal constraint, we arrive at the weak McKean-Vlasov transportation problem. 
	Namely, for a running cost $f$ that may depend, among others, on the controlled characteristic $u_t=(\alpha_t,\sigma_t\sigma_t^\top)$ and on the law $\mu_t$, we consider
\begin{align*}
\Vc_c(\nu_0,\nu_T)
\coloneqq
\inf\bigg\{
\E^{\P}\bigg[\int_0^T f(t,X,u_t,\mu_t) \diff t\bigg]:\;
\P\circ X_0^{-1}=\nu_0,\;
\P\circ X_T^{-1}\stackrel{c}{=}\nu_T
\bigg\},
\end{align*}
where the infimum is taken over all admissible laws $\P$.

\paragraph*{Examples.}
The formulation above encompasses several stochastic transportation problems of independent interest. We briefly discuss three examples that will also help place our results in context.

\begin{enumerate}[label=$(\roman*)$, ref=.$(\roman*)$,wide,  labelindent=0pt]
\item \emph{The martingale Benamou--Brenier problem.} 
	If $\muin$ and $\mufin$ are in convex order (i.e. $\langle \phi,\muin\rangle \le \langle \phi, \mufin\rangle$ for every $\phi:\R^\xdim\to \R$ convex with linear growth), one seeks a martingale $M$
%\[
%    M_t=M_0+\int_0^t \sigma_s\,\d B_s,
%    \qquad M_0\sim \muin,\quad M_T\sim \mufin,
%\]
minimizing a distance between the quadratic variation of $M$ and that of a reference martingale, for instance given $\varepsilon>0$, consider
\[
    \mathrm{MBB}(\muin,\mufin)
    \coloneqq
    \inf
    \bigg\{ 
    \E^{\P}\bigg[\int_0^T \|\sigma_t-I_m\|^{2 +\varepsilon} \diff t\bigg]
    :\;
    M_t=M_0+\int_0^t \sigma_s \diff B_s,
	\P\circ M_0^{-1}=\muin,\;
	\P\circ M_T^{-1}=\mufin
    \bigg\}.
\]
This problem was introduced in \citet*{backhoff2017martingale}, where its optimizer is characterized as the stretched Brownian motion. 
	Its connection with Bass martingales has since been revealed to have further variational and geometric structures; see \citet*{backhoff2023existence} and \citet*{backhoff2025geometric}.
	This problem is recovered from our formulation by taking the drift control equal to zero and choosing the weak cost
$c(x,\rho)=\int\1{x\neq y}\rho(\d y)$.
The weak transport formulation also allows one to relax the terminal constraint. For instance, the choice
$c(x,\rho)=\int_{\R^\xdim}|x-y|^2\rho(\d y)-\operatorname{Var}(\rho)/\theta $
yields a distributionally robust constraint: $\Wc_c(\P\circ X_T^{-1},\mufin)=0$ implies that $\P\circ X_T^{-1}$ belongs to a Wasserstein ball centered at $\mufin$, with radius $r$ satisfying $r^2=\operatorname{Var}(\mufin)/\theta$; see \citet[Section 3]{hernandez2025marginal} for this and other choices of weak transport costs.

\item \emph{The Schr\"odinger--Bass bridge.}
 In this case, one seeks a semimartingale
%\[
%    M_t=M_0+\int_0^t\alpha_s\d s + \int_0^t \sigma_s\,\d B_s,
%    \qquad M_0\sim \muin,\quad M_T\sim \mufin,
%\]
attaining
\[
    \mathrm{SBB}(\nu_0,\nu_T)
    \coloneqq
     \inf
    \bigg\{
    \E^{\P}\bigg[\int_0^T\frac{1}{2}\|\alpha_t\|^2 + \theta \|\sigma_t-I_m\|^{2 +\varepsilon}\diff t \bigg]:\; 
    M_t=M_0+\int_0^t\alpha_s\diff s + \int_0^t \sigma_s \diff B_s,\;
    \P\circ M_0^{-1}=\muin,\;
    \P\circ M_T^{-1}=\mufin
    \bigg\}.
\]
This interpolates between the classical Schrödinger bridge and the martingale Benamou–Brenier problem, balancing the cost of controlling the drift against that of controlling the volatility. 
	It is recovered from our formulation by taking $b=0$, choosing the running cost as above, and using the weak cost that enforces equality of terminal marginals. 
	The case $\varepsilon=0$, and its relevance to sampling and generative models, was recently studied by \citet*{henry2026bridging} and \citet*{hasenbichler2026weak}.

\item \emph{The mean field Schr\"odinger problem.}
Taking a fixed volatility $\sigma\in\R^{\xdim\times\xdim}$, a drift of the form $b(x,\mu)$, and the quadratic running cost for the drift control, our formulation yields the problem
\begin{align*}
\Vc_e(\muin,\mufin)
\coloneqq
\inf\bigg\{
\E^{\P}\bigg[\int_0^T\frac12\|\alpha_t\|^2\diff t\bigg]:
\P\circ X_0^{-1}=\muin,\;
\P\circ X_T^{-1}=\mufin
\bigg\},
\end{align*}
where the infimum is taken over admissible laws associated with these dynamics.\medskip

Equivalently, see \cite[Lemma 1.1]{backhoff2020mean}, and up to the entropy of the initial law, this problem admits the nonlinear entropy formulation
\[
\Vc_e(\muin,\mufin)
=
\inf\big\{
\Hc(\P\,\vert\,\Upsilon(\P)):
\P\in{\rm Prob}(\Omega),\;
\P\circ X_0^{-1}=\muin,\;
\P\circ X_T^{-1}=\mufin
\big\},
\]
where $\Upsilon(\P)$ denotes the law of
\begin{equation*}
	X_t = X_0 + \int_0^t b(X_s, \P\circ X_s^{-1} ) \diff s+ \sigma  W_t, \,t\in [0,T].
\end{equation*}	
Thus, unlike its classical counterpart, the Schrödinger problem for weakly interacting diffusions has a reference law $\Upsilon(\P)$ that depends on the path measure $\P$ being optimized.
\end{enumerate}
%{\color{red}how do we present this problem in the intro?}\camilo{I don't get it, because we do random? I think it is okay like this here since we clarify the distinction in the manuscript}

\paragraph*{Contributions and related literature.} The main contributions of this paper are threefold. 
	First, we prove the existence of optimizers for a broad class of weak McKean-Vlasov transportation problems, allowing for control of both the drift and volatility.
	Second, we establish a convex dual formulation extending weak Kantorovich duality and stochastic optimal transport duality to this nonlinear controlled setting. 
	Third, for the mean field Schrödinger problem described above, we show that the problem arises as the limit of finite-particle Schrödinger problems with weakly interacting particles. 
	As a byproduct of the finite-particle analysis, we also establish the existence of an optimal Schrödinger potential for the limiting mean field problem. 
\medskip

The first contribution, \Cref{thm:existence}, establishes the existence of optimizers for the weak McKean-Vlasov transportation problem with value $\Vc_c(\muin,\mufin)$. 
	This extends the theory of mean field Schrödinger problems developed in \cite{backhoff2020mean,hernandez2023propagation} in two directions: the characteristics of the state process may themselves be controlled, including its quadratic variation, and the terminal constraint is relaxed to the weak transport sense. 
	Under lower semicontinuity and convexity assumptions on the weak cost $c$, coercivity of the running cost in the controlled characteristics, and continuity and linear growth assumptions on the McKean-Vlasov drift, we show that the infimum is attained whenever the value is finite. 
	A notable feature of the result is that it is obtained in a weak formulation and requires no global Lipschitz condition on the drift.
	The proof relies on the compactness of the controlled laws and the stability of the associated martingale problems, together with the lower semicontinuity of the running cost and of the weak transport value function.
	\medskip
	
The second contribution, \Cref{thm:duality.gen}, is a convex dual representation of the weak McKean-Vlasov transportation problem. 
	The result admits a natural geometric interpretation. 
	The value, viewed as a function of its %initial and 
	terminal marginal, is convex and lower semicontinuous and is therefore determined by its supporting hyperplanes. 
	Our duality theorem identifies these hyperplanes explicitly in terms of weak Kantorovich potentials $\varphi$, while the weak $c$-transform
\[
    Q_c\varphi(x) \coloneqq \inf_{\rho } \big(c(x,\rho)+\langle \varphi, \rho\rangle \big),
\]
encodes the geometry of the terminal weak transport constraint. 
	We obtain
\begin{equation*}
    \Vc_c(\muin,\mufin)
    =
    \sup_{\varphi}
    \big(
    \Psi^\varphi(\muin) -
    \langle \varphi,\mufin\rangle
    \big),
\end{equation*}
where $\Psi^\varphi(\muin)$ is the value of an unconstrained McKean-Vlasov stochastic control problem with terminal cost $Q_c\varphi$. 
Thus, the dual potentials act as Lagrange multipliers for the terminal weak transport relation. 
The formula relies on the weak transport duality of \citet*{gozlan2017kantorovich} and \citet*{backhoff2019existence} and extends the stochastic optimal transport duality developed by \citet*{mikami2006duality} and subsequently in \cite{tan2013optimal}. 
%{\color{cyan}
A distinctive feature of our formulation is the use of exogenous randomization in the McKean-Vlasov law flow. This is essential for convex duality: as we show by an explicit example in \Cref{sub:convex_dual_representation}, the value function associated with deterministic law flows need not be convex in the terminal distribution.
Allowing randomized flows restores this convexity in the terminal distribution, and \Cref{thm:duality.gen} further shows that the resulting value is precisely the lower-semicontinuous convex envelope of the deterministic-flow value.
%}
In contrast with the latter setting, both the dynamics and the costs may depend on the law of the controlled process, while the terminal constraint is imposed in the sense of weak optimal transport.
\medskip

The third contribution, \Cref{thm:parti.conv.N}, is a finite-particle approximation theorem for the mean field Schrödinger problem $\Vc_e(\nu_0,\nu_T)$. 
	We consider $N$ weakly interacting controlled diffusions, started from an i.i.d.~configuration, whose drifts interact through the empirical measure and whose terminal marginal constraints are imposed particle by particle. 
	{See \Cref{assump.Df} for detailed structural assumptions.} 
	If $\Vc^N_e(\nu_0,\nu_T)$ denotes the corresponding $N$-particle cost, we prove the convergence of the cost per particle,
\[
\frac{1}{N}\Vc^N_e(\nu_0,\nu_T)
\longrightarrow
\Vc_e(\nu_0,\nu_T).
\]
	The theorem thus identifies the nonlinear mean field Schrödinger problem as the exact infinite-particle limit of finite-particle Schrödinger problems. 
	It can be viewed as a propagation of chaos result at the level of constrained variational problems, rather than solely at the level of the underlying dynamics. 
	A central difficulty is that the terminal constraint induces strong dependence among the particles at time $T$, while their controls are themselves coupled through the empirical measure.
	To the best of our knowledge, our analysis also provides the first {\color{black}self-contained} proof of the existence of an optimal Schrödinger potential for the mean field Schrödinger problem; see \Cref{cor:dual.exists.char}. 
	That is, there is a function $\varphi\in \L^2(\R^\xdim, \mufin)$ with a $\mufin$-square-integrable weak gradient (in the sense of distributions) such that 
	\[
	\Vc_e(\muin,\mufin) =  \Psi^\varphi(\muin) - \langle \varphi,\mufin\rangle.
	\]
	{\color{black}	In comparison, \cite[Section~1.3.2]{backhoff2020mean} discusses the existence of such a potential under additional assumptions on the so-called ``corrector’’ associated with the optimal mean field Schrödinger bridge.}
	The proof of \Cref{thm:parti.conv.N} relies essentially on the convex duality developed above and proceeds through matching lower and upper bounds. 
	The lower bound follows from compactness and entropy estimates, together with mean field limit results for controlled McKean-Vlasov systems from \citet*{djete2020mckean2,lacker2017limit,budhiraja2012large}. 
	The upper bound requires compactness of the finite-particle Schr\"odinger potentials. 
	To obtain it, we use tools from diffusion calculus, notably Bakry--\'Emery $\Gamma$-calculus, Fisher information estimates, and a Poincar\'e inequality for the target law. 
%	{
%	\color{cyan}
	A second structural ingredient is permutation symmetry.
	Although finite-dimensional duality yields particle-dependent potentials, the dual functional is permutation-invariant and concave; averaging over all labels can therefore only improve the dual value.
	Consequently, the limsup argument depends only on the empirical average of the particle potentials, for which the Sobolev estimates yield compactness.
	This allows us to identify a limiting Schrödinger potential and close the limsup inequality.
	Without the symmetrization step,
%	}
	these uniform Sobolev estimates for the potentials are reminiscent of the gradient estimates in \citet*{chiarini2023gradient}. 
	In this way, our result extends the mean field limit theory for controlled McKean-Vlasov systems \cite{djete2020mckean2,lacker2017limit} to the mean field Schrödinger problem with terminal marginal constraints, and complements the propagation of chaos results for constrained systems in \citet*{hernandez2023propagation}, where convexity assumptions on the target law are imposed. 
	The study of further properties of Schrödinger potentials under McKean-Vlasov dynamics is left for future work.
	In contrast to \Cref{thm:existence} and \Cref{thm:duality.gen}, which are proved for non-Markovian McKean-Vlasov dynamics with controlled volatility and possibly non-quadratic costs, the mean field limit result \Cref{thm:parti.conv.N} is obtained for the case of the mean field Schr\"odinger problem with Markovian and uniformly elliptic state dynamics, echoing that, for infinite-dimensional problems, dual attainment is more delicate than the dual representation itself.
	Indeed, the existence of a dual optimizer cannot be inferred from abstract convex analysis alone, since the usual interior conditions fail in spaces of measures, nor from the regularity of the coefficients. 
	The structure of the mean field Schrödinger problem provides, as previously alluded, exactly what is needed.
	See \Cref{sec:proof_of_the_finite_particle_approximation} for details.
%	 {\color{red}I think we need one or two sentences justifying this without sounding too defencive. This will allow to front run the only clear criticism they will have.}
	\medskip

The paper is organized as follows. 
Section \ref{sec.formulation} rigorously introduces the weak control rules, the weak optimal transport constraint, the value function, and the standing assumptions. 
It also states the main existence, convex duFality, and finite-particle approximation results. 
Section \ref{sec:proofs_of_existence_and_convex_duality_results} proves the existence theorem and the convex dual representation, and Section \ref{sec:proof_of_the_finite_particle_approximation} proves the finite-particle convergence theorem for the mean field Schr\"odinger problem.

\begin{notation*}
\begin{enumerate}[label=$(\roman*)$, ref=.$(\roman*)$,wide,  labelindent=0pt]
\item Throughout $T>0$ is fixed. For a nonnegative integer $\ell$, $\R^\ell$ is endowed with the Euclidean distance $d$, and $\lambda$ denotes the Lebesgue measure.
	$\S^\ell_+$ denotes the set of $\ell\times\ell$ symmetric positive semidefinite matrices.
	By $\Cc^\ell:=C([0,T],\R^\ell)$ we denote the space of continuous paths equipped with the sup norm.
	When $\ell=1$, we drop the superscript.
	$\nabla $ and $\nabla^2$ denote the gradient and Hessian operators, respectively, acting only on the space variable $x\in  \R^{\ell}$ of a smooth function $f:[0,T]\times \R^{\ell} \longrightarrow \R$.
	For a Polish space $E$, $C_b(E)$ denotes the space of bounded continuous functions.
	%For $\bmu\in \Pc(\Cc^\ell)$ a Borel probability measure on the canonical space $\Cc^\ell$ with canonical process $X$, we let $\bmu(t)=\bmu\circ X^{-1}_{t\wedge \cdot}$.
	For $p\geq 1$, we let $C_{p}(\R^\ell)$ denote the space continuous function of a mosth $p$-growth and
\(
C_{b,p}(\R^\ell) \coloneq \{\varphi\in \C_p(\R^\ell):   \exists a, b \in \R: a \leq \varphi(y)  \leq  b( 1+\| y\|^p),\, \forall y\in \R^\ell\}.
\)
	
\item For a metric space $\Omega$ and $\Fc$ its Borel $\sigma$-field, $\Pc(\Omega)$ denotes the set of probability measures on $\Omega$ equipped with weak topology. 
For $\P\in \Pc(\Omega)$ and a countably generated sub-$\sigma$-field $\Gc\subseteq\Fc$, $(\P^{\Gc}_\omega)_{\omega\in\Omega}$ denotes a family of regular conditional probability distribution (r.c.p.d.)~of $\P$ knowing $\Gc$.
Given $\Q\in\Pc(\Omega)$, we define the concatenation $\Q\otimes \P^{\Gc}_\cdot (A) \coloneqq \int_\Omega \P^{\Gc}_\omega (A) \Q(\d \omega), A\in \Fc$.
	For $\P,\Q\in \Pc(\R^\ell)$, $\Hc(\Q| \P)$ denotes the relative entropy of $\Q$ with respect to $\P$, given by 
 $\E^\Q\big[ \log \frac{\d \Q}{\d \P}\big]$, if $\Q\ll\P$ and $\infty$ otherwise.
For $\mu\in\Pc(\R^\ell )$, and $\varphi:\R^\ell \longrightarrow \R$, we set $\langle \varphi,\mu\rangle \coloneqq \int_{\R^\ell} \varphi (x)\mu(\d x)$.
In particular, $M_p(\mu)\coloneqq \int_{\R^\ell}\|x\|^p\mu(\d x)$ denotes the $p^{\mathrm{th}}$-moment of $\mu$.
$\Pi(\mu,\nu)$ denotes the set of all couplings $\pi\in \Pc(\R^\ell\times \R^\ell)$ with marginals $\mu$ and $\nu$ and $(\pi_x)_{x\in X}$ denotes the {$\mu\ae$~unique} regular disintegration with respect to the first marginal.
For $p\geq 1$, $\Pc_p(\R^\ell)$ denotes the family of $\mu\in\Pc(\R^\ell)$ with finite $p^{\mathrm{th}}$-moment endowed with the the $p$-Wasserstein distance
\(
W_p^p (\mu,\nu)\coloneqq \inf_{\pi\in \Pi(\mu,\nu)}\int_{\R^\ell\times \R^\ell} \|x-y\|^p  \pi(\d x, \d y).
\)
%\medskip

\item 
%\begin{definition}\label{def:C12}
$v:[0,T]\times\Pc_2(\R^\xdim)\longrightarrow \R$ belongs to
$\mathcal C^{1,2}_{\rm flat}$ if:
\begin{enumerate*}[label=$(\alph*)$]
\item $v$ is continuous and, for every $\mu \in \Pc(\R^\xdim)$, $t\longmapsto v(t ,\mu)$ is differentiable and $\partial_tv$ is continuous on $[0,T]\times\Pc_2(\R^\xdim)$;
\item there is a jointly continuous
$u:[0,T]\times\Pc_2(\R^\xdim)\times\R^\xdim\longrightarrow \R$ such that $x\longmapsto u(t,\mu)(x)$ is
$C^2$ for every $(t,\mu)\in [0,T]\times\Pc(\R^\xdim)$, $\nabla_xu$ and $\nabla^2_xu$ are jointly continuous, and for all
$\mu,\mu'\in\Pc_2(\R^d)$ we have
\(
v(t,\mu')-v(t,\mu)=\int_0^1\!\!\int_{\R^d}
u\big(t,\lambda\mu'+(1-\lambda)\mu\big)(x)\,(\mu'-\mu)(\d x)\,\d\lambda;
\)
\item $\nabla_x u$ and
$\nabla^2_x u$ are bounded. %, with $M_2(\mu):=\int_{\R^\xdim}|x|^2\mu(dx)$.
\end{enumerate*}
%In the L-derivative notation of [CD-I, Chapter~5],
%$\partial_\mu v(t,\mu)(x)=\nabla_xu(t,\mu)(x)$ and
%$\partial_x\partial_\mu v=\nabla^2_xu$; (b)--(c) say that $v(t,\cdot)$
%is \emph{partially $C^2$} in the sense of [CD-I, Definition~5.83],
%uniformly in $t$, with the growth conditions of [CD-I, Section~5.6].
%The normalization $\int u(t,\mu)\,\d\mu=0$ may be imposed; it plays no
%role.
%\end{definition}
The function $u$ is unique up to an additive constant, and thus, we normalize it by
%We remind the reader that flat derivatives are normalized by 
$\int u(t,\mu)(x) \mu(\d x)=0$.
We call $u$ the flat derivative of $v$, denoted by $\fder v(t,\mu)(x) \coloneqq u(t,\mu)(x)$, or simply $\delta v$, see, e.g., {\rm \cite{carmona2018probabilisticI,carmona2018probabilisticII}} for further details and properties.

\item $v:[0,T]\times\Pc_2(\R^\xdim)\longrightarrow \R$ belongs to
$\mathcal C^{1,2,{\rm full} }_{\rm flat}$ if: 
\begin{enumerate*}[label=$(\alph*)$]
\item $v\in\mathcal C^{1,2}_{\rm flat}$ and, for every $(x,\mu) \in \R^\xdim\times \Pc(\R^\xdim)$, $t\longmapsto \delta v(t ,\mu)(x)$ is differentiable and $\partial_t \delta v$ is continuous on $[0,T]\times \R^\xdim\times \Pc(\R^\xdim)$;
\item $\delta v(t,\mu)(y)$ is flat-differentiable, i.e.,
\(
\delta^2v(t,\mu)(y;x)\!\coloneqq \frac{\delta}{\delta\mu}\big(\delta u(t,\cdot)(y)\big)(\mu)(x)
\);
%,
%\(
%(t,\mu,y,x)\in[0,T]\times\Pc_2(\R^\xdim)\times\R^\xdim\times\R^\xdim,
%\)
\item $(x,y)\mapsto \delta^2v(t,\mu)(y;x)$ is jointly continuous and $C^2$ in $y$ and in $x$ separately;
\item $\psi\in \{\nabla_y\delta^2 v,\nabla^2_y\delta^2 v, \nabla_x\delta^2 v, \nabla^2_x\delta^2 v\}$ is jointly continuous and bounded; % of growth $|\psi(x)|\le C(1+|x|+|y|+M_2(\mu))$; 
\item $\partial_t \delta v $, $\nabla_y \delta v$, $\Delta_y \delta v$ are flat-differentiable with flat derivatives $\partial_t\delta^2v$, $\nabla_y\delta^2v$, $\Delta_y\delta^2v$.
\end{enumerate*}
\end{enumerate}
\end{notation*}

%%%%%%%%%%%%%%%%%%%%%%%%%%%%%%%%%%
%%%%%%%%%%%%%%%%%%%%%%%%%%%%%%%%%%
%%%%%%%%%%%%%%%%%%%%%%%%%%%%%%%%%%
%%%%%%%%%%%%%%%%%%%%%%%%%%%%%%%%%%

\section{Problem formulation and main results}\label{sec.formulation}

We now rigorously present the {path-dependent McKean-Vlasov optimal transportation} problem considered in this article and state the main results.
Let $\xdim$ and $\bmdim$ be two positive integers and let $\Uc\subseteq \R^\xdim\times \S^{\xdim}_+$ be a nonempty closed set containing $0$.
	We will consider two mappings
\begin{align}\label{eq.data}
\begin{split}
b:[0,T]\times\Cc^\xdim\times \Pc(\Cc^\xdim)\longrightarrow \R^\xdim,\; f:[0,T]\times\Cc^\xdim\times \Uc \times \Pc(\Cc^\xdim)\longrightarrow \R,
\end{split}
\end{align}
assumed to be Borel-measurable and nonanticipative in the sense that for all $(t, x ,u ,   \mu) \in [0,T]\times\Cc^\xdim\times \Uc \times \Pc(\Cc^\xdim)$
\begin{align*}
f(t,x, u, \mu):=f(t,x _{t\wedge \cdot} , u ,  \mu(t)) \mbox{ and } b(t, x ,  \mu):=b(t, x_{t\wedge \cdot}  ,  \mu(t)) \mbox{ where } \mu(t) = \mu\circ x_{\cdot\wedge t}^{-1}.
\end{align*}

To streamline the presentation, we defer all comments on the formulation to \Cref{rem:defcomments}.

\paragraph*{Controlled McKean-Vlasov dynamics.}
%\camilo{Note that in Mikami's paper the probability space is never fixed}
	Formally, in the spirit of \cite{mikami2006duality} we will consider the class of continuous semimartingales $X$, defined on some filtered probability space carrying a $\bmdim$--dimensional Brownian motion $W$, for which there is a $\Uc$--valued measurable process $u=(\alpha,\sigma\sigma^\t)$, %\footnote{For any $a\in \S^m_+$ we let $\sigma=a^{1/2}$} 
	such that
\begin{equation}\label{eq.dyn}
%\label{eq:control.SDE.MF}
	``X_t^u = X_0^u + \int_0^t\big( \alpha_r+b(r, X^u, \mathcal{L}(X^u))\big) \diff r + \int_0^t  \sigma_r  \diff W_r ,\; t\in [0,T]".
\end{equation}
Rigorously, the above dynamics will be specified in the weak sense as in the following definition:%\medskip

\begin{definition}\label{def.controlrule.weak}
Let $p\geq 1$ and $\nu\in \Pc(\R^\xdim)$. A weak control rule with initial condition $\nu$ consists of a tuple
\[
\gamma:=\big(\Omega^\gamma,\Fc^\gamma,\F^\gamma=(\Fc^\gamma_t)_{t\in [0,T]},\G^\gamma=(\Gc^\gamma_t)_{t\in [0,T]}, \P^\gamma, X^\gamma,W^\gamma,\mu^\gamma,  u^\gamma\big)
\]
satisfying,
\begin{enumerate}[label=$(\roman*)$, ref=.$(\roman*)$,wide,  labelindent=0pt]
\item \label{def.controlrule.weak.i} $(\Omega^\gamma,\Fc^\gamma,\P^\gamma)$ is a probability space and $\F^\gamma,\G^\gamma$ are filtrations satisfying,
%
% , $t\in[0,T]$, satisfying the $(H)$--hypothesis type condition
\begin{equation*}%\label{eq:Hcond}
 \Gc^\gamma_t\subseteq\Fc^\gamma_t,\; \mbox{and }\; \E^{\P^\gamma} [\mathbf 1_D\, |\,\Gc^\gamma_t] =\E^{\P^\gamma} [\mathbf 1_D\, |\,\Gc^\gamma_T],\; \P^\gamma\as \mbox{ for all } D\in\mathcal F^\gamma_t\vee\sigma(W^\gamma),\ t\in[0,T];
\end{equation*}

\item \label{def.controlrule.weak.ii} $X^\gamma=(X^\gamma_t)_{t\in[0,T]}$ is an $\R^\xdim$--valued, $\F^\gamma$--adapted continuous process, $u^\gamma=(u^\gamma_t)_{t\in[0,T]},\, u_t^\gamma=(\alpha^\gamma_t,\sigma^\gamma_t(\sigma^\gamma_t)^{\!\t}), %\camilo{(\alpha_t^\gamma,\beta^\gamma_t), \, (\beta_t^\gamma)^{1/2}=\sigma^\gamma_t }, 
t\in [0,T]$ is an $\Uc$--valued, $\F^\gamma$--predictable process, and
\[
	\E^{\P^\gamma}\bigg [\sup_{t\in [0,T]} \|X^\gamma_t\|^p\bigg]+\E^{\P^\gamma}\bigg[\int_0^T  \| u_t^\gamma\|^p \diff t\bigg]<+\infty;%, \text{ for some }{\bu}_0\in U.
\]
\item \label{def.controlrule.weak.iii}
$W^\gamma=(W^\gamma_t)_{t\in[0,T]}$ is an $\R^\bmdim$--valued, $(\F^\gamma,\P^\gamma)$--Brownian motion, 
	$\Fc^\gamma_0\vee\sigma(W^\gamma)$ is $\P^\gamma-$independent of $\Gc^\gamma_T$,
	 and $\mu^\gamma=(\mu^\gamma_t)_{t\in[0,T]}$ is a $\Pc(\Cc^\xdim)$--valued, $\G^\gamma$--predictable process such that
\[
	\mu^\gamma_t  	=
	\Lc^{\P^\gamma}(X^\gamma_{t\wedge\cdot}\, |\,\Gc^\gamma_t\big),
	\d t\otimes\d \P^\gamma\ae\mbox{~on }[0,T]\times\Omega^\gamma;
\]

\item \label{def.controlrule.weak.iv} $ \P^\gamma\circ (X^\gamma_0)^{-1}=\nu$, and it holds that
\[
X^\gamma_t = X_0^\gamma + \int_0^t\big( \alpha_r^\gamma+b(r, X^\gamma,\mu^\gamma_r) \big) \diff r + \int_0^t  \sigma_r^\gamma  \diff W_r^\gamma ,\; t\in [0,T],\, \P^\gamma\as
\]
\end{enumerate}
We let 
$$\Gamma(\nu):=\{\gamma: \gamma\text{ is a weak control rule with initial condition $\nu$} \}.$$
\end{definition}
%{\color{cyan}
Observe that the conditions on $\Gc^\gamma$ prevent the exogenous randomness from revealing future information about $\Fc^\gamma$.\medskip
%}

To alleviate the notation, and following ideas in {\rm \cite{djete2020mckean2}}, in our analysis, we will represent a $\Uc$--valued process $u^\gamma$ by means of an $\R$-valued process $U^\gamma$.
Indeed, since $\Uc$ is a Polish space, it is isomorphic to a Borel subset of $[0,1]$, and thus we denote by $\pi$ an isomorphic bijection between $\Uc$ and $\pi(\Uc)\subset[0,1]$.
Thus, on $[0,T]\times\Omega^\gamma$ we define the process
\begin{align}\label{eq.repofcontrols}
    U_\cdot^\gamma := \int_0^\cdot \pi(u_s^\gamma) \diff s, \text{ and note that,} \;  u_t^\gamma=\pi^{-1}\Big(\limsup_{n\longrightarrow\infty} n\big(U_t^\gamma-U_{(t-1/n)\vee 0}^\gamma\big)\Big), \, \d t\otimes\d \P^\gamma\ae\text{ on }[0,T]\times\Omega^\gamma.
\end{align}
The advantage of this definition is twofold: it will avoid introducing measure-valued processes on the canonical space of measures on $\Uc\times[0,T]$ %{\color{blue}what is $U$?} 
when we introduce the canonical formulation; moreover, it simplifies the notation as we keep track of the control pair via a single canonical process.

%{\color{blue}Names with references or no names with references?} \camilo{ Camilo: Names for everyone the first time they appear, then numbers}

\paragraph*{The weak optimal transportation problem under McKean-Vlasov dynamics.}
In this work, we will consider transportation problems for which the terminal configuration is specified in a broader sense.
	This section presents the weak optimal transport problem (WOT), first introduced in \cite{gozlan2017kantorovich}.\medskip
	%Further details, results relevant to our analysis, and examples are deferred to \Cref{sec:examples}.

Let $\mu,\nu\in \Pc(\R^\xdim)$. 
Given a suitable cost functional $c:\R^\xdim\times\Pc(\R^\xdim)\longrightarrow \R\cup\{+\infty\}$, see \Cref{assump.costs} below.
%assumed to be bounded from below, lower semicontinuous (l.s.c. for short) and with $c(x,\cdot)$ convex, see \Cref{assump.costs} below, WOT is given by
\begin{align}\label{eq.weaktransport}
 \Wc_c (\mu,\nu):=\inf_{\pi\in \Pi(\mu,\nu)}\int_{\R^\xdim} c(x,\pi_x)\mu(\d x).
\end{align}
%Whenever we wish to stress the cost functional $c$, we write $\Wc_c$ which 
	%We refer to $\Wc_c$ as the \emph{{\rm WOT} induced by $($the cost function$)$ $c$}, or the $c$-transport cost.
	Notice that taking $c(x, \rho) =\langle \mathfrak{c}(x,\cdot), \rho\rangle$, for $\mathfrak{c}(x,y)$ a \emph{classical} cost leads to the classical optimal transport problem. 
	By a classical cost, we mean that $\textup{c}:\R^\xdim\times\R^\xdim\longrightarrow \R\cup\{+\infty\}$ is l.s.c.~and bounded from below. 
	In particular, with $\mathfrak{c}(x,y)= \|x-y \|^p $, for $p\geq 1$, we arrive at the $p$-Wasserstein distance $W_p$.\medskip
%\begin{align*}
%W_p^p (\mu,\nu)\coloneqq \inf_{\pi\in \Pi(\mu,\nu)}\int_{\R^\xdim\times \R^\xdim} \|x-y\|^p  \pi(\d x, \d y).
%\end{align*}

	Recalling that, for $p\geq 1$, $W_p$ is a metric on $\Pc_p(\R^\xdim)$, i.e. $\mu=\nu \Longleftrightarrow  W_p(\mu,\nu)=0$, it seems natural to introduce the binary relation induced by $\Wc_c$ over the elements of $\Pc(\R^\xdim)$.
	%, and consequently over the family of random variables with given distributions. 

\begin{definition}\label{def.msim}
Let $\mu,\nu\in\Pc(\R^\xdim)$. We say $\mu \=c \nu \Longleftrightarrow  \Wc_c(\mu,\nu)=0$. %For a r.v.~$\xi$ such that $\xi \sim\mu$, we write $\xi \msim\nu\Longleftrightarrow  \Wc_c(\mu,\nu) $.
\end{definition}

Some relevant instances of \eqref{eq.weaktransport} and their associated relation $\overset{c}{=}$ are given, e.g.,~in \cite{hernandez2025marginal}. 
	%These would correspond to different specifications of cost functional $c$. 
	For now, to make \eqref{eq.weaktransport} precise and, more importantly, take values in $\R$, we introduce the standing assumptions on the cost functionals we will be considering.

%\camilo{In \cite[Page 3]{backhoff2019existence}: regarding the Polish space $X$, we fix from now on a compatible \emph{bounded} metric $d_X$.}

\begin{assumption}\label{assump.costs}
We say $c:\R^\xdim\times \Pc_p(\R^\xdim)\longrightarrow \R_+$ satisfies {\rm \Cref{assump.costs}} with $p\geq1$ if:
\begin{enumerate}[label=$(\roman*)$, ref=.$(\roman*)$,wide,  labelindent=0pt]
\item \label{assump.costs.1} It is jointly l.s.c.~with respect to the product topology of $\R^\xdim\times\Pc_p(\R^\xdim)$ and $c(x,\delta_x)=0$.
%\footnote{$\Pc_p(\R^\xdim)$ being equipped with the topology generated by $W_p$.} 
\item \label{assump.costs.2} The map $\rho\longmapsto c(x,\rho)$ is convex for any $x\in \R^\xdim$ in the sense that
\[
c(x,\theta \rho_1+(1-\theta )\rho_2)\leq \theta  c(x,\rho_1)+(1-\theta  )c(x,\rho_2),\,  \rho_1,\rho_2\in \Pc_p(\R^\xdim),\, \theta \in [0,1].
\]
\item \label{assump.costs.3}  There is $C>0$ such that
\[
c(x,\rho)\leq C\big(1+\|x\|^p +M_p(\rho) \big), \, (x,\rho)\in \R^\xdim\times\Pc_p(\R^\xdim).
\]
\end{enumerate}
\end{assumption}

%and
%{\color{cyan}
%For this, we introduce the space
%
%}
%{\color{blue}Do we need the cyan stuff here? I propose to move it to where we first need it.}

Let us introduce the set
\begin{align*}
\Gammac(\muin,\mufin)\coloneqq\big\{ \gamma \in \Gamma(\muin): \P^\gamma \circ (X_T^\gamma)^{-1}\=c\mufin \big\}.
\end{align*}
%{\color{blue}
% \begin{itemize}
% \item careful about the case $\Pc(\mu_0,\mu_T) = \emptyset$
% %\item I find $(\mu_0,\mu_T)$ cumbersome. Can we use $(\mu,\nu)$ and the dummy variable on $\Pc$ be $m$? then the dimension be $d$? \camilo{We will probably save $\mu$ for the law of $X$. You can use $\textbackslash $muin and $\textbackslash $mufin for $\muin, \mufin$ :)}
% %\item I find the usage of the Greek letter $\nu$ below possibly misleading, given that $\mu$ is a measure, $\nu$ might seem like a measure too. how about using $u = (\alpha,\beta)$? \camilo{Agree. $\nu$ is now used for the data of the transportation problem}
% \end{itemize}
%}
We are interested in the McKean-Vlasov optimal transportation problem\footnote{As usual, we use the convention $\inf\emptyset \coloneqq + \infty$.}
\begin{equation}
\label{eq:problem}
	V_c(\muin, \mufin)
	\coloneqq
	\inf_{\gamma \in \Gammac(\muin,\mufin)} J(\gamma), 
	\; \text{where } 
	J(\gamma)\coloneqq \E^{\P^\gamma}\bigg[ \int_0^T f(t,X^\gamma,u^\gamma_t , \mu^\gamma_t )\diff t\bigg].
\end{equation}

\begin{remark}\label{rem:defcomments}
%\begin{enumerate}[label=$(\roman*)$, ref=.$(\roman*)$,wide,  labelindent=0pt]
%\item  
Under the weak formulation in {\rm \Cref{def.controlrule.weak}}, following {\rm \cite[Definition~2.1]{djete2019mckean}}, the coefficients are evaluated along a possibly random process $\mu^\gamma$.
	The special case in which $\G^\gamma$ is trivial is compatible with {\rm \Cref{def.controlrule.weak}}.
	In that scenario {\rm \Cref{def.controlrule.weak}\ref{def.controlrule.weak.i}} and the independence in {\rm \Cref{def.controlrule.weak}\ref{def.controlrule.weak.iii}} are vacuous and $\Lc(X^\gamma_{t\wedge\cdot}|\Gc^\gamma_t) = \P^\gamma\circ(X^\gamma_{t\wedge\cdot})^{-1}$.
	We refer to rules of this special form as deterministic-flow rules and write $\Gamma^{\rm det}(\muin)$ for their collection and, analogously, we define $\Gammac^{\rm det}(\muin,\mufin)$. With this, we may also consider the McKean-Vlasov transportation problem
	\begin{align}\label{eq:problem.det}
	V_c^{\rm det}(\muin, \mufin)\coloneqq
	\inf_{\gamma \in \Gammac^{\rm det}(\muin,\mufin)} J(\gamma).
	\end{align}
% \item The reader familiar with the literature on stochastic transportation problems might wonder about the need to consider $\Gamma(\muin)$ as the class of weak control rules for the case of transportation via McKean-Vlasov dynamics.
% 	In anticipation of our discussion below, let us note that, in general, the map $\mufin\longmapsto V_c^{\rm det}(\muin, \mufin)$ is not convex; see {\rm \Cref{sec.conx.fails}}.
% 	Ultimately, the exogenous randomization carried into $\G$ endows $\mufin\longmapsto V_c(\muin, \mufin)$ with the convexity needed for our duality analysis, see {\rm \Cref{prop.convex.V}}
%\end{enumerate}
\end{remark}

% {\color{blue}are we crazy enough to consider the problem
% \begin{equation*}
% 	V_c(\muin, \mufin) := \inf_{\gamma \in \Gamma(\muin)} J(\gamma), \; \text{where }  J(\gamma):=\E^{\P^\gamma}\bigg[ \int_0^T f(t,X^\gamma,u^\gamma_t , \mu^\gamma_t )\diff t\bigg] + g(\Wc_c(\P^\gamma\circ (X^\gamma_T)^{-1}, \mufin)).
% \end{equation*}
% for some convex, lsc and proper function $g:\R\to \R\cup\{\infty\}$?}
%
%
%
%
% and $\P$ denotes the unique weak solution to the McKean--Vlasov SDE {\color{blue}double check below. Is this the correct equation? don't we want $\Lc_\alpha(X_s) := \P^{\alpha}\circ X_s^{-1}$?}
% \begin{equation}\label{eq.mkvsde.convergence}
% 	X_t = X_0 + \int_0^t b(X_s, \Lc(X_s)) \d s+ \int_0^t\sigma_s\d W_s, t\in [0,T],\P\as
% \end{equation}	
% Let us put $\P\circ X_0^{-1}= \mu$.
% % where $\bar\mu\in \Pc(\R^\xdim)$ denotes its unique stationary distribution, and a $\P$--Brownian motion $W$.
% Then it follows by Girsanov's theorem that $V$ is the value of the entropy minimization problem
% \begin{equation*}
% 	V(\muin, \mufin) = \inf\big\{ \Hc(\Q|\P): \Q\ll \P, \;\; \Q\circ X_0^{-1} = \muin, \, \Q\circ X_T^{-1} = \mufin\big\}  - \Hc(\muin|\mu).
% \end{equation*}
%

\subsection{Main results}
We will now present the assumptions made on the data and state our main results, which are existence and duality for Problem \eqref{eq:problem}, and propagation of chaos in the case of the {Markovian} mean field Schr\"odinger problem discussed in the introduction.
	Furthermore, for the {Markovian} mean field Schr\"odinger problem, we will prove the existence of an optimal Sch\"odinger potential.
We assume the following conditions on the data $\muin,\mufin$ and the functions $b$ and $f$:

\begin{assumption}\label{assump.bf}
{Let $\muin,\mufin\in \Pc_p(\R^\xdim)$} and the functions $b$ and $f$ satisfy
	%The functions $b$ and $f$ satisfy {\rm \Cref{assump.bf}} if
	\begin{enumerate}[label=$(\roman*)$, ref=.$(\roman*)$,wide,  labelindent=0pt]
	\item \label{assump.bf.i}  {For every $t \in [0,T]$ the function $({x}, {u}, \mu)\longmapsto f(t, x , u,\mu)$ is continuous, nonnegative and convex in u {\rm ;} and $f$ is uniformly continuous in $t$.
 That is,
 \begin{align*}
 	\Delta_tf(\varepsilon) := \sup \frac{|f(t, x,u,\mu) - f(s, x, u,\mu)|}{1 + f(s, x, u,\mu)} \longrightarrow 0, \mbox{ as }  \eps\longrightarrow 0,  
 \end{align*}
where the supremum is taken over all $(t,s)\in [0,T]^2$, such that $|t - s|\le \varepsilon$ and $(x, u , \mu) \in  \Cc^{ \xdim}\times \Uc\times \Pc(\Cc^\xdim)$;}

% \todo[inline]{In \cite{mikami2002optimal, tan2013optimal} the limit uniform in $(x,u,\mu)$. One takes the sup over $(x,u)$ and I guess $\mu$ for us}
% \camilo{we dont use convexity of $f$ to get convexity of $V$.}
	\item \label{assump.bf.2} There is $\hat p > 1$ and $C_1 \in \R, C_2>0$, such that for all $(t, x, u , \mu) \in [0,T]\times \Cc^\xdim \times \Uc\times \Pc(\Cc^\xdim)$
	\[
	    f(t,x,u,\mu) \geq C_1+C_2\| u\|^{\hat p};
	\]\vspace{-17pt}
	\item \label{assump.bf.3} The function $(t,x, \mu)\longmapsto b(t,x,\mu)$ is measurable and continuous in $\mu$, and 
% {\color{blue}we probably don't need continuity in $\bf x$}
 there is $C>0$ such that for all $(t, x,\mu) \in [0,T]\times \Cc\times \Pc(\Cc^\xdim)$,
%{\color{blue}shouldn't we assume the bound in terms of $(\int_{\Cc^\xdim}\|\bx^\prime\|^2 \bmu(d\bx^\prime))^{1/2}$?}
%\camilo{Make sure that $\|u\|^p<C$ guarantees $\|\sigma\sigma^T\|^p<C^\prime$. finite-dimensional space is enough, no?}
	\[
	\| b(t,x,\mu) \|\leq C\big (\|x\|+ M_2(\mu)^{1/2} \big).
	\]
	\end{enumerate} 
\end{assumption}

\begin{remark}
{\rm Assumption \ref{assump.bf}} collects standard structural conditions ensuring compactness and stability of admissible weak control rules, while remaining flexible enough to accommodate McKean-Vlasov interactions. 
	The convexity of $f$, its continuity in $( x,u,\mu)$, together with its uniform continuity in time, are tailored to guarantee lower semicontinuity of the cost functional under weak convergence of controls and states, which is essential for the existence argument.
	The coercivity condition in {\rm$(ii)$} provides uniform integrability of the control processes and plays a key role in establishing tightness of admissible laws. 
	Finally, the linear growth condition on the drift $b$ in {\rm $(iii)$}, imposed only in terms of the current state and the marginal flow, is sufficient to ensure tightness of the state process.
	Notably, we do not impose global Lipschitz assumptions on $b$, as these are not needed for the variational analysis developed here;\footnote{Indeed, these are not necessary conditions for the existence of a weak control rule, see, for instance, \cite[Assumption S]{carmona2015probabilistic} or \cite[Assumption B2]{lacker2018strong}.}~instead, the formulation relies on weak solutions and martingale problem techniques, which allow for greater generality.
\end{remark}

The results of this section are organized as follows. 
We first state, in \Cref{thm:existence}, the existence of optimal controls, followed by a convex dual representation in \Cref{thm:duality.gen}.
These results follow from {\rm Assumptions} {\rm\ref{assump.costs}} and {\rm\ref{assump.bf}}.
Then, after introducing the mean field entropic Schr\"odinger problem, we introduce its associated finite-particle approximation and state convergence to the mean field problem, see \Cref{thm:parti.conv.N}.
This result will require more structure on the problem, see \Cref{assump.Df}.

%???} {\color{cyan}to be fixed when we decide where to put the example}
%{\color{red}In our generality (because $\sigma$ is controlled) we should have $\Pc(\muin,\mufin)\neq\emptyset$. Let's check this.}
%\camilo{Can we take $\alpha_t=-b(t,X_t,\mu_t)$ and then do Skorokhod embedding?}
\begin{theorem}
\label{thm:existence}
	Let {\rm Assumptions} {\rm\ref{assump.costs}} and {\rm \ref{assump.bf}} hold and let $(\muin,\mufin)$ be such that $V_c(\muin,\mufin)<+\infty$.
    Then, there exists $\gamma \in {\Gamma_c}(\muin,\mufin)$ such that $V_c(\muin,\mufin) = J(\gamma)$.
\end{theorem}
\begin{remark}\label{rmk.exist.unique}
A key feature of this result is that no Lipschitz condition is imposed on the drift. The proof relies entirely on the compactness of controlled laws and the stability of the associated martingale problems, rather than on pathwise arguments.\medskip

Let us also remark that the uniqueness of the optimizer remains a much more subtle question. 
	In particular, in the above weak formulation, the problem may have multiple solutions even if we strengthen our assumptions on $b$. 
	In fact, even for specific instances of $V_c(\muin,\mufin)$, such as the mean field Schrödinger problem described in the introduction, uniqueness remains open in general, owing to the nonlinear nature of the underlying entropy minimization problem, see {\rm \cite[Remark 1.1]{backhoff2020mean}}.
\end{remark}

   % \camilo{It is open since the control $\P$ also determines the reference measure in the entropy minimization formulation $\P\mapsto H(\P\vert\Upsilon (\P))$. For the control formulation, the only case I could think of was the (boring) linear in the law case, $b({\rm Law}(X_t))=\E[X_t]$.}\\
    
%    
%{\color{blue}I would propose to move the first cyan text here.}\camilo{The definitions of $C_{b,p}$ and $Q_c\varphi$ are used in the dual formula of the weak transport problem. To me, it makes sense to keep the paragraph section ``the weak transport problem'' self-contained with the results we are gonna cite.}\\

In addition to this existence result, we have the following convex dual representation.
To state this dual representation, for any $\varphi\in C_{b,p}(\R^\xdim)$, we introduce the \emph{weak $c$-transform} of a real-valued function $\varphi$, given by
\begin{align*}%\label{eq.cconjugate}
    Q_c\varphi(x)\coloneqq \inf_{\rho \in \Pc_p(\R^\xdim)} \big( c(x,\rho)+ \langle \varphi,\rho\rangle\big).
\end{align*}

\begin{theorem}
\label{thm:duality.gen}
	If {\rm Assumptions }{\rm \ref{assump.costs}} and {\rm \ref{assump.bf}} are satisfied, then it holds
	\begin{equation*}
		V_c(\muin,\mufin) 
		= 
		\sup_{\varphi \in C_{b,p}(\R^\xdim)}
		\big( \Psi^\varphi( \muin)   - \langle \varphi,\mufin\rangle 
		\big),
	\end{equation*}
	where for every $\varphi \in C_{b,p}(\R^\xdim)$, the function $\Psi:\Pc_p(\R^\xdim)\longrightarrow \R$ is given by
	\begin{equation}
	\label{eq:def.Phi}
		\Psi^\varphi(\nu) \coloneqq \inf_{\gamma \in \Gamma(\nu)}J^\varphi(\gamma), \; \mbox{ where } J^\varphi(\gamma)\coloneqq \E^{ \P^\gamma}\bigg[\int_0^Tf(s, X^\gamma, u_s^\gamma, \mu_s^\gamma)\diff s+ Q_c\varphi(X_T^\gamma) \bigg].
	\end{equation}
Moreover, we have that \(\Psi^\varphi(\nu) =\inf_{\gamma \in \Gamma_{\rm det}(\nu)}J^\varphi(\gamma)\), and,
\begin{align*}
(V_c^{\rm det}(\muin,\cdot)^\ast)^\ast=V_c(\muin,\cdot).
\end{align*}	
%\review{We have not defined the set $\Pc$ at this moment.
%\begin{equation}
%	\label{eq:def.Phi}
%		\Psi^\varphi(\nu) \coloneqq \inf_{\P \in \Pc(\nu)}\E^{ \P}\bigg[\int_0^Tf(s, X, u_s, \mu_s)\diff s+ Q_c\varphi(X_T) \bigg].
%	\end{equation}
%Moreover, we have that
%\begin{align}
%\Psi^\varphi(\nu)
%=\inf_{\mathbb P\in \Pc_{\rm det}(\nu)}\E^{\P}\bigg[\int_0^T
%f(s,X,u_s,\mu_s)\diff s+Q_c\varphi(X_T) \bigg].
%\end{align}}
	%\camilo{do you mean just $\muin$-measurable? our proof holds even if both sides are $+\infty$}
%{\color{cyan}
Thus, as a function of the terminal law, $V_c$ is the lower-semicontinuous convex envelop of $V_c^{\mathrm{det}}$.
%}
\end{theorem}

%Theorem \ref{thm:existence} ensures the existence of an optimal weak control rule under Assumptions \ref{assump.costs} and \ref{assump.bf}, relying on compactness and stability of the martingale problem formulation.
\Cref{thm:duality.gen} provides a convex dual representation of the value function, extending weak optimal transport duality to the present controlled setting. 
	At the geometric level, this result identifies the value function as a lower semicontinuous convex functional of the endpoint law and represents it via its supporting hyperplanes. 
	While the terminal constraint $\mathbb P^\gamma\circ (X_T^\gamma)^{-1}\stackrel{c}{=}\nu_T$ is nonlinear at the level of admissible controls, it is convex when viewed through the weak transport functional $\Wc_c$. 
	The dual variable may therefore be interpreted as separating hyperplanes in the space of terminal measures: a test function $\varphi$ separates the prescribed target law $\mufin$ from the terminal laws generated by controlled McKean-Vlasov dynamics, while the weak $c$-transform $Q_c$ incorporates the geometry of the weak transport constraint into the intercept $\Psi^\varphi(\muin)$. 
	In this way, the weak marginal condition is converted into a terminal penalty in a stochastic control problem.
%It reveals that the value function is the convex envelope induced by stochastic control with terminal penalty given by the $c$-transform.
	Unlike previous results on stochastic optimal transport, see e.g., \cite{mikami2006duality,tan2013optimal}, the stochastic optimal transport problem in the dual formulation is of McKean-Vlasov type. 
	In fact, both the drift and the running cost depend on the law of the state process.
	Such problems naturally appear in the study of cooperative stochastic differential games, see, e.g., \citet*{carmona2018probabilisticI,carmona2018probabilisticII}, \citet*{djete2019mckean,djete2020mckean2}, \citet*{carmona2026conditional}, \citet*{lacker2017limit}, \citet*{soner2024viscosity} for a small list of references on the topic.
\medskip

We will now leverage \Cref{thm:duality.gen} to obtain a dual characterization of primal optimizers of McKean--Vlasov transportation problems.
This can be seen as a verification theorem, and it is essentially the mean field analog of the standard Schr\"odinger equation \cite{leonard2014survey}.	
In complete analogy with the classic literature, we introduce the following definition.

\begin{definition}\label{def.potential}
A Borel measurable function $\varphi:\R^\xdim\longrightarrow \R$  is a mean field potential whenever $V_c(\muin,\mufin)=\Psi^\varphi(\muin)-\langle \varphi,\mufin\rangle$.
\end{definition}

In fact, we will exploit the verification argument for the McKean-Vlasov control problem \eqref{eq:def.Phi} with state-dependent data and assume enough regularity to derive the characterization.
To this end, we introduce the Hamiltonian, for $(t,x,\mu,p, Q)\in[0,T]\times\R^\xdim \times\Pc_2(\R^\xdim )\times\R^\xdim\times \mathbb{S}^\xdim$,
\[
	H(t,x,\mu,p, Q)\coloneqq \inf_{{\bf u} \in \Uc } h(t,x,\mu,p, Q, {\bf u} ) , \mbox{ where, }  h(t,x,\mu,p, Q, {\bf u} )\coloneqq {\rm a}\cdot p+f(t,x,{\bf u} ,\mu)+ \tfrac12 {\rm Tr} [ Q\,  {\bf s} {\bf s}^{\!\top} ],
\]
%finite and upper semicontinuous, the infimum being attained by
%coercivity of $f$ in $a$ \eqref{eq:growth}. 
and consider the HJB equation with datum $\varphi:\R^\xdim\longrightarrow \R$ given by
\begin{equation}\label{eq:HJB}
\begin{cases}
	&\partial_tv(t,\mu)
	+
	\displaystyle \int_{\R^\xdim }\!
	\Big( 
	b(t,x,\mu)\cdot\nabla \delta v(t,\mu)(x)
	+
	H \big(t,x,\mu,\nabla \delta v (t,\mu)(x), \nabla^2 \delta v (t,\mu)(x)\big)\Big) 
	\mu(\d x)
	=
	0, \mbox{ on } [0,T)\times\Pc_2(\R^\xdim) \\
	&v(T,\mu)=\langle Q_c \varphi, \mu\rangle .
\end{cases}
\end{equation}

A function $v:[0,T]\times\Pc(\R^\xdim)\longrightarrow \R$ is said to be a solution of \eqref{eq:HJB}
% with datum $\varphi$, 
whenever $v \in\mathcal C^{1,2}_{\rm flat}$ satisfies \eqref{eq:HJB}. 
%with $\varphi$ a terminal condition.
Whenever we want to stress the dependence in the terminal condition $\varphi$ of a solution $v$ to \eqref{eq:HJB}, we denote it $v^\varphi$.

\begin{theorem}\label{theorem.verification}
Let {\rm \Cref{assump.bf}} hold.
Let $v$ be a solution to \eqref{eq:HJB} with datum $\varphi$. Moreover, let $\gamma^\star\in \Gamma(\muin)$.\footnote{We denote the elements of $\gamma^\star$ with a $\star$ superscript.}
The following are equivalent
\begin{enumerate}[label=$(\roman*)$, ref=.$(\roman*)$,wide,  labelindent=0pt]
\item 
\(
H (t,X^{\star}_t,\mu^{\star}_t, \nabla \delta v (t,\mu^{\star}_t )(X^{\star}_t),  \nabla^2 \delta v (t,\mu^{\star}_t )(X^{\star}_t))= h (t,X^{\star}_t,\mu^{\star}_t, \nabla \delta v  (t,\mu^{\star}_t )(X^{\star}_t),  \nabla_x^2 \delta v (t,\mu^{\star}_t )(X^{\star}_t), u^{\star}_t), 
\; \d t\otimes\d\P^{\star}\text{-a.e.}
\)
\item $\gamma^\star$ is optimal for $\Psi^\varphi(\muin)$ and $v(0,\muin)=\Psi^\varphi(\muin)$.
\end{enumerate}
	Moreover, if $c \equiv \int_{\R^\xdim} \1{x\neq y}\rho (\d y)$, then the following statements hold:
	\begin{enumerate}[label=$(\alph*)$, ref=.$(\alph*)$,wide,  labelindent=0pt]
	\item Let $\varphi^\star$ be a mean field potential and $v^{\varphi^\star\!\!}$~be a solution to \eqref{eq:HJB}. 
	Let $\gamma^\star\in \Gamma_{\rm det} (\muin,\mufin)$ be such that $(i)$ holds, then $\gamma^\star$ is optimal for $V_c(\muin,\mufin)$. 
\item Let $\gamma^\star\in \Gamma_{\rm det} (\muin,\mufin)$ be optimal for $V_c(\muin,\mufin)$. 
	Let $\varphi^\star$ be such that $(i)$ holds with $v^{\varphi^\star\!\!}$ solution to \eqref{eq:HJB}, then $\varphi^\star$ is a mean field potential.
\end{enumerate}
\end{theorem}

The mean field optimal transport problem discussed so far is intrinsically infinite-dimensional, since $V_c$ is the value function of a controlled McKean-Vlasov dynamics. 
	This problem is the natural infinite-particle limit of an associated finite-particle problem.
	This is referred to in the literature as the so-called convergence problem.
	The second aim of the present work is to make this link rigorous in the context of the mean field entropic Schr\"odinger problem, introduced in \cite{backhoff2020mean}, which corresponds to the choice $f\equiv \frac12 \|\sigma^\dagger a\|^2$, $c\equiv \int_{\R^\xdim} \1{x\neq y}\rho (\d y)$ and $\Uc\equiv  \R^{\xdim}\times \{\sigma\}$ in \eqref{eq:problem}, where $\sigma^\dagger$ denotes the Pseudoinverse of $\sigma$.\footnote{The choice of $c$ corresponds to the total variation distance, see, e.g., \cite[Remark 3.2]{hernandez2025marginal}, and $\Uc$ is chosen so that the volatility is uncontrolled.}
%	{
%	\color{cyan}
	However, the distinction between deterministic and randomized law flows is important here.
	In the formulation of \cite{backhoff2020mean}, the interaction is evaluated at a deterministic flow. 
	Our weak formulation enlarges this class by allowing an exogenous randomization and evaluating the coefficients at the corresponding law flow.
	This enlargement is immaterial for the unconstrained dual control problems, which can be reduced to deterministic-flow rules, but it is essential for the constrained problem: the deterministic-flow value need not be convex in the target law, but the randomized formulation is.
%	}
	For readability, we introduce
\begin{equation}
\label{eq:e.problem}
	V_e(\muin, \mufin) \coloneqq \inf_{\gamma \in \Gamma(\muin,\mufin)} J_e(\gamma), \; \mbox{where } J_e(\gamma)\coloneqq\E^{\P^\gamma}\bigg[ \int_0^T \frac12 \|\sigma^\dagger  \alpha^\gamma_t\|^2 \diff t\bigg],
\end{equation}
and $\Gamma(\muin,\mufin)\coloneqq\big\{ \gamma \in \Gamma(\muin): \P^\gamma \circ (X_T^\gamma)^{-1}=\mufin \big\}$.
We remind the reader that {\color{black}for deterministic-flow rules,} thanks to \cite[Lemma 1.1]{backhoff2020mean}, $V_e(\muin, \mufin)$ corresponds to the nonlinear entropy minimization problem example discussed in \Cref{sec:intro}. \medskip

Let us now introduce the finite particle optimal transport problem associated with $V_e(\muin, \mufin)$.

\paragraph*{The $N$-particle Schr\"odinger problem.}
	To motivate the following formulation, in complete analogy to \Cref{sec.formulation}, for any $N\geq 1$, we consider a continuous semimartingale $\X=(X^{1,N}, \dots, X^{N,N})$, defined on some filtered probability space carrying an $N$-tuple $(W^1, \dots,W^N)$ of $\bmdim$--dimensional Brownian motions independent from $\X_0$, for which there is an $N$-tuple $\balpha=(\alpha^1,\dots,\alpha^N)$ of $\R^\xdim$--valued measurable processes such that
	%\todo[inline]{properly write $\sigma$ in the dynamics. Make sure it is consistent in the rest of the paper.}
	\begin{equation*}
	``X^{i,N}_t = X_0^{i,N} + \int_0^t\big( \alpha^i_s+ b( X^{i,N}_s, m^N_{\X_s})\big) \diff s + \sigma \diff W^{i}_t,\; i=1,\dots,N,\;  \Lc(\X_0)=\muin^{\otimes N}, \text{ and, } m^N_\X= \frac1N\sum_{i=1}^N\delta_{X^{i,N}}."
\end{equation*}
In particular, here $b$ is a function of the state at a fixed time.
Unlike the general setting of \Cref{thm:duality.gen}, it is natural to impose additional conditions on $b$, see \Cref{assump.Df}\ref{assump.Df.2}, so that the limit theory for standard McKean-Vlasov control problems holds.
	In fact, under \Cref{assump.Df} the above dynamics will have a unique strong solution and the previous formulation will be independent of the background probability space.
	Consequently, we may, and will, formulate the $N$-particle problem on a fixed canonical space.\medskip
	
Let $N\geq 1$ and consider, respectively, the canonical space and canonical process
\[
\Omega^N\coloneqq (\R^\xdim \times \Cc^\xdim)^N, \; (\X_0, \X)\coloneqq ((X^{1}_0, \dots, X^{N}_0),(X^{1,N}, \dots, X^{N,N})),
\]
with canonical filtration  $\F^N\coloneqq (\Fc_t^N)_{t\in [0,T]}$, $\Fc_t^N\coloneqq \sigma(X^i_0,X^{i,N}_s, s\in [0,t], i = 1,\dots,N ), t \in [0,T]$.
	%{\color{red}I was thinking we should rather have ... the canonical space $\Omega^N\coloneqq \Cc^{\xdim N}$ with canonical process $\X\coloneqq (X^{1,N}, \dots, X^{N,N})$. We denote by $\F^N$ the canonical filtration, i.e., $\Fc_t^N\coloneqq \sigma(X^i_s, s\in [0,t], i = 1,\dots,N )$ for all $t \in [0,T]$.}
	\medskip

\begin{definition}\label{def.weak.N.sol}
Let $\muin\in\Pc(\R^\xdim)$, $N\geq 1$ be fixed.
	$\Gamma^N(\muin^{\otimes N})$ denotes the collection of $\P^{N\!,\balpha}\in {\rm Prob}(\Omega^N)$ for which there is a $\R^{N\xdim}$--valued, $\F^N$--predictable $\boldsymbol{\alpha} = (\alpha^1,\dots,\alpha^N)$ such that
\begin{equation*}
	\sum_{i=1}^N\E^{\P^{N\!,\balpha}} \bigg[\int_0^T\|\sigma^\dagger \alpha^i_t\|^2\diff t\bigg]<+\infty, %\quad \mbox{for some } {\rm u_0} \in \Uc
\end{equation*}
%{\color{blue}What is $\sigma^\dagger$? why not $\sigma^\top$?}
%\begin{equation*}
%	\dfrac{\d\P^{N,\boldsymbol{\alpha}}}{\d\P^N} \coloneqq \Ec\bigg(\sum_{i=1}^N\int_0^\cdot\sigma^{-1}_s \alpha^i_s \cdot \d W^i_s \bigg)_T,
%\end{equation*}
and a $(\P^{N\!,\boldsymbol{\alpha}},\F^N)$--Brownian motion, $(W^{1,\P^{N\!,\boldsymbol{\alpha}}},\dots,W^{N,\P^{N\!,\boldsymbol{\alpha}}})$ such that
\begin{equation*}
%\label{eq:sde.girsanov}
	X^{i,N}_t = 
	X_0^{i} 
	+ 
	\int_0^t\big( \alpha^i_s+ b( X^{i,N}_s, m^N_{\X_s})\big) \diff s 
	+ 
	%\int_0^t
	\sigma %\diff 
	W^{i,\P^{N\!,\boldsymbol{\alpha}}}_t,\; t\in [0,T], \;  \P^{N,\balpha}\as,\; i=1,\dots, N,\;  \P^{N,\balpha}\circ (\X_0)^{-1}=\muin^{\otimes N}.
\end{equation*}
\end{definition}
Given $\mufin\in\Pc(\R^\xdim)$, the $N$-particle stochastic optimal transport problem under consideration is thus
\begin{equation}
\label{eq:N-part.schroedinger}
	V^N_e(\muin, \mufin)
	\coloneqq 
	\inf\Big\{
	J^N_e\big(\P^{N\!,\balpha}\big):\;  
	\P^{N\!,\balpha}\in \Gamma^N(\muin^{\otimes N}) ,\;
	 \P^{N\!,\balpha} \circ (X^{i,N}_T)^{-1} = \mufin, \; i= 1,\dots, N 
	 \Big\},
\end{equation}
with %{\color{blue}(Not sure about the equality below. Are you missing $\Hc(\muin^{\otimes N}|\mu^N)$)?}
\begin{equation}\label{eq.rewardentropy}
	J^N_e\big(\P^{N\!,\balpha}\big) \coloneqq \E^{\P^{N\!,\balpha}}\bigg[\sum_{i=1}^N\int_0^T\frac12 \|\sigma^\dagger   \alpha^i_t\|^2 \diff t\bigg ].
	%=\Hc\big( \P^{N,\balpha} \big|\P^N\big).
\end{equation}

\medskip

%{\color{red}we can discuss the $N$-particle problem here when we are ready (move the description below to here).}

Notice that in the above problem, the terminal constraint is imposed at the level of \emph{individual particles}.
This means that we are interested in the most likely path of a cloud of \emph{interacting particles} starting from an i.i.d.~configuration and each having terminal marginal $\mufin$.
This \emph{particle-by-particle} formulation serves a dual purpose: on the one hand, it offers a probabilistic interpretation of the mean field problem in terms of large systems of controlled particles; on the other hand, it enables a rigorous passage to the mean field limit.\medskip 

To present our result, it remains to introduce $\P^N$ a probability measure on $(\Omega^N, \Bc(\Omega^N))$ under which 
there is an $(\F^N,\P^N)$--Brownian motion $(W^{1,\P^N},\dots, W^{N,\P^N})$ independent of $\X_0$, and 
%{\color{red}such that} the process $\X\coloneqq (X^1, \dots,X^N)$ satisfies
\begin{equation}\label{eq.Np.sde.ref}
	X^{i,N}_t 
	= 
	X_0^{i} 
	+ 
	\int_0^t 
	b(X_s^{i,N}, m^N_{\X_s})\diff  s
	+ 
	%\int_0^t
	\sigma
	%\diff 
	W^{i,\P^N}_t,\, t\in [0,T], \, \P^N\as, \, i=1,\dots, N,\text{ and, } \P^N\circ\X_0^{-1}=\mu^N\in \Pc(\R^{\xdim N}). %\; m^N_\X\coloneqq \frac1N\sum_{i=1}^N\delta_{X^{i,N}},
\end{equation}	
%, which we will take to be the stationary distribution of the above SDE. 

Let us now introduce the time reversal mapping $\overleftarrow{\Tc}:\Omega^N\longrightarrow \Omega^N$ given by $\overleftarrow{\Tc}(\omega)\coloneqq \overleftarrow{\omega}$ where $\overleftarrow{\omega}_t\coloneqq \omega_{T-t}$ and, for any $\Q\in {\rm Prob}(\Omega^N)$, the time-reversed measure $\overleftarrow{\Q} \coloneqq \Q\circ \overleftarrow{\Tc}$.
	Recall that a measure $\Q$ is said to be reversible if  $\overleftarrow{\Q} = \Q$.

\begin{assumption}\label{assump.Df}
	{\rm \Cref{assump.Df}} is satisfied if $m=1$ %\todo{Ask Ludo} 
	and it holds
\begin{enumerate}[label=$(\roman*)$, ref=.$(\roman*)$,wide,  labelindent=0pt]
	\item \label{assump.Df.1} {$\sigma\sigma^\top \ge \Lambda I_\xdim$ for some $\Lambda>0$}.
	\item \label{assump.Df.2} {The function $b$ is time-independent, continuously differentiable and $\ell_b$--Lipschitz continuous with linear growth, with $\Pc_2(\R^\xdim)$ equipped with the $2$-Wasserstein metric.
		 % $(b(x^1,m_{\bx}),\dots, b(x^N, m_{\bx})) = -\nabla H_N(\bx)$.
		Moreover, 
		%the function $H_N$ satisfies the convexity condition $\langle \xi, \nabla^2H_N\xi\rangle\ge \rho |\xi|^2$ for some $\rho\ge0$ and every $\xi \in \R^{\xdim N}$.
		the $\R^{\xdim N}$-valued function $\Bc(\bx) \coloneqq (b(x^1,m^N_{\bx}),\dots, b(x^N, m^N_{\bx}))^\t$ satisfies 
		%$-\langle \xi, \nabla \Bc(\bx) \xi\rangle\ge  \rho |\xi|^2$ for some $\rho\ge0$ and every $\xi \in \R^{N}$ and $\bx \in \R^N$.\medskip
		the dissipativity condition: there is $\rho\geq 0$ such that
			\[
			-\left\langle \Sigma^\top \xi,\,\Sigma^\top \nabla\Bc(\bx)\xi \right\rangle
			\ge
			\rho\|\Sigma^\top\xi\|^2, \mbox{ for all }  \bx,\xi\in\R^{\xdim N},
			\]
where $\Sigma\coloneqq I_N\otimes \sigma \in \R^{\xdim N\times \bmdim N}$ denotes the Kronecker product of $I_N\in\R^{N\times N}$ with $\sigma\in \R^{\xdim\times\bmdim}$.}

%\todo[inline]{Equivalently, in matrix form,
%\(
%\frac12\Big(
%A\nabla\Bc(\bx)
%+
%\nabla\Bc(\bx)^\top A
%\Big)
%\preceq
%-\rho A.
%\)
%Your old condition was
%\[
%-\langle \xi,\nabla\Bc(\bx)\xi\rangle
%\ge \rho|\xi|^2.
%\]
%This is equivalent to
%\[
%\frac12\big(\nabla\Bc+\nabla\Bc^\top\big)
%\preceq -\rho I.
%\]
%
%By itself, this does not imply the weighted condition above for a general non-scalar $A=\Sigma\Sigma^\top$. The reason is that
%\(
%p^\top A\nabla\Bc\,p
%\)
%cannot in general be controlled by
%\(
%p^\top\nabla\Bc\,p
%\).
%If
%$\Sigma\Sigma^\top=\lambda I_{mN}$,
%for example $ \sigma\sigma^\top=\lambda I_m$, then
%\(
%-p^\top A Jp
%=
%-\lambda p^\top Jp
%\ge
%\rho\lambda |p|^2
%=
%\rho\|\Sigma^\top p\|^2,\)
%}
	\item \label{assump.Df.3} %The particle system \eqref{eq.Np.sde.ref} admits a unique invariant measure denoted $\mu^N$ and 
		For any $N\geq 1$, the measure $\P^N$ is reversible. 
	\item \label{assump.Df.4} There is $c_0>0$ such that $\int_{\R^\xdim}e^{c_0\|x\|^2}\nu_0(\d x)<+\infty$.
	\item \label{assump.Df.5} $\mufin$ satisfies Poincar\'e inequality, is equivalent to Lebesgue measure $\lambda$ and $w_T \coloneqq \d\mufin/\d\lambda$ satisfies the lower-envelope condition: there is $C_\nu>0$ and $\kappa>0$ such that 
\begin{equation*}
	\inf_{y\in B_x}w_T(y) \ge C_{\nu} e^{-\kappa(1 + \|x\|^2)}, \mbox{ for all $x\in \R^\xdim$,}
\end{equation*}
where $B_x$ is the ball centered at $x$ of radius $1$, and $\kappa <\min (c_0/c_1, (2T\lambda_{\rm max}(\sigma\sigma^\top) c_1)^{-1})$ with $c_1\coloneqq 2\e^{8\ell_b T^2}\big(1 + 8\ell_b T^2 e^{16\ell_b T^2} \big)$.
%{\color{cyan}Camilo, can you please write something to "defend" the condition on the bounds?}	

%\camilo{and $\mufin\sim\lambda$ with $\d\mufin=w\, \d \lambda$ and $w, w^{-1} \in \L^\infty_{\rm loc}(\lambda)$}.
	\end{enumerate}
\end{assumption}

\begin{remark}
Let us comment on the previous conditions on the data.
\begin{enumerate}[label=$(\roman*)$, ref=.$(\roman*)$,wide,  labelindent=0pt]	
\item The uniform ellipticity condition guarantees that the infinitesimal generator $L_N$ of the diffusion $\X$ defined by \eqref{eq.Np.sde.ref} is a strictly elliptic operator.
	This, together with the assumption on the dimension of the state process, will provide the regularizing effects that allow us to derive a key Sobolev compactness of a sequence of Schr\"odinger potentials in the proof of the convergence result.
	The assumption of constant $\sigma$ simplifies the analysis involving the associated carr\'e du champ operator, see {\rm \Cref{sec.upperbound}}.
	We expect that non-constant, uniformly elliptic $\sigma$ could also be treated, at the cost of substantial additional technicalities.
	We stress that the assumption on $m$ is needed only at one step in the argument, namely, {\rm \Cref{prop.conv.phi^N}}.
	This is the reason we keep the notation for the general case in the proof.
	See {\rm \Cref{remark.assump.dim.den}} for further comments on this.
	
\item In {\rm \Cref{assump.Df}\ref{assump.Df.2}--\ref{assump.Df}\ref{assump.Df.3}} we have in mind models of the form
\begin{equation*}
	b(x,\mu) = -\nabla F(x) - \int_{\R^\xdim} \nabla_x G(x,y)\, \mu (\d y) ,
\end{equation*}
for a $C^2$ and convex confining potential $F$, and a $C^2$, bounded from below, symmetric and jointly convex interaction potential $G$ {and $ \sigma\sigma^\top=\Lambda I_\xdim$}.
	For our purpose, in {\rm\Cref{assump.Df}\ref{assump.Df.2}} it is enough to assume $\rho\ge0$.
	For the above form of $b$, if the (joint) convexity coefficient $\rho$ is strictly positive and $\sigma$ is time-independent, then {\rm \Cref{assump.Df}\ref{assump.Df.2}--\ref{assump.Df}\ref{assump.Df.3}} are already satisfied, see {\rm \cite{cattiaux2023time}}.
\item The reversibility assumption on $\P^N$ will be key to derive the regularity of a sequence of potentials of $V^N$.
	It follows from {\rm\Cref{assump.Df}\ref{assump.Df.2}} that $\P^N$ is the unique weak solution to \eqref{eq.Np.sde.ref}.

% and the invariant measure is given in Gibbs form (when $\sigma$ is constant).
% by $\mu^N(d\bx) = e^{-\frac{1}{\sigma}H_N(\bx)}/Z^N$.
\item Assuming that $\mufin$ satisfies Poincar\'e inequality can be seen as a regularity constraint on the target distribution. 
This condition is well known to be weaker than the log-concavity of $\mufin$.

\item The conditions in {\rm \Cref{assump.Df}\ref{assump.Df.4}--\ref{assump.Df}\ref{assump.Df.5}} are technical and required in our analysis.
	Given the controlled dynamics in \eqref{eq:N-part.schroedinger}, we believe the square exponential form of the lower envelope of the density $w_T$ is compatible with Gaussian-type tails.
	The bound on $\kappa$ is, ultimately, a restriction relating $T$ and $\ell_b$.
	This condition arises from the need to derive uniform exponential moment estimates in our proof of the upper bound.
	If one is willing to assume a stronger form of the lower envelope, this condition can be weakened or removed altogether.
	%These model parameters represent the time horizon and the interaction coefficient, respectively.
	Indeed, see {\rm \Cref{remark.assump.dim.den}} for other settings under which these assumptions can be weakened or dropped.
\end{enumerate}
\end{remark}
% \begin{assumption}\label{assump.cf}
% 	The function $b$ satisfies Assumption {\rm \ref{assump.bf}}, $f(t,\bx,u,\mu) = \frac12\|a\|^2$ for all $u = (a,\sigma)$.
% 	In particular, $(\sigma_t)$ is uncontrolled 
% 	{\color{red}if we stick with $\Uc$ Polish here, then we should replace $|\alpha|^2$ by $d^2_{\Uc}(\alpha, u_0)$ everywhere in second part of the proof. I don't like it.}
% 	\begin{enumerate}[label=$(\roman*)$, ref=.$(\roman*)$,wide,  labelindent=0pt]
% 	\item {\color{red}For the quadratic function $f(a) = |a|^2$ do we have $f(a) \le C - C_L|a|^p$ with $C_L>0$ and $p>2$?}
% %	\item Here, add the conditions of Theorem \cite[Theorem 3.6 and Proposition 3.7]{djete2020mckean2}.
% 	\end{enumerate}
% \end{assumption}

The following result establishes that the mean field optimal transport problem arises as the asymptotic limit of a system of controlled interacting diffusions.
Intuitively, we understand 
%$V_c$ as the optimal ``cost per particle'' in the infinite system problem, whereas 
$V_N$ as the optimal cost to steer $N$ weakly interacting particles to a distribution $\mufin$.
Theorem \ref{thm:parti.conv.N} says that the optimal ``cost per particle'' $N^{-1}V_e^N$ of the finite particle problem converges to $V_c$, the effective ``cost per particle'' in the infinite limit. 
%making $1/N$ the right scaling. 

\begin{theorem}
\label{thm:parti.conv.N}
	Let Assumptions {\rm \ref{assump.bf}} and {\rm \ref{assump.Df}} be satisfied, and further assume that $\muin\in \Pc_p(\R)$ for some $p>2$.
	Let $V_e$ and $V^N_e$ be respectively given by \eqref{eq:e.problem} and \eqref{eq:N-part.schroedinger}.
	If we have
	\begin{align}\label{assump.parti.conv.N}
	\sup_{N\ge1}N^{-1}V^N_e(\muin,\mufin) <\infty,\mbox{ and, } \sup_{N\ge1}N^{-1} \Hc(\muin^{\otimes N}|\mu^N) <\infty,
	\end{align}
%	Let the weak cost $c$ be given by $c(x,\rho) = \int_{\R^\xdim}\1{x\neq y}\rho(dy)$.
	then it holds
	\begin{equation*}
		\lim_{N\to \infty}N^{-1}V^N_e(\muin,\mufin) = V_{e}(\muin,\mufin).
	\end{equation*}
\end{theorem}
Establishing this convergence is a difficult task due to the lack of independence at the terminal time and the nonlinear interaction through the controls. 
	The proof of this convergence result makes use of the duality theorem stated earlier, along with gradient estimates of Schr\"odinger potentials.
	The result can also be interpreted as a variational convergence statement linking finite-dimensional stochastic control problems to optimal transport and highlighting the robustness of the entropic formulation under mean field scaling.
	In fact, the entropic case is particularly suited to derive the convergence result because it allows to characterize the (near) optimal $N$-particle bridge $\Q^N$ as
 \begin{align*}\label{eq.QN.charac.main}
		\dfrac{\d   \Q^N}{\d  \P^N} 
		=
		\dfrac{\d    \muin^{\otimes N}}{\d   \mu^N}( \X_0) 
		\exp\big(
		- \varphi^{N}(\X_T) 
		- \psi^N(\X_0)
		\big)
\end{align*}
for some (near) optimal potentials $\varphi^N$ and $\psi^N$.
We do not know whether the result extends to the case of stochastic optimal transport problems with generic cost $f$.

\medskip

 An important consequence of our proof is the following existence result of an optimal potential for the mean field Schr\"odinger problem.
 This potential admits a weak gradient.
\begin{corollary}
\label{cor:dual.exists.char}
	Under the conditions of {\rm \Cref{thm:parti.conv.N}}, there exists an optimal mean field Schr\"odinger potential.
	That is, there is a function $\varphi\in \L^2(\R,\mufin)$ admitting a weak gradient $($in the sense of distributions$)$ in \(\L^2(\R,\mufin)\)
	%{\color{blue}I think we can say that this function has a weak gradient}
	such that 
	\[
	V_e(\muin,\mufin) =   \Psi^\varphi(\muin) - \langle \varphi,\mufin\rangle.
	\]
\end{corollary}

\begin{remark}\label{remark.assump.dim.den}
\begin{enumerate}[label=$(\roman*)$, ref=.$(\roman*)$,wide,  labelindent=0pt]
\item It is natural to ask whether the above potential is unique, at least up to an additive constant, as is the case in the classical, i.e., non mean field, Schrödinger problems.
Based on the discussion in {\rm \Cref{rmk.exist.unique}}, we leave this for future work.
%	  This is unclear {\color{cyan}is this unique ?}
\item The growth requirement on the density of $w_T$ of $\mufin$ and the condition $\xdim=1$ are technical conditions that may be avoided using different proof techniques.
They are needed here only to be able to use Sobolev embedding with respect to Lebesgue measure (embed $H^1$ into $\L^\infty$), which allows to derive uniform in $N$ growth properties on Schr\"odinger potentials of the $N$ particle problem.
These conditions can be avoided for instance if we can, a priori, show that the $N$ particle problem admits optimal potentials $\varphi^{i,N}$ that have polynomial growth: $|\varphi^{i,N}(x)| \le C(1 + |x|^p)$ for some $p\ge0$.
Moreover, these conditions on the dimension as well as the tail estimates and lower-envelope assumptions on the target density, can be dropped if the problem is set on a torus as for instance in {\rm \citet*{berman2020sinkhorn,baradat2021regularized,lavenant2019new,yao2025convergence,soner2024viscosity}}. In this case, the state space would be compact and the growth properties needed on the Schr\"odinger potential would be easier to derive.
\end{enumerate}
\end{remark}

Let us conclude this section by specializing \eqref{eq:HJB} to the mean field entropic Schrödinger problem considered by \cite{backhoff2020mean}.

%\begin{assumption}\label{ass:R}
%\begin{enumerate}[label=$(\roman*)$, ref=.$(\roman*)$,wide,  labelindent=0pt]
%Let $\varphi$ be an optimal potential constructed in \Cref{cor:dual.exists.char}.
%We assume that:  
%\item  AA \label{ass:R.2} 
%%\item 
%% solve \eqref{eq:HJBBCGL},
%%In addition to $v\in\mathcal C^{1,2}_{\rm flat}$: the function
%%and 
%%$u(t,\mu)(y)$ is 
%\end{enumerate}
%\end{assumption}

%\medskip\noindent\textbf{Standing situation.} Let
%$v\in\mathcal C^{1,2}_{\rm flat}$ solve \eqref{eq:HJBBCGL}, 
%%let $(t_0,\mu_0)\in[0,T)\times\Pc_2(\R^d)$, and 
%let $\gamma^*\in\Gamma_{\rm det}(\mu_0)$ be an optimal feedback rule as in
%\Cref{thm:verification}(ii) (by \Cref{lem:reduction} it may be taken
%with deterministic flow). Write $\mu^*_t:=\mu^{\gamma^*}(t)$ for
%$t\in[0,T]$ and
%\[
%U(t,x):=u(t,\mu^*_t)(x),\qquad
%\mathcal A_tg(y):=\tfrac12\Delta g(y)-\big(\nabla W*\mu^*_t(y)+\nabla U(t,y)\big)\cdot\nabla g(y).
%\]
%By \Cref{thm:verification}(iii), $\alpha^{\gamma^*}_t=-\nabla U(t,X^{\gamma^*}_t)$
%a.e., so $\mu^*$ solves the Fokker--Planck equation with generator
%$\mathcal A_t$: for every $g\in C^2(\R^d)$ with $|g|+|\nabla g|+|\nabla^2g|$
%of growth at most $|x|$,
%\begin{equation}\label{eq:FP}
%\frac{\d}{\d t}\int g\,\d\mu^*_t=\int\mathcal A_tg\,\d\mu^*_t,
%\qquad t\in[t_0,T].
%\end{equation}
%

\begin{corollary}\label{cor.HJB.linearization}
	Let $V\in C^2(\R^\xdim)$ be symmetric with bounded Hessian, $b(t,x,\mu)\equiv -\nabla V\ast \mu(x)$, $\xdim=\bmdim$, and $\sigma=I$.
	Let $v\in \Cc^{1,2,{\rm full}}_{\rm flat}$ be a solution to \eqref{eq:HJB} with datum a mean field potential $\varphi^\star$, as constructed in {\rm \Cref{cor:dual.exists.char}}.
	Let $\gamma^\star\in \Gamma_{\rm det} (\muin,\mufin)$ be such that 
%with $U(t,x)\coloneqq u(t,\mu_t^{\star})(x)$, 
	$\alpha_t^{\star}=-\nabla \delta v (t,\mu_t^\star )( X^\star_t),  \; \d t\otimes \d \P^{\star}\ae$ 
	Then, there is %$\psi^\star:[0,T]\times\R^\xdim\longrightarrow \R$, 
	$\psi^\star\in C^{1,2}([0,T)\times\R^\xdim)\cap  C([0,T]\times\R^\xdim)$, satisfying
\begin{equation}\label{eq:BCGL17}
	\partial_t\psi^\star
	+
	\tfrac12\Delta\psi^\star
	+
	\tfrac12|\nabla\psi^\star |^2
	-
	 \int_{\R^\xdim}
	 \nabla V(x-\tilde x)\cdot\big(\nabla\psi^\star(t,x)-\nabla\psi^\star(t,\tilde x)\big)
	 \,\mu^\star_t(\d\tilde x)
	 =
	 0,\;
	\psi^\star(T,\cdot)=-\varphi^\star
%\end{cases}
\end{equation}
In particular, $\psi^\star (t,x)\coloneqq -  \delta v (t,\mu_t^\star )( x) +\int_t^T\varrho(s)\,\d s+\int\varphi^\star\,\d\mu^\star_T$, for a function $\varrho:[0,T)\longrightarrow \R$.\medskip
\end{corollary}

\begin{remark}
\begin{enumerate}[label=$(\roman*)$, ref=.$(\roman*)$,wide,  labelindent=0pt]
\item	It is not hard to show that $\mu^\star$ satisfies $\partial_t\mu^\star_t-\tfrac12\Delta\mu^\star_t+\nabla\cdot\big((-\nabla V\ast\mu^\star_t+\nabla\psi^\star)\mu^\star_t\big)=0$ in the weak sense.
	Consequently, the pair $(\psi^\star,\mu^\star)$ is a solution of {\rm \cite[Corollary 1.1]{backhoff2020mean}}, and moreover, $\psi^\star$ is a solution of the corresponding equation in {\rm \cite[Corollary 1.2]{backhoff2020mean}}, although with specified terminal condition $\psi^\star(T,\cdot)=-\varphi^\star$ where $\varphi^\star$ is a mean field potential for $V_e(\muin,\mufin)$.

\item Let us also note that, in contrast to {\rm \Cref{theorem.verification}}, {\rm \Cref{cor.HJB.linearization}} imposes additional regularity on $v$.
	To see why, recall that the terminal condition of the {\rm HJB} equation \eqref{eq:HJB} is $\langle \varphi,\mu\rangle$, echoing the fact that the ambient domain of this equation is $[0,T] \times \Pc_2(\R^\xdim)$.
	Consequently, arriving at \eqref{eq:BCGL17} requires additional regularity in the measure variable for $v$.
	This is precisely the idea behind {\rm \Cref{cor.HJB.linearization}}, which identifies \eqref{eq:HJB} with \eqref{eq:BCGL17} after \emph{linearizing}, i.e., taking an additional linear derivative, alongside the optimal flow $\mu^\star_\cdot$.
\end{enumerate}
\end{remark}

\section{Proofs of existence and convex duality results}
\label{sec:proofs_of_existence_and_convex_duality_results}
For the proofs it will be more convenient to set the problem on an enlarged canonical space which we now introduce.
\paragraph*{The Canonical Space.}

Let us reformulate the set of weak control rules in the canonical space.
	%\footnote{\camilo{As an isomorphism of metric spaces it follows that $\pi$ is continuous}}
	Letting $\iota$ denote a cemetery point and $\overline \Uc:=\Uc\cup\{\iota\}$, we extend $\pi^{-1}$ to $\pi^{-1}:\R\cup\{-\infty,\infty\}\longrightarrow \overline \Uc $, setting $\pi^{-1}(x):=\iota$ for all $x\notin \pi(\Uc)$.
	With this, as mentioned in \eqref{eq.repofcontrols}, we can represent a control $u=(\alpha,\sigma\sigma^\top)$ %{\color{blue}$u=(\alpha,\sigma\sigma^\top)$?} 
	via the real-valued continuous process $\int_0^\cdot \! \pi(u_t)\d t$ and take $\Cc$ as associated canonical space. 
	\medskip

Consider first the canonical space $\bar \Omega \coloneqq \Cc^\xdim\times\Cc\times\Cc^\bmdim$ with canonical process $(\bar X,\bar U,\bar W)$.
	We denote by $C([0,T],\Pc(\Cc^\xdim))$ the space of all continuous paths on $[0,T]$ with values in $\Pc( \Cc^\xdim)$ topologized by uniform convergence so that it is a Polish space, see \cite[Section 3.19]{aliprantis2006infinite}.
	We consider the canonical space
\[
\Omega\coloneqq\Cc^\xdim \times \Cc \times \Cc^\bmdim\times C([0,T],\Pc(\bar\Omega)), \text{ with canonical process } (X,U,W,  \bar \mu),  \text{ and generic elements }\omega\coloneqq ({\rm \bf x},  {\bf u},{ \bf w},{\rm \bf m}).
\]
%{
%	\color{red} should $\bar\Omega$ above rather be $\Cc^m$?
%}
Let $\Fc$ denote the Borel $\sigma$-field on $\Omega$ and $\F:=(\Fc_t)_{t\in [0,T]} $ be the canonical filtration.
	The process $\bar \mu\coloneqq (\bar  \mu_t)_{t\in [0,T]}$ carries information and induces an element of $\Pc(\Cc^\xdim)$ via $\bar \mu_t$ for any $t\in [0,T]$. 
	We thus define
	\[
	\bar\G\coloneqq (\overline{\Gc}_t)_{t\in [0,T]},\;
\overline{\Gc}_t\coloneqq \sigma\big(\langle\phi,\bar \mu_s\rangle :\, s\in[0,t],\, \phi\in C_b(\bar \Omega)\big),\;  t\in[0,T],
	\]
	and the process $\mu\coloneqq(\mu_t)_{t\in [0,T]}$ on $\Omega$ as follows: for every $\omega  =({\rm \bf x},  {\bf u},{ \bf w},{\rm \bf m}) \in \Omega$, we put
% via {\color{blue}how is the $\omega$ dependence on the rhs below?} \camilo{on the coordinate ${\bf m}$ of $\omega$}
\begin{align}\label{eq.lawcanonical}
\mu_t(\omega) \coloneqq \bar \mu_t \circ \bar X_{t \wedge \cdot}^{-1} (\omega) =\E^{{\bf m}_t}\big [\delta_{\bar X_{t\wedge \cdot}}\big],\, (t,\omega)\in [0,T]\times\Omega .
\end{align}
%{\color{blue}what is $\bar X$?}	\camilo{One can think of $\bar X$ as a copy of $X$ and allows us to define $\mu$ independently from $\P$ (important for l.s.c.). This is in contrasts to defining the law of $X$ (without introducing $\bar\mu$ and $\bar X$) as $\mu^\P:=\P\circ X_t^{-1}$}

Moreover, we also define the $\overline \Uc$-valued process $u:=(u_t)_{t\in [0,T]}$ on $\Omega$ given by
\begin{align}\label{eq.Ucanonical}
u_t\coloneqq \pi^{-1}\Big( \limsup_{n\to\infty} n(U_t-U_{(t-1/n)^+})\Big), \, t\in [0,T],
\end{align}
inducing a pair of processes $u=(\alpha,\sigma\sigma^\top)$ %{\color{blue}$(\alpha,\sigma\sigma^\top)$?} 
on $\Omega$. \medskip

We now reformulate the weak control rules in the canonical space. To do so, for $(t,{\bf x},{\bf w},\bmu,{\rm u})\in [0,T]\times \Cc^\xdim \times \Cc^\bmdim\times \Pc(\Cc^\xdim)\times \Uc$ with ${\bf u}=({\bf a},{\bf s}{\bf s}^\top)$, %{\color{blue}${\bf u}=({\bf a},{\bf s}{\bf s}^\top)$?} 
we introduce 
\[
\bar b(t,{ \bf x},{ \bf w}, \bmu,{\bf u})\coloneqq \begin{pmatrix}
{\bf a}+ b(t,{\bf x} , \bmu)\\
 0_d
 \end{pmatrix},\;
 { \overline \Sigma}(t,{\bf x},{\bf w} ,{\bf u})\coloneqq \begin{pmatrix}
{\bf s} & 0_{d\times d}\\
I_{d\times d} & 0_{d\times d} 
 \end{pmatrix} 
 \begin{pmatrix}
{\bf s} & 0_{d\times d}\\
I_{d\times d} & 0_{d\times d} 
 \end{pmatrix}^{\!\t},
\]
and, for $\phi \in C_b^2(\R^{m+d})$, the generator $\Lc$ given by
\[
(\Lc\phi)(t,{\bf x} ,{\bf w},\bmu,{\bf u}) \coloneqq \bar b(t,{ \bf x},{ \bf w}, \bmu,{\bf u})\cdot \nabla \phi ({\bf x}_t,{\bf w}_t)+\frac12 {\rm Tr}\big [\overline \Sigma (t,{\bf x},{\bf w} ,{\bf u})\nabla^2 \phi ({\bf x}_t,{\bf w}_t)\big].
\]

\begin{definition}\label{def.controlrule}
Let $p\geq 1$ and $\nu\in \Pc(\R^\xdim)$. We denote by $\Pc(\nu)$ the collection of $\P\in {\rm Prob}(\Omega)$ under which\footnote{Recall that by definition of $u$ it can be that $u_t\notin \pi(U)$ on a null set.
	The condition in $(i)$ and the introduction of the cemetery point $\iota$ guarantee $u_t$ is always well defined, even on that null set. Condition $(ii)$ should be understood in the sense of the associated martingale problem, see {\rm \cite{stroock1997multidimensional,djete2019mckean}}} 
\begin{enumerate}[label=$(\roman*)$, ref=.$(\roman*)$,wide,  labelindent=0pt]
\item \label{def.controlrule.1} $\P[u_t\in \Uc]=1$, for $ \lambda \ae\, t\in [0,T]$ and
\[
 \E^\P\bigg[  \sup_{t\in [0,T]} \|X_t\|^p \bigg] +  \E^\P\bigg[\int_0^T \|u_t\| ^p\d t \bigg]<\infty.% \text{ for some } {\bf u}_0\in U.
\]
%{\color{blue}maybe replace $\rho$ by $\delta$, in view of $c(x,\rho)$ coming below} \camilo{Check everywhere and remove}
\item \label{def.controlrule.2} For all $\phi \in C_b^2(\R^{m+d})$, the process $M^\phi: [0,T]\times \Omega \longrightarrow \R $ is an $(\F,\P)$--local martingale, with $M^\phi$ defined by
\[
M_t^\phi\coloneqq \phi(X_t,W_t)-\int_0^t (\Lc \phi) (s,X,W, \mu_s,u_s)\diff s,\, t\in [0,T].
\]
% $W$ is a $(\F,\P)$--Brownian motion and
%\[ X_t = X_0 + \int_0^t \big( \alpha_s+b(r, X,   \mu_s )\big)\diff s + \int_0^t  \sigma_s  \diff W_s ,\; t\in [0,T],\, \P\as \]
%{\color{blue}are we integrating $\beta$ ore $\sigma$?} \camilo{$\sigma$ always, check}
\item  \label{def.controlrule.3} $\P\circ X_0^{-1}=\nu$, $\P\circ W_0^{-1}=\delta_0$, the pair $(X_0,W)$ is $\P$-independent of $\overline{\Gc}_T$, and, for every
$t\in[0,T]$,
%\begin{equation}\label{eq:fixedpoint}
\[
	\bar\mu_t =
	\P^{\overline{\Gc}_T}_{\omega}\circ (X_{t\wedge\cdot},U_{t\wedge\cdot},W)^{-1},
	\P\as
\]

\end{enumerate}
\end{definition}

We now define 
\begin{align}\label{def.prob.can}
\Jc(\P)\coloneqq \E^\P\bigg[\int_0^T f(t,X,u_t, \mu_t)\diff t\bigg],\,  \P\in \Pc(\nu), \text{ and, } \Pc_c(\muin,\mufin)\coloneqq \{\P\in \Pc(\muin):\P\circ X_T^{-1}\=c\mufin\}.
\end{align}

% {\color{cyan}
% So the problem becomes
% \begin{equation*}
% 		\bar V = \inf_{\P\in \Pc(\muin)}\Big( \Jc(\P) + g(\Wc_c(\P\circ X_T^{-1},\mufin)) \Big)
% \end{equation*}
% I am still a bit unlear on the fact that the terminal law $\P\circ X^{-1}_T$ is controlled, but the "running" law $\mu_t$ is not. }

\begin{remark}\label{rmk.equiv} 
\begin{enumerate}[label=$(\roman*)$, ref=.$(\roman*)$,wide,  labelindent=0pt]
\item We remark that the sets $\Gamma(\muin)$ and $\Pc(\muin)$ are in correspondence, see {\rm \cite[Lemma 4.3]{djete2019mckean}}. 
	Indeed, for any $\gamma\in \Gamma(\muin)$ we have that $\overline\P^\gamma:=\P^\gamma\circ(X^\gamma,U^\gamma,W^\gamma,\mu^\gamma)^{-1}\in \Pc(\muin)$, with \( \widehat\mu^{\,\gamma}_t\coloneqq \Lc^{\P^\gamma}  ( (X^\gamma_{t\wedge\cdot},U^\gamma_{t\wedge\cdot},W^\gamma) \, |\,\mathcal G^\gamma_T ),\, t\in[0,T] \).
%see also {\rm \cite[Remark 5]{djete2019mckean}}. 
Conversely, given $\overline\P\in \Pc(\muin)$, there exists a weak control rule $\gamma\in \Gamma(\muin)$ such that $\P^\gamma\circ(X^\gamma,U^\gamma,W^\gamma,\mu^\gamma)^{-1}=\overline\P$. 
In particular, this shows that 
\[
J(\gamma)=\Jc(\P).
\]
\item In complete analogy with \Cref{rem:defcomments}, we denote by $\Pc_{\rm det}(\nu)$ the collection of elements in $\Pc(\nu)$ with trivial $\bar \G$.
\end{enumerate}
\end{remark}

\subsection{Existence}

The proof of the existence claim will follow as a direct consequence of the following lower semicontinuity of the function $V_c$.
\begin{lemma}
\label{Lem:v-LSC}
	Let {\rm Assumptions \ref{assump.costs} and \ref{assump.bf}} hold. The function $\Pc(\R^\xdim)\times \Pc(\R^\xdim)\ni (\mu,\nu)\longmapsto V_c(\mu,\nu) \in \R \cup\{+\infty\}$ is lower semicontinuous for the topology of weak convergence.
\end{lemma}
\begin{proof}
	The proof follows some ideas of the proof of \cite[Lemma 3.1]{mikami2008optimal}.
	We argue in three steps.\medskip
	
	{\bf Step 1.} Let $(\muin^n,\mufin^n)$ be a sequence of elements of $\Pc(\R^\xdim)\times \Pc(\R^\xdim)$ converging to $(\muin,\mufin)$.
	Without loss of generality, we assume that $\liminf_{n\to \infty}V_c(\muin^n,\mufin^n) <\infty$.
	Thus, up to a subsequence, $(V_c(\muin^n,\mufin^n))_{n\in \N}$ is bounded.
	In particular, for each $n\ge1$ there is $\gamma_n\in \Gammac(\muin^n,\mufin^n)$ such that
	\begin{align}
	\label{eq:lsc.eps-optimal}
		V_c(\muin^n,\mufin^n) \ge J(\gamma_n) - \frac1n
	\end{align}
	and combined with the coercivity condition Assumption \ref{assump.bf}\ref{assump.bf.2}, we have 
	%{\color{red}If I understand the setting properly, then $(\alpha^n,\beta^n, u^n,\sigma^n)$ below should be $(\alpha,\beta, u,\sigma)$. But the proof should be the same because $A^n$ below will still depend on $n$ via $\mu^n$.}
	\begin{equation*}
		\sup_{n\in \N}\E^{\P^{\gamma_n}}\bigg[\int_0^T\|u^{\gamma_n}_s \|^{\hat p} \diff s \bigg] <\infty.
	\end{equation*}
	Thus, on $(\Omega^{\gamma_n},\Fc^{\gamma_n},\F^{\gamma_n},\P^{\gamma_n})$, with $u^{\gamma_n}:= (\alpha^{\gamma_n},\sigma^{\gamma_n}{\sigma^{\gamma_n}}^\t)$, thanks to \Cref{assump.bf}\ref{assump.bf.3}, we have that
	\begin{align}\label{eq.assZheng}
	X^{\gamma_n}_0,\; \int_0^T \| \alpha_s^{\gamma_n} + b(s, X^{\gamma_n}, \mu^{\gamma_n}_s)\|^{\hat p} \diff s,\; \text{and,} \; 	\int_0^T \|\sigma^{\gamma_n}_s\sigma^{\gamma_n \t}_s\|^{\hat p} \diff s,
	\end{align}
	are uniformly bounded in probability. 
	Recall \eqref{eq.repofcontrols} and note that thanks to \Cref{rmk.equiv}, we have that 
	\[
	\P^n\coloneqq \P^{\gamma^n}\circ(X^{\gamma^n}, U^{\gamma^n}, W^{\gamma^n}, \bar \mu^{\gamma^n})^{-1}\in \Pc(\muin^n).
	\]
	Note that by definition $(\P^n)_{n\geq 1}\subseteq \Pc_c(\muin^n,\mufin^n)$.\medskip
	
	{\bf Step 2.} We claim that $(\P^n)_{n\geq 1}\subseteq {\rm Prob}(\Omega)$ is tight and that any of its limit laws is in $\Pc(\muin,\mufin)$.
	We first note that for all $n\geq 1$, $X= X_0+ A + M^n, \P^n\as$, and $\langle X\rangle=\langle M^n\rangle=\Sigma$, $M^n$ denoting the $\P^n$-local martingale part,
	%dependence of the stochastic integral on $\P^n$, 
	for absolutely continuous processes $A$ and $\Sigma$ on $\Omega$ defined by \(
A_\cdot \coloneqq \int_0^\cdot \big( \alpha_s + b(s, X,  \mu_s )\big) \diff s, \; \Sigma_\cdot:=\int_0^\cdot \sigma_s \sigma^\top_s \diff s.\)\medskip
%	we recall \eqref{eq.AS} and the superscript on .
%	Moreover, $A$ and $\langle M^n\rangle=\Sigma$ are absolutely continuous with respect to the Lebesgue measure.
%	\medskip
	
Notice now that
	\[
		\sup_{n\in \N}\E^{\P^{n}}\bigg[\int_0^T\|u_s \|^{\hat p} \diff s \bigg] = \sup_{n\in \N}\E^{\P^{\gamma_n}}\bigg[\int_0^T\|u^{\gamma_n}_s \|^{\hat p} \diff s \bigg] <\infty,
	\]
	%{\color{blue}it is probably better to use the notation $\mu^\P$ and $u^\P = (\alpha^\P, \beta^\P)$ in the definition in page 2 if we mean that $\alpha,\beta$ depend on $\P$.}\medskip
	%\camilo{Let $\widehat \P^n:= \P^{\gamma^n}\circ(X^{\gamma^n}, A^{\gamma^n}, M^{\gamma^n}, \Sigma^{\gamma^n})^{-1}$ and identify it with $\P^n:= \P^{\gamma^n}\circ(X^{\gamma^n}, U^{\gamma^n}, W^{\gamma^n}, \bar \mu^{\gamma^n})^{-1}$}
	and, thanks to \Cref{assump.bf}\ref{assump.bf.3}, the analogous version of \eqref{eq.assZheng} holds on $\Omega$.
	By \cite[Theorem 3]{zheng1985tightness}, $\widehat \P^n:= \P^{n}\circ(X, A, M, \Sigma)^{-1}$, where $M=X-X_0-A$, is tight for the uniform topology on $\widehat\Omega:=(\Cc^\xdim)^3\times \Cc^{\xdim\times \xdim}$.\footnote{the process $(X, A, \Sigma, M)$ is tight on $(\widehat\Omega,\P^n)$.}		
	With this, we are ready to argue the claim. 
	First note that, by definition of $\P^n\in\Pc(\nu_0)$, $\P^n\circ W^{-1}$ is tight and $U_t=\int_0^t \pi(u_t)\diff t, \P^n\as$, and thus, $\P^n\circ U^{-1}$ is tight.
	Now note that, by the tower property and \Cref{def.controlrule}\ref{def.controlrule.3}, $\E^{\P^n}[\bar\mu_T]=\P^n\circ (X,U,W)^{-1}$, and thanks to \cite[Lemma A.2]{lacker2015mean}, $\P^n\circ \bar \mu_T^{-1}$ is tight in $\Pc(\Pc(\bar\Omega))$.
	Thanks to \Cref{lem:flowdetermined}, we have that $\bar\mu=\mu_T\circ (\bar  X_{t\wedge \cdot},\bar  U_{t\wedge \cdot},\bar W)^{-1}$, $\P^n\as$, for every $n$.
	Hence, $\P^n\circ \bar \mu^{-1}$ is tight in $\Pc(C[0,T],\P(\bar\Omega))$.
	Thanks to \cite[Lemma A.2]{lacker2015mean} again, we deduce that $(\P^n)_{n\geq1}$ is tight.\medskip

Let $\P$ be a limit law of $(\P^n)_{n\geq1}$. 
	Notice that $\P$ induces $\widehat\P\in {\rm Prob}(\widehat\Omega)$ which is a limit law of $(\widehat \P^n)_{n\geq1}$.
	Letting $(\widehat X, \widehat A, \widehat M, \widehat \Sigma)$ denote the coordinate process on $\widehat\Omega$, \cite[Theorem 3]{zheng1985tightness} ensures that under $\widehat \P$, 
	%{\color{blue}first point to discuss: It might be easier to have $M$ as canonical rather than W?}
	%{\color{cyan}I still don't understand. What is the canonical process? why is it different from our defined extended space? }
	$\widehat X_t = \widehat X_0 + \widehat A_t + \widehat M_t, \widehat\P\as,$ $\widehat M$ is a continuous local martingale, and, $\widehat A_t$ and $\widehat \Sigma_t: = \langle \widehat M\rangle_t$ are absolutely continuous, with densities $ (\hat a , \hat \beta)$ such that
	\begin{equation*}
		\int_0^T\big( \|\hat a_s\|^{\hat p}+\|\hat \beta_s\|^{\hat p}\big) \diff s <\infty,\; \widehat  \P\as
	\end{equation*}
	
%{\color{cyan}Also, I was thinking that what we need is to have tighness of the probability measures $\P^n:= \P^{\gamma^n}\circ(X^{\gamma^n}, A^{\gamma^n}, \Sigma^{\gamma^n}, M^{\gamma^n}, \mu^{\gamma^n})^{-1}$ on the extended space.}	
	Let us now define $\hat \mu_t:= \widehat \P\circ (\widehat X_{t\wedge\cdot})^{-1}$, $\hat u_t:=(\hat \alpha_t,\hat \sigma_t)$, where $\hat \sigma  =\hat \beta^{1/2}$, $\hat \alpha_t:= \hat a_t - b(t, \widehat X,  \hat \mu_t )$ and note that $\hat u$ is in correspondence with the process $u$ defined on $\Omega$.
	Since $\Uc$ is closed we get from the Portmanteau theorem that $n \int_{( t-1/n)\vee 0}^t  u_s \diff s \in \Uc,  \P\as$ for every $t\in [0,T]$ and $n\geq 1$.
	Thus, $ \P[  u_t \in \Uc]=1$ for $\lambda \ae\, t\in [0,T]$.
	Using \Cref{assump.bf}\ref{assump.bf.3}, we deduce that $\P$ satisfies \Cref{def.controlrule}\ref{def.controlrule.1}.
	Moreover, it follows from Itô's formula that for any $\hat \phi \in C_b^2(\R^m)$, the process
	\[
	\widehat M_t^{\hat \phi}=\hat \phi(\widehat X_t)-\int_0^t \big(  ( \hat\alpha_s+b(s,\widehat X,\hat \mu_s) )\cdot \nabla \hat\phi (\widehat X_s)+\frac12 {\rm Tr}[\hat \beta_s \nabla^2\hat\phi(\widehat X_s)]\big) \diff s, \; t\in [0,T],
	\]
	is an $(\widehat  \F,\widehat\P)$-martingale, where $\widehat\F$ denotes the filtration generated by $(\widehat X, \widehat A, \widehat M, \widehat \Sigma)$.
	By \cite[Theorem 4.2.1]{stroock1997multidimensional}, it follows that $ \widehat \P$ is a solution to the martingale problem associated with \eqref{eq.dyn} on $ \widehat\Omega$.
	Consequently, recall \Cref{rmk.equiv} and $\widehat \P=  \P\circ ( X, A, M,\Sigma)^{-1}$, $\P$ solves the martingale problem associated with \eqref{eq.dyn} on $\Omega$.
	As in the proof of \cite[Theorem 4.5.2]{stroock1997multidimensional}, we deduce that on $\Omega$, $W$ is a $(\F,\P)$--Brownian motion and \Cref{def.controlrule}\ref{def.controlrule.2} holds.
	Thus, for any $\phi \in C_b^2(\R^{m+d})$ the process $M^\phi$ is an $(\F,\P)$-local martingale. 
	That is, \Cref{def.controlrule}\ref{def.controlrule.2} holds.
	Verifying the independence of $(X_0,W)$ and $\bar\Gc_T$ as well as the identity in \Cref{def.controlrule}\ref{def.controlrule.3} follows as in \cite[Section 4.3.2.2]{djete2020mckean2}.

	\medskip
	
	Finally, we claim that $\P \in \Pc_c(\muin,\mufin)$
	%Note that by martingale representation, see \cite[Theorem 3.4.2]{karatzas1991brownian}, on a possibly enlarged space, there is $\widehat W$ such that $\widehat M=\int_0^\cdot \hat \sigma_t \d \widehat W_t$. With a slight abuse of notation, we denote $(\widehat \P, \widehat X,\widehat U,\widehat W,\hat \mu )$ the version in the enlarged space.
	%\medskip
	%We now let $\P:=\widehat \P\circ(\widehat X,\widehat U,\widehat W,\hat \mu)^{-1} $ and claim that $\P\in \Pc(\mu)$. \camilo{To do.} In particular, in the space $\Omega$ it holds that
	%with $\sigma$ such that $\beta = \sigma\sigma^\top$.
	indeed, since $\P^{n} \circ (X)^{-1}\longrightarrow \P\circ X^{-1}$ and $\P^{n}\in \Pc_c(\muin^n,\mufin^n)$, it follows on the one hand that $\P\circ X^{-1}_0  = \muin$ and on the other hand, since by \cite[Theorem 2.9]{backhoff2019existence} the function $\Wc_c$ is lower semicontinuous, we have that
%	{\color{blue} I hope the function $\Wc$ is continuous. any reference?)}
	\begin{equation*}
		0 = \liminf_{n \to \infty} \Wc_c(\P^n\circ X^{-1}_T, \mufin^n) \geq \Wc_c(\P\circ X^{-1}_T, \mufin),
	\end{equation*}
	which shows that $\P \in \Pc_c(\muin,\mufin)$.
	%{\color{cyan}The two parts of the proof are disconnected. Fixes are still needed.}

	\medskip

	{\bf Step 3.} % THANKS
	Let us now argue the statement of the Lemma. 
	Let $(\P^n)_{n\geq 1}\subseteq {\rm Prob}(\Omega)$ and $\P$ be as in the previous step. 
	Let us now show that $\liminf_{n\to \infty} \Jc(\P^n) \ge \Jc(\P)$.
	Arguing as in the proof of \cite[inequality (3.17)]{mikami2002optimal}, we obtain from \Cref{assump.bf}\ref{assump.bf.i} that
	\begin{align*}
		\int_0^Tf(s, {\bf x}, u_s, \mu) \diff s \ge \frac{1}{1 + \Delta_tf(\varepsilon)}\bigg( \int_0^{T - \varepsilon} f(s,{\bf x}, \eps^{-1}(\eta_{s+\eps}-\eta_s), \mu )\diff s - T \Delta_tf(\varepsilon)\bigg)
	\end{align*}
	for every $(t,{\bf x}, \mu)\in [0,T]\times \R^\xdim\times\Pc(\R^\xdim)$, and $\R^\xdim\times\S^{\xdim}_+$-valued process $u$, with $\eta_t = \int_0^t u_s\diff s $.
	% \camilo{I think the inequality is missing a braket
	% \begin{align*}
	% 	\int_0^Tf(s, {\bf x}, u_s, \mu_s) ds \ge \frac{1}{1 + \Delta_tf(\varepsilon)} \bigg( \int_0^{T - \varepsilon} f(s,{\bf x}, \eps^{-1}(\eta_{s+\eps}-\eta_s), \mu_s )ds - \Delta_tf(\varepsilon)\bigg)
	% \end{align*}
	% }
	Recall that for any $\P\in \Pc(\muin)$, $\P[u_t \in \Uc]=1$ for $\lambda \ae\,  t\in [0,T]$ and $\eta_\cdot \coloneqq \int_0^\cdot u_s\diff s$ is well defined on $\Omega$.
	Thus, 
	%\camilo{this is fixed in the enlarged space} \sout{putting $\mu^n:= \P^n\circ X^{-1}$}, 
	since $\P^n\in \Pc(\muin)$, we have 
	%{\color{red} \sout{here also $(A^n,U^n)$ might be fixed?}}\camilo{Yes}
	\begin{align*}
		\Jc(\P^n) &\ge \frac{1}{1 + \Delta_tf(\varepsilon)}\bigg(  \E^{\P^n}\bigg[\int_0^{T - \varepsilon} f(s,X , \varepsilon^{-1}(\eta_{s+\eps}-\eta_s) \big) , \mu ) \diff s \bigg] - T \Delta_t f(\varepsilon)\bigg).
				%& \ge \frac{1}{1 + \Delta_tf(\varepsilon)}\E^{\bar \P^n}\bigg[\int_0^{T - \varepsilon} f(s,X , \frac{1}{\varepsilon}(A_{s+\varepsilon} - A_s), \frac{1}{\varepsilon}(U_{s+\varepsilon} - U_s), \mu_s )ds \bigg] - \Delta_Tf(\varepsilon)
	\end{align*}
	%where $\bar\P^n:= \P^n\circ(X, A^n, U^n,\mu^n)^{-1}$. {\color{blue}the second inequality might not be okay. We might need stronger regularity condition $f$ in $\mu$ to get the limit.}
	Now, since $f$ is continuous and nonnegative, Fatou's Lemma gives that $(X,\eta,\mu)\longmapsto \int_0^{T-\eps}  f(s,X , \varepsilon^{-1}(\eta_{s+\eps}-\eta_s) \big) , \mu ) \diff s $ is lower-semicontinuous, and thus, for every $\varepsilon>0$,
	\begin{equation*}
		\liminf_{n\longrightarrow\infty} \Jc (\P^n)  \ge \frac{1}{1 + \Delta_tf(\varepsilon)}\bigg(  \E^{\P}\bigg[\int_0^{T - \varepsilon} f(s,X,   \varepsilon^{-1}(\eta_{s+\eps}-\eta_s)  ,\mu_s )ds \bigg] - T \Delta_tf(\varepsilon)\bigg)
 	\end{equation*}
	Now, taking the limit as $\varepsilon\longrightarrow 0$, which is justified again by Fatou's lemma, we have thanks to \Cref{assump.bf}\ref{assump.bf.i} that
 	\begin{equation*}
 		\liminf_{n\longrightarrow \infty}\Jc(\P^n) \ge \E^{\P}\bigg[\int_0^Tf(s, X, u_s, \mu_s)ds\bigg] = \Jc(\P).
 	\end{equation*}
 	Using \eqref{eq:lsc.eps-optimal}, that $J(\gamma^n)=\Jc(\P^n)$, and $\P \in \Pc(\mu,\nu)$ yields
 	\begin{equation*}
 		\liminf_{n\longrightarrow \infty}V_c(\muin^n,\mufin^n)\ge V_c(\muin,\mufin),
 	\end{equation*}
 	from which we deduce lower semicontinuity of $V_c$.
\end{proof}

\begin{proof}[Proof of {\rm Theorem \ref{thm:existence}}]
	Since $\Pc_c(\muin,\mufin) \neq \emptyset$, then for every $n\ge1$ there is $\P^n\in \Pc_c(\muin,\mufin)$ such that
	\begin{equation*}
			V_c(\muin,\mufin) \ge \Jc(\P^n) - \frac1n.
	\end{equation*}	
	As argued in the proof of Lemma \ref{Lem:v-LSC} above, there is $\P \in \Pc(\muin,\mufin)$ such that $\Jc(\P) \le \liminf_{n\to\infty}\Jc(\P^n)$,
	%Thus, we have $V_c(\mu,\nu) \ge \Jc(\P)$.
	%This shows that $\P$ is optimal.
	and, by \Cref{rmk.equiv}, there is a weak control rule $\gamma\in \Gamma(\muin,\mufin)$ such that $J(\gamma)=\Jc(\P)$.
	Thus, $\gamma$ is optimal for Problem \eqref{eq:problem}.
\end{proof}

\subsection{Convex dual representation}
\label{sub:convex_dual_representation}

Before proving convexity of the function $\mufin\mapsto V_c(\muin, \mufin)$, let us discuss the nonconvexity situation announced in \Cref{sec.formulation}.
%\subsection{On the convexity of $V_c^{\rm det}(\muin,\cdot)$}\label{sec.conx.fails}
In fact, the McKean-Vlasov transportation problem over deterministic-flow rules is not necessarily convex as observed from the example below.

\begin{example}\label{prop:noconvexdet}
	Let $\xdim =d=1$, $\sigma\equiv1$, $\muin\coloneqq \mathcal N(0,1)$, $f\equiv  \tfrac12|a|^2$, \( b(t,x,\mu)\coloneqq \int_{\mathbb R}y\,\mu(t)(\mathrm dy)-x_t,\, c(x,\rho)\coloneqq \int_\R \1{x\neq y}\rho(\d y)\).
	Moreover, let $ \upsilon >0$, $\theta^2>0$, and \( \mufin^{+}\coloneqq \mathcal N( \upsilon,\theta^2),\; 
 	\mufin^{-}\coloneqq \mathcal N(- \upsilon,\theta ^2),\;
	\bar{\nu}_T \coloneqq \tfrac12\mufin^{+}+\tfrac12\mufin^{-}\).
Then, {\rm Assumptions~A and~B} hold, and for a constant $c_{\theta,T}>0$ independent of $\upsilon$, we have that
\[
V^{\rm det}_c(\nu_0,\mufin^{+})\leq 
	\frac{\upsilon^2}{2T}+c_{\theta,T},\; 
	V^{\rm det}_c(\nu_0,\mufin^{-})
	\leq 
	\frac{\upsilon^2}{2T}+c_{\theta,T} ,\mbox{ and, } 
	V^{\rm det}_c(\nu_0,\bar{\nu}_T) 
	\geq C_{\upsilon,\theta,T},
	\mbox{ where }
	C_{\upsilon,\theta,T}\coloneqq \frac{(\sqrt{\upsilon^2+\theta^2}-\e^{-T}-1/\sqrt2\big)_+^2}{1-\e^{-2T}}.
\]
In particular, there is $\upsilon$ large such that the map $\nu_T\longmapsto V^{\rm det}_c(\nu_0,\nu_T)$ is not convex.
\end{example}

\begin{proof}
Let us first verify the assumptions. 
	Note $c$ is bounded, affine, and thus convex in $\rho$, satisfies $c(x,\delta_x)=0$, and is jointly l.s.c.
	%\ since $(x,\rho)\mapsto\rho(\{x\})$ is u.s.c.\ (if $x_n\to x$ and $\rho_n\to\rho$ weakly, then $\limsup_n\rho_n(\{x_n\}) \le\limsup_n\rho_n(\bar B_\varepsilon(x))\le\rho(\bar B_\varepsilon(x)) \to\rho(\{x\})$ as $\varepsilon\downarrow0$).
	Moreover, $|b(t,x,\mu)|\le|x_t|+(\int|x'_t|^2\mu(t)(\d x^\prime))^{1/2}$. 
	%The requirement $0\in\mathcal U$ can be met by taking $\mathcal U=\mathbb R\times\{0,1\}$ without affecting any of the bounds.
	%Indeed, the upper bounds use $s\equiv1$; the lower bound only improves if $s\in\{0,1\}$).
	Note also that the coercivity condition holds. %on this $\mathcal U$ with $\hat p=2$ since the volatility component is bounded there
\medskip

	In the remainder, given a rule $\gamma$ with deterministic flow and	$J(\gamma)<\infty$, we write $m_t\coloneqq\mathbb E[X_t]$, $\bar\alpha_t\coloneqq \mathbb E[\alpha_t]$, and $Z\coloneqq X-m$.
	To streamline the presentation, we drop the superscript $\gamma$. 
	Since $b(t,X,\mu_t)=m_t-X_t$, we have
\begin{equation}\label{eq:decomp}
	m_t = \int_0^t\bar\alpha_r\,\mathrm dr,
	\; 
	\d Z_t
	=
	\big(\alpha_t-\bar\alpha_t-Z_t\big)\diff t+\d W_t,
	\mbox{ and, }
	\E \bigg[ \int_0^T |\alpha_t|^2\diff t\bigg] =\int_0^T\!|\bar\alpha_t|^2\mathrm dt +\E\bigg[  \int_0^T |\alpha_t-\bar\alpha_t|^2\diff t\bigg],
\end{equation}
	where the first identity follows from the state equation and the cost decomposition follows from the variance formula.\medskip

	Let us first obtain the upper bounds.
	Let $k\in \R$, to be given, and put $\alpha_t^+ \coloneqq \frac \upsilon T+k (X_t-m_t)$. 
	Note that the McKean-Vlasov SDE with control $\alpha^+$ has a unique strong solution $X^+$, whose flow is deterministic, and thus induces an element of $\Gamma_{\rm det}(\nu_0)$. 
	By \eqref{eq:decomp}, $m_t^+ \coloneqq \E[X^+_T] = \frac \upsilon T  t$, hence $m_T= \upsilon $, and $Z^+ = X^+-m^+$ satisfies $\d Z_t^+=(k-1)Z_t^+ \d t+\d W_t$, so it is an Ornstein--Uhlenbeck process.
	Thus, $v_k(t)\coloneqq \mathbb E[Z_t^2]$ solves $v^\prime_k=2(k-1)v_k+1$, $v_k(0)=1$. 
	Note that $k\longmapsto v_k(T)$ is continuous with $v_k(T)\longrightarrow 0$ as $k\longrightarrow -\infty$ and $v_k(T)\longrightarrow\infty$ as $k\longrightarrow  +\infty$.
	Thus, there is $k=k_{\theta,T}$ with $v_k(T)=\theta ^2$. 
	Then, $X_T=m_T^+ +Z_T\sim\mathcal N(\upsilon ,\theta^2)=\mufin^{+}$, and, using $\mathbb E[Z_t]=0$ and \eqref{eq:decomp},
	\[
		J(\gamma^+)
		=\frac12\E\bigg[ \int_0^T \big |\frac \upsilon T+k Z_t \big|^2 \diff t \bigg]
		=
		\frac{\upsilon ^2}{2T}+\frac{k^2}{2}\int_0^T v(t)\diff  t \eqqcolon \frac{\upsilon^2}{2T}+c_{\theta,T}.
	\]
	The fact that $J(\gamma^-)= \frac{\upsilon^2}{2T}+c_{\theta,T}$ follows similarly with $\alpha_t^- \coloneqq  -\frac \upsilon T+k (X_t-m_t)$.\medskip

	We now obtain the lower bound on $V^{\rm det}_c(\nu_0,\bar{\nu}_T)$. 
	Let $\gamma\in\Gamma^{\rm det}_{c}(\nu_0,\bar{\nu}_T)$ with $J(\gamma)<\infty$, otherwise the result holds. 
	Note the terminal constraint implies $m_T=0$ and $\E[Z_T^2]=\mathrm{Var}(\bar{\nu}_T)=\upsilon^2+\theta^2$ and that we might write \eqref{eq:decomp} as
\[
Z_T=e^{-T}Z_0+\int_0^T e^{-(T-r)}\big(\alpha_r-\bar\alpha_r\big)\mathrm dr
+\int_0^T e^{-(T-r)}\,\mathrm dW_r .
\]
Moreover, by definition of $\muin$, $\E[ |Z_0|^2]=1$, and $\int_0^T\! \e^{-2(T-r)}\d r =(1-e^{-2T})/2\leq  1/2$.
	By Cauchy--Schwarz,
\[
\sqrt{\upsilon^2+\theta^2}
	 \leq \e^{-T} + \frac1{\sqrt2}
	+ \sqrt{\frac{1-e^{-2T}}{2}}
	\bigg(\E\bigg[  \int_0^T |\alpha_r-\bar\alpha_r|^2\diff r\bigg] \bigg)^{1/2},
\]
and consequently, by \eqref{eq:decomp},
\[
J(\gamma)\geq 
	\dfrac12 \E\bigg [  \int_0^T |\alpha-\bar\alpha|^2 \diff t\bigg]
	\geq
	C_{\upsilon,\theta,T}.
\]
Since $\gamma$ was arbitrary, this proves the lower bound $V^{\rm det}_c(\nu_0,\bar{\nu}_T)  \geq C_{\upsilon,\theta,T}$.\medskip

To conclude, note that since $1-\e^{-2T}< 1/(2T)$, for $\upsilon$ large enough, the leading coefficient of $\upsilon^2$ in $C_{\upsilon,\theta,T}$ exceeds that of the averaged upper bounds. That is
\[
V^{\rm det}_c(\muin,\bar{\nu}_T )
	 >
	\dfrac12 V^{\rm det}_c(\muin,\mufin^{+})+\dfrac12V^{\rm det}_c(\muin,\mufin^{-}).\qedhere
\]
\end{proof}

We now turn to the proof of the convex dual representation result.
This will follow as a consequence of the Fenchel--Moreau theorem.
To this aim, we will show that the function $\mufin\longmapsto V_c(\muin,\mufin)$ is convex and lower semicontinuous, then compute its convex conjugate.
\begin{proposition}
\label{prop.convex.V}
	Let $\muin\in \Pc_p(\R^\xdim)$. Then, $\Pc_c(\muin,\mufin)$ is convex and the mapping $\Pc_p(\R^\xdim)\ni\mufin\longmapsto V_c(\muin,\mufin)$ is convex.
\end{proposition}
\begin{proof}
	Let us first note that $\Pc(\muin)$ is convex.
	Indeed, this follows from {\rm \cite[Theorem $1.(i)$]{djete2020mckean2}}, and \Cref{rmk.equiv}.
	Let $\theta \in (0,1)$, $\mufin^i\in \Pc_p(\R^\xdim)$, $i=1,2$ be given, and $\P^i \in \Pc(\muin,\mufin^i)$ be given from $\gamma^i\in \Gamma(\muin^i,\mufin^i)$ as in \Cref{rmk.equiv}.
	Let $\overline \P = \theta \P^1 + (1 - \theta)\P^2$ and define ${\overline \nu}_T$ similarly.
	We claim that $\overline \P\in \Pc_c(\muin,{\overline\nu}_T)$.
	Since $\Pc(\muin)$ is convex, it only remains to show %$\overline \P\circ X^{-1}_0$ and 
	that $\overline \P\circ X^{-1}_T\=c{\overline\nu}_T$. 
	% and notice that $\overline \P\circ X^{-1}_0={\overline\nu}_0$ is clear.
	Let us recall from \cite[Theorem 1.3]{backhoff2019existence} the duality formula for $\Wc_c$
	\begin{align}\label{eq.dualitywot}
		\Wc_c(\P\circ X^{-1}_T , \mufin) = \sup_{\varphi \in C_{b,p}(\R^d)}(\langle Q_c\varphi,\P\circ X^{-1}_T\rangle -\langle \varphi,\mufin\rangle  ).
	\end{align}
	Thus, since $\overline \P\circ X^{-1}_T=\theta \P^1\circ X^{-1}_T+(1-\theta) \P^2\circ X^{-1}_T$, and $c$ is nonnegative, $\overline \P\circ X^{-1}_T\=c{\overline\nu}_T$ follows since we have that 
	\[
	0\leq \Wc_c(\overline \P\circ X^{-1}_T,{\overline \nu}_T ) \le \theta \Wc_c(\P^1\circ X^{-1}_T, \mufin^1) + (1-\theta)\Wc_c(\P^2\circ X^{-1}_T, \mufin^2 ) = 0.
	\]
	
	Moreover, we have
	\begin{align*}
		V_c({\overline\nu}_0,{\overline\nu}_T) \le \Jc(\overline \P) = \E^{\bar\P}\bigg[ \int_0^T f(t,X, u_t,\mu_t)\diff t\bigg] 
		& = \theta \E^{\P^1}\bigg[ \int_0^T f(t,X, u_t,\mu_t)\diff t\bigg] + (1 - \theta)\E^{\P^2}\bigg[ \int_0^T f(t,X, u_t, \mu_t )\diff t\bigg]\\
		& = \theta  \Jc(\P^1) + (1 - \theta) \Jc( \P^2) .
	\end{align*}
	This yields the convexity of $\mufin\longmapsto V_c(\muin,\mufin)$ since $\P^1$ and $\P^2$ were chosen arbitrarily.
\end{proof}
	
\begin{lemma}\label{lem:disint}
Let $\nu\in \Pc_p(\R^\xdim)$ and $\P\in\Pc(\nu)$. %and let $(\P_\omega^{\overline{\Gc}_T})_{\omega\in\Omega}$ be a family of r.c.p.d.\ of $\P$ knowing $\overline{\Gc}_T$. 
	Then, % for $\P\ae$ $\omega$:
	$\P_\omega^{\overline{\Gc}_T}\in\Pc_{\rm det}(\nu),\,\P\as$, and, $\bar\mu=\bar\mu(\omega), \, \P_\omega^{\overline{\Gc}_T}\as$ 
	Consequently,
	\(
	\Jc(\P)=\langle \Jc(\P_\cdot^{\overline{\Gc}_T}), \P\rangle,\;
	\P\circ X_T^{-1} =\P\otimes (\P_\cdot^{\overline{\Gc}_T}\circ X_T^{-1}).
	\)
%Conversely, if $(\P_\omega)_{\omega}$ is a measurable family in $\Pc_{\rm det}(\nu)$ and $\P\in \Pc(\Omega)$ is such that
%	\begin{equation}\label{eq:mixmoment}
%	\int_\Omega \Big(\E^{\P_\omega}\big[\sup_{t\in[0,T]}\|X_t\|^p\big]
%	+
%	\E^{\P_\omega}\Big[ \int_0^T \|u_t\|^p\,\d t\Big]\Big)
%	\P(\d \omega )<+\infty,
%	\end{equation}
%then, $\bar\P\coloneqq \P\otimes \P_\cdot\in \mathcal P(\nu)$.
In particular, for every $\varphi\in C_{b,p}(\R^\xdim)$, \(\Psi^\varphi(\nu)
=\inf_{\gamma \in \Gamma_{\rm det}(\nu)}J^\varphi(\gamma)\).
%\begin{align}
%\Psi^\varphi(\nu)
%=\inf_{\mathbb P\in\mathcal P_{\rm det}(\nu)}\mathbb E^{\mathbb P}\bigg[\int_0^T
%f(s,X,u_s,\mu_s)\,\mathrm ds+Q_c\varphi(X_T) \bigg].
%\end{align}
\end{lemma}
\begin{proof}
Throughout we write $\xi_t:=(X_{t\wedge\cdot},U_{t\wedge\cdot},W)$ and $\mathbb F^\circ\coloneqq(\mathcal F^\circ_t)_{t\in[0,T]}$, $\mathcal F^\circ_t\coloneqq \sigma(X_{t\wedge\cdot},U_{t\wedge\cdot},W_{t\wedge\cdot})$, the filtration generated by the first three coordinates.
	Let us first verify that $\P_\omega$ satisfies \Cref{def.controlrule} for $\P\ae \, \omega\in \Omega$. We argue in four steps. %For this, on each of the next steps we prove that:\medskip
\medskip

{\bf Step 1}. There is $N\subseteq \Fc, \P(N)=0$, such that \(\bar\mu_t(\omega)=\mathbb P_\omega\circ\xi_t^{-1},\;\mbox{for all } t\in[0,T],\, \omega \in  N^c\). 
	Since $\bar\mu$ is $\overline{\Gc}_T$--measurable, there is a
$\mathbb P$--null set $N_1$ with $\bar\mu=\bar\mu(\omega)$,
$\mathbb P_\omega$--a.s., for $\omega\notin N_1$. 
	Next, let $t\in[0,T]$. 
	By \Cref{def.controlrule}\ref{def.controlrule.3} there is a $\P$--null set $N(t)$ such that $\bar\mu_t(\omega)=\P_\omega\circ\xi_t^{-1}$ for
$\omega\notin N(t)$. Let $N_2:=\bigcup_{t\in\mathbb Q_T}N(t)$ with $\mathbb Q_T:=(\mathbb Q\cap[0,T])\cup\{0,T\}$, a $\P$--null set. 
	The statement holds with $N_1\cup N_2$ by continuity in time on both sides of $\bar\mu_t(\omega)=\P_\omega\circ\xi_t^{-1}$.
	\medskip
	
{\bf Step 2}. $\P_\omega\circ(X_0,W)^{-1} =\P\circ(X_0,W)^{-1}=\nu\otimes\gamma_{\rm w}$, $\P\as$, where $\gamma_{\rm w}$ denotes the Wiener measure on $\Cc^\bmdim$. 
	In particular, $\P_\omega\circ X_0^{-1}=\nu$ and $\P_\omega\circ W_0^{-1}=\delta_0$.
	Since $(X_0,W)$ and $\overline{\mathcal G}_T$ are independent, for any bounded Borel $h$, $\E^{\P}[h(X_0,W)]$ is constant and a version of $\E^{\P}[h(X_0,W)|\overline{\mathcal G}_T]$.
	Taking $h$ in a countable measure-determining family and since $\omega\longmapsto\P_\omega\circ(X_0,W)^{-1}$ is a kernel, we get $\P_\omega\circ(X_0,W)^{-1} =\P\circ(X_0,W)^{-1}=\nu\otimes\gamma_{\rm w}$, $\P\as$
	\medskip
	
{\bf Step 3}. $\P\as$, $W$ is an $(\F^\circ,\P_\omega)$--Brownian motion and $\E^{\P_\omega}\big[\sup_t\|X_t\|^p
+\int_0^T\|u_t\|^p\d t \big]<+ \infty.$
	Note that the second statement follows from the tower property.
	Let $0\le s<r\le T$ be rationals, $\phi:\R^\bmdim\longrightarrow \R$, $\varphi:\bar\Omega\longrightarrow \R$ be bounded continuous, and define the bounded Borel functions on $\bar\Omega$, 
	$ h(\bar{\mathbf x},\bar{\mathbf u},\bar{\mathbf w})\coloneqq \phi(\bar{\mathbf w}_r-\bar{\mathbf w}_s) \varphi(\bar{\mathbf x}_{s\wedge\cdot},\bar{\mathbf u}_{s\wedge\cdot},\bar{\mathbf w}_{s\wedge\cdot})$,
	$h_2(\bar{\mathbf x},\bar{\mathbf u},\bar{\mathbf w})\coloneqq \varphi(\bar{\mathbf x}_{s\wedge\cdot},\bar{\mathbf u}_{s\wedge\cdot},\bar{\mathbf w}_{s\wedge\cdot})$.
	Note that, since the third component of $\xi_s$ is the path $W$, $h(\xi_s)=\phi(W_r-W_s)\varphi(X_{s\wedge\cdot}, U_{s\wedge\cdot},W_{s\wedge\cdot})$. 
	We first show
	\begin{equation}\label{eq:condindepW}
	\langle h,\bar\mu_s\rangle
	=\E^{\P} [\phi(W_r-W_s)] \langle h_2,\bar\mu_s\rangle, \; \P\as
	\end{equation}
	Since both sides are $\overline{\Gc}_s$--measurable, it suffices to test against $\mathbf 1_G$, $G\in\overline{\Gc}_s$. 
	By \Cref{def.controlrule}\ref{def.controlrule.3}, $\E[\langle h,\bar\mu_s\rangle\mathbf 1_G] =\E[h(\xi_s)\mathbf 1_G] =\E[\phi(W_r-W_s) h_2(\xi_s)\mathbf 1_G]$. 
	Note $h_2(\xi_s)\mathbf 1_G$ is $\Fc_s$--measurable and $W_r-W_s$ is $\P$--independent of $\Fc_s$.
	Thus, $\E[\langle h,\bar\mu_s\rangle\mathbf 1_G] =\E[\phi(W_r-W_s)]  \E[h_2(\xi_s)\mathbf 1_G] =\E[\phi(W_r-W_s)] \E[\langle h_2,\bar\mu_s\rangle\mathbf 1_G]$, and since $\phi$ and $\varphi$ were arbitrary, we obtain \eqref{eq:condindepW}.
	\medskip
	
By {\bf Step 1} at time $s$, for $\omega\in N^c$, $\E^{\P_\omega}[h(\xi_s)]=\langle h,\bar\mu_s(\omega)\rangle$ and $\mathbb E^{\mathbb P_\omega}[h_2(\xi_s)]=\langle h_2,\bar\mu_s(\omega)\rangle$, so \eqref{eq:condindepW} gives that $\P\as$,
	\[
	\E^{\P_\omega}[\phi(W_r-W_s) \varphi(X_{s\wedge\cdot},U_{s\wedge\cdot},W_{s\wedge\cdot})] 
	=\E^{\P}\big[\phi(W_r-W_s)\big] \E^{\P_\omega}\big[\varphi(X_{s\wedge\cdot},U_{s\wedge\cdot},W_{s\wedge\cdot})\big].
	\]
Taking $\varphi\equiv1$ shows that the law of $W_r-W_s$ under $\P_\omega$ agrees with that of $W_r-W_s$ under $\mathbb P$, i.e., $\Nc(0,(r-s)\mathrm I_d)$.
	Thus, the previous equation shows $W_r-W_s$ is $\P_\omega$--independent of $\Fc^\circ_s$. 
	Choosing countable separating families of $(\phi,\varphi)$ and intersecting the null sets over rational $s<r$, and subsequently extending to arbitrary $s<r$ by continuity of $W$ and dominated convergence, we conclude that $W$ has independent, stationary Gaussian increments with respect to $\mathbb F^\circ$ under $\mathbb P_\omega$, $\P\as$ 
	Together with {\bf Step 2}, the claim follows.\medskip

{\bf Step 4}. Arguing as in \cite[Lemma 4.3$(ii)$]{djete2019mckean}, we may enlarge the space to $\widehat \Omega\coloneqq \Omega\times\Omega^\prime$, $\widehat{\P}\coloneqq \P\otimes\P^\prime$, where $(\Omega^\prime,P^\prime)$ carries a Brownian motion independent of everything to account for degenerate $\sigma_r\sigma_r^\top$, such that, with $\mu_r=\mu_r(\omega)$ as in {\bf Step 1},
\[
X_t=X_0+\int_0^t\big(\alpha_r+b(r,X,\mu_r(\omega))\big)\,\mathrm dr
+\int_0^t\sigma_r\,\mathrm dW_r,\qquad t\in[0,T], \; \P_\omega\otimes\P^\prime\as.
\]

Let $N$ collect all the null sets of {\bf Steps~1--4}.
	It follows that for $\omega\in N^c$, under $\P_\omega\otimes\P^\prime$ the process $X$ solves the controlled SDE driven by the Brownian motion $W$, with control $u$, deterministic flow $\mu(\omega)$, initial law $\nu$, and the
integrability in \Cref{def.controlrule}\ref{def.controlrule.1} holds.
	By It\^o's formula, $M^\phi$ is a continuous $(\F^\circ,\P_\omega\otimes\P^\prime)$--local martingale for every $\phi\in C^2_b(\R^{\xdim_\bmdim})$, hence an $(\F^\circ,\P_\omega)$--local martingale, since $M^\phi$ does
not involve the $\Omega^\prime$--coordinate, and ultimately an $(\F,\P_\omega)$--local martingale, because every set in $\Fc_t$ coincides, up to a $\P_\omega$--null set, with a set in $\Fc^\circ_t$ since $\bar\mu$ is $\P_\omega\as$ constant.  
	We conclude that $\P_\omega\in\Pc_{\rm det}(\nu)$.\medskip

The two equalities
	\(
	\Jc(\P)=\langle \Jc(\P_\cdot^{\overline{\Gc}_T}), \P\rangle
	\),
	and,
	\(
	\P\circ X_T^{-1} =\P\otimes (\P_\cdot^{\overline{\Gc}_T}\circ X_T^{-1})
	\),
	follow from the tower property.
	Let us also note that \(\Psi^\varphi(\nu)
=\inf_{\gamma \in \Gamma_{\rm det}(\nu)}J^\varphi(\gamma)\) follows from these two equalities together with the inclusion $\mathcal P_{\rm det}(\nu)\subseteq\mathcal P(\nu)$.
\end{proof}

\begin{proof}[Proof of {\rm \Cref{thm:duality.gen}}]
	For any $\muin,\mufin\in \Pc_p(\R^m)$, with $\Psi^\varphi$ defined in the statement of \Cref{thm:duality.gen}, let us put
	\begin{equation*}
			D(\muin,\mufin) \coloneqq \sup_{\varphi\in C_{b,p}(\R^m)}\big(  \Psi^\varphi(\muin) -  \langle \varphi, \mufin\rangle   \big).
	\end{equation*}
%	{\color{blue}
%	\begin{equation*}
%		D(\muin,\mufin) := \sup_{\varphi\in C_{b,p}(\R^m)}\bigg(\inf_{\P\in \Pc(\muin)} \E^\P\bigg[\int_0^Tf(s, X, u_s,\mu_s)\diff s +  Q_c\varphi(X_T)\bigg] -  \langle \varphi, \mufin\rangle   \bigg)
%	\end{equation*}}

Assume $V_c(\muin, \mufin) = + \infty$ for all $\mufin \in \Pc_p(\R^\xdim)$.
	In this case, 
	 since
	 \begin{equation*}
		D(\muin,\mufin) \ge  \inf_{\P\in \Pc(\muin)} \E^\P\bigg[\int_0^Tf(s, X, u_s,\mu_s)\diff s +   Q_c\varphi(X_T)\bigg] - \langle \varphi, \mufin\rangle,
	\end{equation*}
	for all $\varphi\in C_{b,p}(\R^d)$, 
	it follows from the convex dual representations of $\Wc_c(\P\circ X_T^{-1},\mufin)$, namely \eqref{eq.dualitywot}, that $D(\muin, \mufin)\ge \Jc(\P)  = \infty$. 
	In particular, $V_c(\muin, \mufin) = D(\muin, \mufin)$.\medskip

	The rest of the proof of the convex dual representation will be done in two steps:\medskip

	{\bf Step 1.} \emph{Weak duality}.
	Let $\mufin\in \Pc_p(\R^\xdim)$ be such that $V_c(\muin,\mufin)<\infty$.
	Then, for every $\varepsilon>0$ there is $\varphi \in C_{b,p}(\R^m)$ such that
	\begin{align*}
		D(\muin,\mufin) &\leq \E^{\P}\bigg[\int_0^Tf(s, X, u_s, \mu_s)\diff s + Q_c\varphi(X_T) \bigg]  + \varepsilon 
	\end{align*}
	for any $\P \in \Pc(\muin,\mufin)$.
	Thus, we have 
	\begin{align*}
		D(\muin,\mufin) &\le \E^\P\bigg[\int_0^Tf(s, X, u_s, \mu_s)\diff s \bigg] + \langle Q_c\varphi, \P\circ X^{-1}_T\rangle - \langle \varphi, \mufin\rangle  + \varepsilon\\
		 &\le \Jc(\P) +	\sup_{\varphi \in C_{b,p}(\R^m)}\bigg(\int_{\R^d}Q_c\varphi(x)\P\circ X^{-1}_T(\d x) - \langle \varphi, \mufin\rangle \bigg)   + \varepsilon\\
		&=\Jc(\P)  + \varepsilon .
	\end{align*}
	Thus it follows that $D(\muin,\mufin) \le \Jc(\P) + \varepsilon$.
		Given that $\P \in \Pc(\muin,\mufin)$ and $\varepsilon>0$ were chosen arbitrarily, we deduce that $D(\muin,\mufin) \le V_c(\muin,\mufin)$.

	\medskip

	{\bf Step 2.} \emph{Strong duality}.
	Let again $\mufin\in \Pc_p(\R^\xdim)$ be such that $V_c(\muin,\mufin)<\infty$.
	Thanks to Lemma \ref{Lem:v-LSC} and Proposition \ref{prop.convex.V}, since the topology on $\Pc_p(\R^\xdim)$ is stronger than $\Pc(\R^\xdim)$, the function $\Pc_p(\R^\xdim)\ni \mufin \longmapsto V_c(\muin,\mufin)$ is convex and lower semicontinuous.
	{Denote by $\Mc_p(\R^\xdim)$ the space of Borel signed measures on $\R^\xdim$ with finite $p$ moment}, and equipped with the weak topology.
	%Let $\mufin\in \Pc_2(\R^m)$ be fixed.
	Notice that $\Mc_p(\R^\xdim)$ and $C_{p}(\R^\xdim)$ are in separating duality and thus $\sigma(\Mc_p(\R^\xdim), C_{p}(\R^\xdim))$ is a locally convex space whose dual is $\sigma(C_{p}(\R^\xdim),\Mc_p(\R^\xdim))$, see \cite[Theorem 3.10]{rudin1991functional}.
	We extend the function $V(\muin,\cdot)$ to be $+\infty$ outside $\Pc_p(\R^\xdim)$, and still denote by $V$ this extension.
	It is standard that the function $V(\muin,\cdot)$ remains convex and lower semicontinuous on $\sigma(\Mc_p(\R^\xdim), C_{p}(\R^\xdim))$.
	It is also proper by the assumption above.
	% $f$ is bounded from below and XXXX {\color{blue}we probably need to know that there is at least one $\nu$ such that $\Pc(\mu,\nu)\neq \emptyset$.}
	Thus, by the Fenchel-Moreau theorem, see \cite[Theorem 2.3.3]{zalinescu2002convex}, we have
	\begin{equation*}
		V_c(\muin,\mufin) = \sup_{\varphi \in C_{p}(\R^\xdim)}(\langle \varphi, \mufin \rangle - V^\ast_c(\varphi)), \; \mbox{ where, } V^\ast_c(\varphi) \coloneqq \sup_{\mufin \in \Mc_p(\R^\xdim)}(\langle \varphi, \mufin\rangle - V_c(\muin,\mufin)),\; \varphi \in C_{p}(\R^\xdim),
	\end{equation*}
	is the convex conjugate of $V_c$.
	Since $V^\ast_c$ is monotone, a monotone convergence argument gives
	\begin{equation*}
		V_c(\muin,\mufin) = \sup_{\varphi \in C_{b,p}(\R^\xdim)}(\langle \varphi, \mufin \rangle - V^\ast_c(\varphi)).
	\end{equation*}	
%	({\color{red}on which space? I dont have the book. We proved lsc on $\Pc(\R^\xdim)$ for the weak topology. We need to explain that this gives lsc on $\Mc_p$ for $\sigma(\Mc_p(\R^\xdim), C_{b,p}(\R^\xdim))$ \camilo{C:I think this is clear now, we prove l.s.c. for a weaker topology}. Not yet, we need to justify that the dual space of $\Mc_p$ is $C_{b,p}$}) 
%	{\color{blue}apparently we need to strengthen this to the Wasserstein topology to apply  [3, Theorem 2.3.3] in Camilo's note.}
	% In order to conclude the proof it is enough to show that 
	% \begin{equation}
	% \label{eq:need.to.finish}
	% 	V^*_c(-\varphi) \le - \langle \Phi^\varphi, \muin\rangle + g^*(1)\quad \text{for all}\quad \varphi \in C_{b,p}(\R^d).
	% \end{equation}

	To conclude the proof, we compute the convex conjugate $V^\ast_c$.
	We will show that $V^*(-\varphi) = -\Psi^\varphi(\muin)$ for all $\varphi\in C_{b,p}(\R^\xdim)$ %, which thanks to \Cref{lem:Analytical.measure} 
	which will give the result. 
	We first show $-V^\ast(-\varphi)\geq \Psi^\varphi(\muin)$. 
	By definition of $V^\ast_c$, we have
	\begin{align}\label{eq.rules.det.random}
		-V^\ast_c(-\varphi) &= -\sup_{\eta \in \Pc_p(\R^\xdim)}  \big(-\langle \varphi,\eta\rangle -V_c(\muin,\eta)\big)\notag \\ %&= \inf_{\mu\in \Pc_p(\R^\xdim)}  \big\{\langle \mu,v(\mu)+\mu(\varphi)\big\}\\
		&=\inf_{\eta \in \Pc_p(\R^\xdim)} \bigg(\inf_{\P\in \Pc_c(\muin,\eta)} \E^\P \bigg[ \int_0^Tf(s, X, u_s, \mu_s)\diff s\bigg]+\langle \varphi,\eta\rangle\bigg)\notag \\
		&=\inf_{\eta\in \Pc_p(\R^\xdim)} \inf_{\P\in \Pc_c(\muin,\eta)} \bigg( \E^\P \bigg[ \int_0^Tf(s,X,u_s,\mu_s)\diff s \bigg]+\langle \varphi,\eta \rangle\bigg).
\end{align}
Now, we recall that by definition $\P\in \Pc_c(\muin,\mu)$ means $\P\in \Pc(\muin)$ and $0=\Wc_c(\P_T,\mu) \ge \langle Q_c\varphi,\P\circ X^{-1}_T\rangle -   \langle\varphi,\mu\rangle$ for all $\varphi\in C_{b,p}(\R^\xdim)$. 
%Where the inequality holds since $\varphi\in\Phi_{b,p}$. That is,
Then, we have
\begin{align*}
-V^\ast_c(-\varphi)&\ge \inf_{\mu\in \Pc_p(\R^\xdim)} \inf_{\P\in \Pc_c(\muin,\mu)} \bigg( \E^\P \bigg[ \int_0^T f(s, X, u_s, \mu_s)\diff s\bigg]+\langle Q_c\varphi,\P\circ X^{-1}_T \rangle\bigg)\\
%&\geq \inf_{\mu\in \Pc_p(\R^\xdim)} \inf_{\P\in \Pc_c(\muin,\mu)}  \bigg\{ \E^\P \bigg[ \frac12 \int_0^T \|\alpha_t\|^2\d t\bigg]+\P_T(Q_c \varphi)\bigg\}\\
&= \inf_{\P\in\Pc(\muin)} \E^\P \bigg[ \int_0^T f(s, X, u_s, \mu_s)\diff s + Q_c\varphi(X_T)\bigg] = \Psi^\varphi(\muin)
\end{align*}
where in the first equality, we used the fact that, thanks to \Cref{assump.costs}, we have
\[
\Pc(\muin)=\bigcup_{\mu\in\Pc_p(\R^\xdim)} \Pc_c(\muin,\mu).
\]
We now argue that $-V^\ast(-\varphi)\leq \Psi^\varphi(\muin)$. 
To this end, some definitions and notation are needed.
For any probability measure $P\in \Pc(\Pc(\R^\xdim))$, recall that the intensity $I(P)\in\Pc(\R^\xdim)$ of $P$ is the measure defined as
% for all $f\in C_b(\R^\xdim)$
\[
I(P)(A) :=\int_{\Pc(\R^\xdim)} p(A)P(\d p),\quad \text{for a Borel set }A \subseteq \R^\xdim.
\]
Similarly, the intensity of a measure $\tilde P\in  \Pc(\R^\xdim \times\Pc_p(\R^\xdim))$ is the measure $\hat I(\tilde P)\in \Pc(\R^\xdim \times\R^\xdim)$ defined as
% for all $f\in C_b( \R^\xdim \times \R^\xdim)$
\[%\int_{\R^\xdim \times \R^\xdim} f(x,y)  
\hat I(\tilde P) (A)=\int_{\R^\xdim \times\Pc_p(\R^\xdim)} \int_{ \Pc_p(\R^\xdim)}    \mathbf{1}_A(x,y )p(\d y) \tilde P(\d x,\d p) \quad\text{for a Borel set } A \subseteq \R^\xdim \times \R^\xdim.\]
For $\pi\in \Pc(\R^\xdim \times\Pc(\R^\xdim))$, let ${\rm Proj}_\bullet (\pi)$ denote the projection onto the $\bullet$-coordinate, e.g., ${\rm Proj}_2 (\pi)\in\Pc(\Pc(\R^\xdim))$ is given by ${\rm Proj}_2 (\pi)(\d p)= \pi(\R^\xdim\times \d p )$.
For any $\mu,\nu\in \Pc(\R^\xdim)$ we define
 $$\Lambda(\mu,\nu)\coloneqq \big\{\tilde P \in \Pc(\R^\xdim \times\Pc_p(\R^\xdim)) : {\rm Proj}_1 \tilde P=\mu,\, I({\rm Proj}_2(\tilde P))=\nu\big\}.$$

 \medskip

Let $\eps>0$, and recall that by \cite[Proposition 7.47]{bertsekas1978stochastic} the function $Q_c \varphi$ is lower semianalytic.
Thus, it follows from \cite[Proposition 7.50]{bertsekas1978stochastic} that there is an analytically measurable kernel $(\pi_x)_{x\in \R^\xdim}$, such that for all $x\in\R^\xdim$, $\pi_x \in \Pc_p(\R^\xdim)$ and, 
\begin{align}\label{eq.ms}
\langle \varphi,\pi_x\rangle +c(x,\pi_x)\leq Q_c\varphi(x)+\eps.
\end{align}
Let $\nu\in \Pc_p(\R^\xdim)$ and $\P\in \Pc_c(\muin,\nu)$. 
	The measure $\P\in {\rm Prob}(\Omega)$ %$\P\circ X^{-1}_T$}
	and $(\pi_x)_{x\in \R^\xdim}$ define the measure $\tilde P^{\P}$ in $ \Pc(\R^\xdim \times\Pc_p(\R^\xdim))$
% characterized by {\color{red}rather than saying "characterized by" can we just define
by $\tilde P^{\P} (\d x,\d p) \coloneqq \delta_{\pi^x}(\d p)\P\circ X_T^{-1} (\d x)$.
	%$\tilde P(dx,dp) := \delta_{\pi^x}(dp)\P\circ X^{-1}_T(dx)$.
%\[\int_{\R^\xdim \times\Pc_p(\R^\xdim)} f(x,p)\tilde P(\d x,\d p)=\int_{\R^\xdim \times\Pc_p(\R^\xdim)}  f(x,p ) \delta_{\pi_x} (\d p) \P\circ X^{-1}_T(\d x)= \int_{\R^\xdim } f(x,\pi_x) \P\circ X^{-1}_T(\d x)\]
Its intensity and projections are respectively given by $\hat I(\tilde P^{\P})(\d x, \d y) =\pi_x(\d y)\P\circ X_T^{-1} (\d x)$,
	%$\hat I(\tilde P)(\d x, \d y) :=\pi_x(\d y)\P\circ X^{-1}_T(\d x)$ 
	${\rm Proj}_1 (\hat I(\tilde P^{\P} ) ) (A)=\P\circ X_T^{-1} (A )$,
	and,
\[
{\rm Proj}_2 (\hat I(\tilde P^{\P} ) ) (A)  = \int_{\R^\xdim} \int_{\R^\xdim}  \mathbf{1}_{A}(y)  \hat I(\tilde P^{\P}) (\d x, \d y) =   \int_{\R^\xdim} \int_{\R^\xdim}  \mathbf{1}_A(y)  \pi_{x}( \d y)\P\circ X_T^{-1}(\d x), 
\]
%\[ {\rm Proj}_2 (\hat I(\tilde P) ) (A)  = \int_{\R^\xdim} \int_{\R^\xdim}  \mathbf{1}_{A}(y)  \hat I(\tilde P) (\d x, \d y) =   \int_{\R^\xdim} \int_{\R^\xdim}  \mathbf{1}_A(y)  \pi_{x}( \d y)\P\circ X^{-1}_T(\d x),  \]
for any Borel set $A \subseteq \R^\xdim$.
Moreover, $I({\rm Proj}_2\tilde P^{\P})\in\Pc(\R^\xdim)$ is given, for any Borel set $A \subseteq \R^\xdim$, by
\begin{align*}
I({\rm Proj}_2\tilde P^\P)(A)=\int_{\Pc(\R^\xdim)} p( A) {\rm Proj}_2\tilde P^\P(\d p)&= \int_{\R^\xdim\times \Pc(\R^\xdim)}p( A) \tilde P^\P(\d x, \d p)\\
& = \int_{\R^\xdim\times \Pc(\R^\xdim)}\mathbf{1}_A( x,y) \delta_{\pi_x}(\d p ) \P\circ X^{-1}_T(\d x)={\rm Proj}_2 (\hat I(\tilde P^\P) ) (A).
\end{align*}
Thus, letting $\hat\mu_2(\P):={\rm Proj}_2 (\hat I(\tilde P^\P) )\in\Pc_p(\R^\xdim)$, we have
\begin{equation}
\label{eq:def.Lambda}
 \tilde P^\P \in \Lambda( \P\circ X_T^{-1} , \hat\mu_2(\P) ) \eqqcolon \Lambda_T.
\end{equation}
Back in \eqref{eq.ms}, recall $\P\in \Pc_c(\mu_0,\nu)$ and $\hat\mu_2(\P)\in \Pc_p(\R^\xdim)$, we find that
\begin{align*}
&  \E^\P \bigg[ \int_0^Tf(s,X,u_s,\mu_s)\diff s\bigg]+\langle Q_c\varphi, \P\circ X_T^{-1} \rangle +\eps\\  
& \geq \E^\P \bigg[ \int_0^Tf(s,X,u_s,\mu_s)\diff s\bigg] + \int_{\R^\xdim}\int_{\R^\xdim} \varphi(y)\pi_x(\d y)\P\circ X_T^{-1}(\d x) + \int_{\R^\xdim}\int_{\Pc(\R^\xdim)} c(x,p)\tilde P^\P(\d x,\d p)  \\
 &\ge \E^\P \bigg[  \int_0^T f(s,X,u_s,\mu_s)\diff s\bigg]  + \langle \varphi,\hat \mu_2(\P)\rangle + \langle c,\tilde P^\P\rangle\\ 
 &\ge \E^\P \bigg[  \int_0^T f(s,X,u_s,\mu_s)\diff s\bigg] + \langle \varphi,\hat \mu_2(\P)\rangle+ \Wc_c(\P\circ X_T^{-1}, \hat\mu_2(\P)) ,
\end{align*}
where the latter inequality follows by \eqref{eq:def.Lambda} and \cite[Lemma 2.1]{backhoff2019existence} which asserts that $\Wc_c(\P_T\circ X_T^{-1} , \hat\mu_2(\P )) =\inf_{\tilde P\in \Lambda_T}\langle c, \tilde P\rangle$.\medskip

Thus, since $\nu$ and $\P\in \Pc_c(\muin,\nu)$ were taken arbitrary, we find that
\begin{align*}
\Psi^\varphi(\muin)+\eps & = \inf_{\nu\in\Pc_p(\R^\xdim)} \inf_{\P \in \Pc_c(\muin, \nu ) }    \bigg(\E^\P \bigg[ \int_0^T f(s,X,u_s,\mu_s)\diff  s\bigg] + \langle Q_c\varphi,\P\circ X_T^{-1}  \rangle\bigg) + \eps\\
& \geq   \inf_{\nu\in\Pc_p(\R^\xdim)}  \inf_{\P \in \Pc_c(\muin, \nu ) }  \E^\P \bigg[  \int_0^T f(s,X,u_s,\mu_s)\diff  s\bigg] + \langle \varphi,\hat\mu_2(\P) \rangle + \Wc_c\big( \P\circ X_T^{-1}  ,  \hat\mu_2(\P) \big)  \\
& \geq   \inf_{\nu\in\Pc_p(\R^\xdim)}  \inf_{\P \in \Pc_c(\muin, \nu ) }  \E^\P \bigg[  \int_0^T f(s,X,u_s,\mu_s)\diff s\bigg] + \langle \varphi,\hat\mu_2(\P) \rangle  ,
\end{align*}
where the last step follows since $\Wc_c\geq 0$.
	Finally, since $\hat\mu_2(\P)\in \Pc_p(\R^\xdim)$, the choice $\nu=\hat \mu_2(\P)$ leads to $\Pc_c(\nu_0,\hat \mu_2(\P))$.
	Thus, we deduce that
\begin{equation*}
\Psi^\varphi(\muin)+\eps  \geq \inf_{\mu\in\Pc_p(\R^\xdim)} \inf_{\P \in \Pc_c(\muin, \mu ) }  \bigg( \E^\P \bigg[  \int_0^T f(s,X,u_s,\mu_s)\diff s\bigg] +  \langle \varphi, \mu\rangle\bigg) =-V^\ast_c(-\varphi).%\qedhere
\end{equation*}	

	That is, $V^\ast_c(-\varphi)=-\Psi^{\varphi}(\muin)$.\medskip

	This shows the dual formulation for $V_c(\muin,\mufin)$ with $\Psi^\varphi(\nu)$ defined by \eqref{eq:def.Phi}. 
	That $\Psi^\varphi$ can be expressed as an infimum over $\Pc_{\rm det}(\nu)$ is given by \Cref{lem:disint}.
	To conclude, note that thanks to \Cref{rmk.equiv} and \Cref{lem:disint}, we can deduce from \eqref{eq.rules.det.random} that $V_c^\ast=(V_c^{\rm det})^\ast$.
	Taking conjugates again and using Fenchel-Moreau theorem we have that $V_c=(V_c^\ast)^\ast=(V_c^{\rm det})^{\ast\ast}$ as desired.
	\qedhere
%Let us now show that for any $\varphi\in C_{b,p}(\R^\xdim)$, the function $x\longmapsto \Phi^\varphi(\delta_x)$ is $\muin$-integrable.
%	First, note that by definition of $V_c^\ast(\varphi)$ that $-\Psi^{\varphi}(\muin)=V_c^\ast(-\varphi)\geq -\langle \varphi, \mufin \rangle -V_c(\muin,\mufin)$.
%	That is, $\Psi^{\varphi}(\muin)\leq V_c(\muin,\mufin)+\langle \varphi, \mufin \rangle <+\infty$
%{\color{blue}In particular $\Psi^\varphi(\muin)<\infty$ and therefore by \Cref{lem:Analytical.measure} the function $x\mapsto \Psi^\varphi(\delta_x)$ is $\muin$-integrable and it holds $\Psi^\varphi(\muin) = \langle\Phi^\varphi,\muin\rangle $.}
\end{proof}

%{\color{cyan} The result follows taking infimum on the left hand side of the inequality.
%	}{\color{red}I don't get this.}

\section{Proof of the finite particle approximation}
\label{sec:proof_of_the_finite_particle_approximation}

This section is dedicated to the proof of Theorem \ref{thm:parti.conv.N}.
To ease the presentation, we split the proof into two steps, providing lower and upper bounds for the limit of the finite particle problem.

%\camilo{what is the version of for each $k\in \N$, $\mu^{(k)}_N\longrightarrow \mu_\infty^{(k)}$ for this result?}

\subsection{The lower bound}
\begin{proposition}
\label{prop:lower.bound.particle}
	Under the conditions of {\rm \Cref{thm:parti.conv.N}}, 
	it holds
	\begin{equation*}
		\liminf_{N\to \infty}N^{-1}V^N_e(\muin,\mufin) \ge %\sup_{\varphi\in C_{b}(\R^d)} (\langle \Phi^\varphi(\delta_\cdot), \muin\rangle - \langle \varphi, \mufin\rangle) = 
		V_e(\muin,\mufin).
	\end{equation*}
\end{proposition}
\begin{proof}
%	%\emph{Step 1: Dual representation.}
%	Given that $(X^{1}_0,\dots, X^N_0)$ are i.i.d. under any $\P^{N,\boldsymbol{\alpha}}$, the problem takes the form
%	\begin{align*}
%		V^N(\muin,\mufin ) 
%		&= 
%		\inf\Big\{ 
%		J^N(\boldsymbol{\alpha}) :  
%		\boldsymbol{\alpha}\in \Gamma^N(\muin^{\otimes N}),\;
%		\P^{N,\boldsymbol{\alpha}}\circ (\X_0)^{-1} = \muin^{\otimes N},\, 
%		\P^{N,\boldsymbol{\alpha}}\circ (X^{i,N}_T)^{-1} = \mufin\mbox{ for all } i=1,\dots,N 
%		\Big\}
%	\end{align*}

The proof of the convergence result is based on our convex dual representation.
	The $N$-particle problem with value $V^N_e(\muin,\mufin)$ does not directly fall into the setting of \eqref{eq:problem} because the random variables $(X^{1}_T, \dots, X^{N}_T)$ are not necessarily independent with respect to $\P^{N,\boldsymbol{\alpha}}$.
	However, applying the same strategy as in the proof of Theorem \ref{thm:duality.gen}, (but rather relying on Fenchel--Moreau theorem applied to the function $\Mc_p(\R^\xdim)^N\ni (\mufin^1,\dots,\mufin^N)\longmapsto \tilde V^N_e(\muin, \mufin^1,\dots,\mufin^N)$) with
	\begin{equation*}
		\tilde V^N_e(\muin, \mufin^1,\dots,\mufin^N)
	\coloneqq 
	\inf\Big\{
	J^N_e\big(\P^{N,\balpha}\big):\;  
	\P^{N,\balpha}\in \Gamma^N(\muin^{\otimes N}) ,\;
	 \P^{N\!,\balpha} \circ (X^{i,N}_T)^{-1} = \mufin^i, \; i= 1,\dots, N 
	 % \Wc_c(\P^{N,\balpha} \circ (X^{i,N}_T)^{-1} , \mufin^i) = 0, \; i= 1,\dots, N 
	 \Big\},
	\end{equation*}
	it follows that $\tilde V^N_e(\muin,\mufin^1,\dots,\mufin^N)$ is proper, convex and lower-semicontinuous on $\sigma(\Mc(\R^\xdim)^N, C_{p}(\R^\xdim)^N)$ and thus admits the dual representation
	\begin{equation}
	\label{eq:raw.dual.N}
		\tilde V^N_e(\muin,\mufin^1,\dots,\mufin^N) 
		= 
		\sup_{\bvarphi \in C_{p}(\R^\xdim)^N} 
		\bigg( \Psi^{N,\boldsymbol{\varphi}}(\muin^{\otimes N})  
		- 
		\sum_{i=1}^N\langle \varphi^i,\mufin^i\rangle \bigg)
		= 
		\sup_{\bvarphi \in C_{b,p}(\R^\xdim)^N} 
		\bigg( \Psi^{N,\boldsymbol{\varphi}}(\muin^{\otimes N})  
		- 
		\sum_{i=1}^N\langle \varphi^i,\mufin^i\rangle \bigg),
	\end{equation}
%\camilo{
%\[
%\sup_{\varphi  \in C_b(\R^\xdim)} 
%		\Big( \langle \Phi^{N,\boldsymbol{\varphi}}, \muin^{\otimes N}\rangle 
%		- 
%		 \langle \varphi,\mufin\rangle \Big)
%\]
%}
where $\boldsymbol{\varphi}\coloneqq (\varphi^1,\dots,\varphi^N)$ and 
\begin{align}\label{eq:psi-phi.N}
		\Psi^{N,\boldsymbol{\varphi}}(\muin^{\otimes N})
		\coloneqq 
		\inf_{\P^{N,\boldsymbol{\alpha}} \in \Gamma^N(\muin^{\otimes N})}
		\E^{\P^{N,\boldsymbol{\alpha}}}\bigg[ 
		\sum_{j=1}^N\int_0^T\frac12\| \sigma^\dagger \alpha^j_t\|^2\diff t 
		+ 
		\sum_{j=1}^N  %Q_c
		\varphi^j(X_T^{j,N}) \bigg] .
\end{align}

Using $V^N_e(\muin, \mufin) = \tilde V^N_e(\muin, \mufin,\dots,\mufin)$ gives the dual representation for $V^N_e$.
	With the representation \eqref{eq:raw.dual.N} of $V^N_e$, we can now prove the lower bound as follows:
\medskip

	Let $\varphi \in C_{b,p}(\R^d)$ be arbitrary.
	By \eqref{eq:raw.dual.N}, letting $\boldsymbol{\varphi} = (\varphi,\dots,\varphi) \in C_{b,p}(\R^\xdim)^N$ we have 
	\begin{align*}
		\frac1NV^N_e(\muin,\mufin) 
		&\ge 
		\frac1N\Psi^{N,\boldsymbol{\varphi}}(\muin^{\otimes N}) 
		- 
		\langle\varphi, \mufin\rangle\\ 
		& = \inf_{\P^{N,\boldsymbol{\alpha}} \in \Gamma^N(\muin^{\otimes N})}
		\frac1N\sum_{j=1}^N\E^{\P^{N,\boldsymbol{\alpha}}}\bigg[
		\int_0^T\frac12\|\sigma^\dagger \alpha^j_t\|^2\diff t 
		+ %Q_c
		\varphi(X_T^{j,N}) \bigg] 
		- \langle \varphi, \mufin\rangle.
	\end{align*}

Noticing that $\sigma$ is constant, by \cite[Theorems 1 and 2 and Proposition 1]{djete2020mckean2} the first term on the right-hand side above converges to the value $\Psi^\varphi(\muin)$ given by 
%{\color{red}I don't get the identity below. Did we omit $Hc(\muin|\bar \mu)$ in the statement of \Cref{thm:duality.gen} and \Cref{lem:Analytical.measure}?}
	\begin{equation*}
%		\langle \Phi^\varphi,\muin\rangle
%		=
		\Psi^\varphi(\muin) 
		= 
		\inf_{\P \in \Pc(\muin)}\E^{ \P}\bigg[
		\int_0^T\frac12\| \sigma^\dagger \alpha_t\|^2\diff t
		+ 
		%Q_c
		\varphi(X_T) \bigg]  
		.
	\end{equation*}
	Since $\varphi\in C_b(\R^d)$ was taken arbitrary, recall \Cref{thm:duality.gen}, this shows that
	\begin{equation*}
		\liminf_{N\to \infty}\frac1NV^N_e(\muin,\mufin) \ge \sup_{\varphi\in C_{b,p}(\R^d)} ( \Psi^\varphi (\muin) - \langle \varphi, \mufin\rangle) = V_e(\muin,\mufin),
	\end{equation*}
	which is the claimed bound.
\end{proof}

\subsection{The upper bound}\label{sec.upperbound}

The rest of the proof is to show that we also have
\begin{equation*}
		\limsup_{N\to \infty}N^{-1}V^N_e(\muin,\mufin) \le 
		V_e(\muin,\mufin),
	\end{equation*}
	which is significantly more technical.
	The proof has three main ingredients. 
	First, finite-dimensional duality produces near-optimal particle potentials.
	Second, reversibility and Bakry-Emery estimates yield a uniform averaged $H^1$ bound.
	Third, permutation invariance and concavity allow us to symmetrize the dual potentials, so only their empirical average enters the limsup argument.
	Compactness of this average then identifies an optimal limiting mean-field potential.
	% The proof will rely on functional inequalities, compactness arguments based on Sobolev embedding theorems, and properties of time-reversed diffusions.
	% We will further split this part of the proof into intermediate results.
	We begin by constructing an approximating sequence of potentials.
	Throughout this section, unless otherwise stated, the conditions of {\rm \Cref{thm:parti.conv.N}} are satisfied.

\begin{lemma}\label{lemma.approx.potential}
Let $N\in\N$ and $\varepsilon>0$ be fixed.
	There are functions $\varphi^{1,N},\dots, \varphi^{N,N}\in C_{b,p}(\R^\xdim)$ and a probability measure $\Q^N\in\Gamma^N(\muin^{\otimes N})$ with $\Q^{N}\circ (X^{i,N}_T)^{-1} = \mufin$ for all $i=1,\ldots,N$ such that, letting $\boldsymbol{\varphi}^N = (\varphi^{1,N},\ldots, \varphi^{N,N})$, we have that
\begin{equation}
	\label{eq:maximizing.seq.particle}
		\frac{1}{N}J^N_e (\Q^N) 
		-
		\frac{\varepsilon}{N}
		\le 
		\frac{1}{N}V^{N}_e(\muin,\mufin)
		\le 
		\frac1N\Psi^{N,\boldsymbol{\varphi}^N}( \muin^{\otimes N}) 
		- 
		\frac1N\sum_{i=1}^N\langle \varphi^{i,N}, \mufin\rangle 
		+ 
		\frac{\varepsilon}{N}.
	\end{equation}
Moreover, there exists a measurable function $\psi^N: \R^{mN}\longrightarrow\R$ such that
\begin{align}\label{eq.QN.charac}
		\dfrac{\d   \Q^N}{\d  \P^N} 
		=
		\dfrac{\d    \muin^{\otimes N}}{\d   \mu^N}( \X_0) 
		\exp\bigg(
		- \sum_{i=1}^N %Q_c
		\varphi^{i,N}(X_T^{i,N}) 
		- 
		{\color{black} \psi^N(\X_0)} + \eps
		\bigg)
\end{align}
satisfying,
\begin{equation}
	\label{eq:conv.3}
		\E^{ \P^N}\bigg[  
		\dfrac{\d    \muin^{\otimes N}}{\d   \mu^N}( \X_0)  
		\exp\bigg(
		- 
		\sum_{i =1}^N %Q_c
		\varphi^{i,N}(X_T^{i,N}) 
		- 
		\psi^N( \X_0) 
		+ 
		\eps \bigg)\bigg] = 1.
\end{equation}
\end{lemma}

\begin{proof}
The existence of $\varphi^{1,N},\ldots, \varphi^{N,N}\in C_{b,p}(\R^\xdim)$ and $\Q^N$ as in the statement such that \eqref{eq:maximizing.seq.particle} holds follows by definition of $V^{N}(\muin,\mufin)$ and \eqref{eq:raw.dual.N}, recall \eqref{eq.rewardentropy} and \eqref{eq:psi-phi.N}.\medskip

Let us argue the second part of the statement.
	It follows from \eqref{eq:maximizing.seq.particle} that
%	\begin{equation}
%	\label{eq:maximizing.seq.Nparticle.V}
%		\frac{1}{N}V^{N}(\muin,\mufin)
%		\leq
%		\frac1N\inf_{\P^{N,\boldsymbol{\alpha}}  \in \Gamma^N(\muin^{\otimes N})}
%		\E^{\P^{N,\boldsymbol{\alpha}}}\bigg[ 
%		\sum_{j=1}^N\int_0^T\frac12| \alpha^j_t|^2\d t 
%		+ 
%		\sum_{j=1}^N\varphi^{j,N}(X_T^{j,N}) \bigg] 
%		%+\frac1N \Hc(\muin^{\otimes N}| \mu^{N}) 
%		- 
%		\frac1N\sum_{i=1}^N\langle \varphi^{i,N}, \mufin\rangle 
%		+ 
%		\frac{\varepsilon}{N}.
%	\end{equation}	
%	Re-arranging the above inequality and using the extended functions gives 
		\begin{align}
	\label{eq:conv.1}
		\frac1NJ^{N}_e (\Q^N) 
		- 
		\frac{\varepsilon}{N}&
		\le 
		\frac1N\Psi^{N,\boldsymbol{\varphi}^N}(\muin^{\otimes N}) 
		- 
		\frac1N\sum_{i=1}^N\langle \varphi^{i,N}, \mufin\rangle
		+ 
		\frac{\varepsilon}{N}.
	\end{align}

As in \cite[Lemma 3.5]{tan2013optimal}, and using the fact that $X^1_0,\dots,X^N_0$ are i.i.d., we see that 
	\begin{equation*}
	\langle \Phi^{N,\boldsymbol{\varphi}}, \muin^{\otimes N}\rangle = \Psi^{N,\boldsymbol{\varphi}}(\muin^{\otimes N})
	\mbox{ where for }\bx\coloneqq (x^1,\dots, x^N),\,
		\Phi^{N,\boldsymbol{\varphi}}(\bx)
		\coloneqq 
		\inf_{\P^{N,\boldsymbol{\alpha}}  \in \Gamma^N(\delta_{\bx})}
		\E^{\P^{N,\boldsymbol{\alpha}}}\bigg[ 
		\sum_{j=1}^N\int_0^T\frac12\| \sigma^\dagger  \alpha^j_t\|^2\diff t 
		+ 
		\sum_{j=1}^N %Q_c 
		\varphi^j(X_T^{j,N}) \bigg].
	\end{equation*}

Consequently,
\begin{align*}
		\Psi^{N,\boldsymbol{\varphi}^N}(\muin^{\otimes N}) 
		&= 
		\inf_{\P^{N,\boldsymbol{\alpha}}  \in \Gamma^N(\muin^{\otimes N})}
		\sum_{i=1}^N\E^{ \P^{N,\boldsymbol{\alpha}}}\bigg[ \int_0^T\frac12\|  \sigma^\dagger \alpha^i_t\|^2\diff t 
		+ %Q_c
		\varphi^{i,N}\big (X_T^{i,N}\big ) \bigg]\\
		&=
		-\int_{\R^{\xdim N}}\log \bigg( \E^{\P^N_{\bx}}\bigg[\exp\bigg(- \sum_{i=1}^N %Q_c
		\varphi^{i,N}\big(X_T^{i,N}\big) \bigg) \bigg]  \bigg)\muin^{\otimes N} (\d \bx) 
	\end{align*}
	where the latter equation follows by Girsanov's theorem, see, e.g., \cite[Lemma 4.3]{hernandez2025marginal} for details.
\medskip	

%Moreover, because $\Q^{N}\circ (X^{i,N}_T)^{-1} \overset{c}{=} \mufin$ for all $i=1,\dots,N$, it follows by the dual representation of $\Wc$ that $\langle \varphi^{i,N}, \mufin\rangle \ge \langle Q_c\varphi^{i,N}, \Q^N\circ (X^{i,N}_T)^{-1}\rangle$ for all $i=1,\dots,N$. 
Thus, coming back to \eqref{eq:conv.1}, since $\Q^{N}\circ (\X_0)^{-1} = \muin^{\otimes N}$ and $\Q^{N}\circ (X^{i,N}_T)^{-1} %\overset{c}{=} 
=\mufin$ for all $i=1,\ldots,N$, letting $\Q^N_{\bx}$ denote the r.c.p.d.~of $\Q^N$ with respect to $\Fc_0^N$, and thanks to \cite[Corollary 4.4]{hernandez2025marginal}, we deduce from \eqref{eq:conv.1} that
\begin{align*}
	\int_{\R^{\xdim N}}\log\bigg( \E^{\P^N_{\bx}}\bigg[
		\exp\bigg(
		- \sum_{i=1}^N %Q_c
		\varphi^{i,N}(X_T^{i,N})
		+ \varepsilon
		\bigg)\bigg] \muin^{\otimes N} (\d \bx)
%		&\le  
%            \int_{ \R^{\xdim N} } \bigg( 
%        -\sum_{i=1}^N\langle \varphi^{i,N},\mufin\rangle + \varepsilon
%        - \Hc\big(\Q^N_{\bx}|\P^N_{\bx}\big) \bigg)\muin^{\otimes N} (\d \bx)\\  
		& \le \int_{ \R^{\xdim N} } \bigg( 
		\E^{\Q^N_{\bx}}\bigg[- \sum_{i=1}^N%Q_c
		\varphi^{i,N}(X_T^{i,N}) + \varepsilon \bigg] 
		- \Hc\big(\Q^N_{\bx}|\P^N_{\bx}\big) \bigg)\muin^{\otimes N} (\d \bx)  .
\end{align*}
This shows that, for $\muin^{\otimes N}\ae\, \bx$, $\Q^N_{\bx}$ is optimal for the variational problem
	\begin{equation*}
	%\label{eq:var.probl.particle}
		\sup_{\Q\ll   \P^N_{\bx}}  \bigg(
		\E^{\Q}\bigg[-  \sum_{i=1}^N %Q_c
		\varphi^{i,N}(X_T^{i,N}) 
		+ 
		\eps \bigg] 
		- \Hc\big(\Q| \P^N_{\bx} \big) \bigg)
		=\log\E^{\P^N_{\bx}}\bigg[\exp\bigg(-  \sum_{i=1}^N %Q_c
		\varphi^{i,N}(X_T^{i,N})+ \eps\bigg)\bigg].
	 \end{equation*} 
The existence of $\psi^N: \R^{mN}\longrightarrow\R$ as in the statement follows as in the proof of \cite[Theorem 3.8]{hernandez2025marginal}.
%\footnote{We remind the reader that under any feasible measure the canonical process $\W$ is a standard Brownian motion, i.e., $\W_0=0$, hence it does not appear in \eqref{eq.QN.charac}.}
In fact, it is given by
\begin{equation*}
	\psi^N(x) = \log\E_x^{\P^N}\Big[\exp\Big(-\sum_{i=1}^N \phi^{i,N}(X_T^{i,N})+\epsilon\Big)\Big];
\end{equation*} 
and its measurability follows from the Markov kernel.
\end{proof}

Let us now recall the time reversal mapping $\overleftarrow{\Tc}:\Omega^N\longrightarrow \Omega^N$ given by $\overleftarrow{\Tc}(\omega)\coloneqq \overleftarrow{\omega}$ where $\overleftarrow{\omega}_t\coloneqq \omega_{T-t}$ and, for any $\P\in {\rm Prob}(\Omega^N)$, the time-reversed measure $\overleftarrow{\P} \coloneqq \P\circ \overleftarrow{\Tc}$.\medskip
%	Recall that the measure $\P^N$ is assumed reversible. 
%	That is, $\overleftarrow{\P}^N = \P^N$.\medskip

We also introduce $L^N$, the infinitesimal generator of the diffusion $\X$, which for a regular function $f:\R^{\xdim N} \longrightarrow \R$ is given by
\[
	L^Nf \coloneqq   \Bc \cdot \nabla f + \frac12\mathrm{Tr}(\Sigma\Sigma^\top\nabla^2  f).
\]
	The following result gives a PDE representation of $\psi^N$. 

%{\color{blue}Is $\Delta$ the Laplacian or the Hessian for us?Usually $\nabla^2$ = Hessian}

\begin{proposition}
	 \label{prop:rQN=Palpha}
Let $N\geq 1$ be fixed and $\P^N,\Q^N\in{\rm Prob}(\Omega^N)$ as in {\rm \Cref{lemma.approx.potential}}.
	Then, there are square-integrable and $\F^N$-progressively measurable processes $\boldsymbol{\hat \alpha} = (\hat \alpha^1,\dots, \hat \alpha^N)$ such that $\P^{\hat \balpha,N} = \overleftarrow{\Q}^N$.
	 	The processes $\hat \alpha^i$ satisfy the identity
	 	\begin{equation}
        \label{eq:rep.boldalpha.N}
	 		\hat \alpha^i_t = \sigma \sigma^\top\nabla_{x_i}\log \overrightarrow u^N(t,\X_{t}),\; \d t\otimes \d  \overleftarrow{\Q}^N\ae \text{ on } [0,T]\times\Omega^N
	 	\end{equation}
	 	where the function $\overrightarrow{u}^N$, given by
	 	\begin{equation*}
	 		\overrightarrow u^N(t,\boldsymbol{x}) 
			\coloneqq 
			\E^{\P^N}\big[
			\exp\big( 
			-
			\psi^N(\X_T) 
			+ 
			\eps 
			\big) 
			\big| \X_{t} = \bx \big],
	 	\end{equation*}
satisfies the Hamilton-Jacobi-Bellman equation % reversed sense of time
\begin{align}\label{pde.forwardun}
\partial_t \overrightarrow u^N + L^N \overrightarrow u^N=0, \; (t,\bx)\in [0,T)\times \R^{\xdim N}, 
	\; \overrightarrow u^N(T,\bx)=\exp(-\psi^N(\bx)+\eps).
\end{align}
\end{proposition}
\begin{proof}
    First, note that given \eqref{eq.QN.charac} as in the proof of \cite[Theorem 3.11]{hernandez2025marginal} we deduce that
\begin{align}\label{eq.rQN.charac}
\dfrac{\d \overleftarrow {  \Q}^N}{\d\overleftarrow{  \P}^N} = \dfrac{\d    \Q^{N}_T}{\d   \mu^N}(\overleftarrow \X_T) 
		\exp\Big(
		- \sum_{i=1}^N %Q_c
		\varphi^{i,N}(\overleftarrow X_T^{i,N}) 
		- 
		{\color{black} \psi^N(\overleftarrow \X_0)} + \eps
		\Big)
\end{align}

which, since $\P^N$ is reversible implies, see \cite[Remark 3.12]{hernandez2025marginal}, that
\begin{align}\label{eq.FBpotentials}
    \sum_{i=1}^N %Q_c
    \varphi^{i,N}(x^{i})
	=
	\log\Big( 
	\E^{\P^N}\Big[ 
	\exp\big( 
	-
	\psi^N(\X_0) 
	+ 
	\eps\big)
	\big| 
	 \X_T= \bx
	 \Big]\Big), \;    \Q^{N}_T\ae\, \bx.
\end{align}

From \eqref{eq.rQN.charac} we have that 
\begin{align*}
\dfrac{\d \overleftarrow {  \Q}^N}{\d{  \P}^N } = \dfrac{\d    \overleftarrow\Q^{N}_0}{\d   \mu^N}(  \X_0) 
		\exp\Big(
		- 
		 \psi^N( \X_T)
		- 
		\sum_{i=1}^N
		%Q_c
		\varphi^{i,N}(  X_0^{i,N}) 
		 + \eps
		\Big).
\end{align*}

We now claim that $\overleftarrow{\Q}^N_{\bx}$, the r.c.p.d. of $\overleftarrow{\Q}^N$, is optimal for the variational problem

\begin{equation}
	\label{eq:rvar.probl.particle}
		\sup_{\Q\ll   \P^N_{\bx}}  \big(
		\E^{\Q}\big[
		- 
		\psi^N( \X_T) 
		+ 
		\eps \big] 
		- \Hc\big(\Q| \P^N_{\bx} \big) \big).
\end{equation} 

Indeed, it follows from classical results, for instance \cite[Section 4.9]{follmer2025stochastic}, that \eqref{eq:rvar.probl.particle} 
% has a unique solution of the form
% \begin{align*}
% \dfrac{\d \overleftarrow {  \Q}^N}{\d{  \P}^N} = \dfrac{\d    \overleftarrow\Q^{N}_0}{\d   \mu^N}(  \X_0) 
% 		\exp\Big(
% 		- 
% 		 \psi^N( \X_T)
% 		- 
% 		\hat \varphi^{N}( \X_0) 
% 		 + \eps
% 		\Big)
% \end{align*}
% where $\hat \varphi^{N}$ is uniquely defined by
% \[
% \E^{\P_{\bx}}\Big[ 
% 	\exp\big( - 
% 		 \psi^N( \X_T)
% 		- 
% 		\hat \varphi^{N}( \X_0) 
% 		 + \eps\big)
% 		 \Big]=1, \; \overleftarrow\Q^{N}_0\ae\, \bx
% \]
% but this is precisely \eqref{eq.FBpotentials}. 
% 	Therefore, $\hat \varphi^{N}(\bx)=\sum_{i=1}^N\varphi^{i,N}(x^{i}), \Q^{N}_T\ae\, \bx $ as desired.
% 	This proves the claim.
 has a unique solution $\mathds{Q}^N_{\bx}$ given by
    \begin{align*}
\dfrac{\d  \mathds{Q}^N_{\bx}}{\d{  \P}^N_{\bx}} = %\dfrac{\d    \mathds{Q}^N_0}{\d   \mu^N}(   \boldsymbol{x}) 
        \exp\Big(
        - 
         \psi^N( \X_T)
        - 
        \hat \varphi^{N}( \bx) 
         + \eps
        \Big)
\end{align*}
where $\hat \varphi^{N}$ is uniquely defined by
\[
\E^{\P^N_{\bx}}\Big[ 
    \exp\big( - 
         \psi^N( \X_T)
        - 
        \hat \varphi^{N}( \bx) 
         + \eps\big)
         \Big]=1.%, \;    \mathds{Q}^N_0 \ae\, \bx.
\]
Putting $\mathds{Q}^N(\cdot) = \int_{\R^\xdim }\mathds{Q}^N_{\bx}(\cdot)\overleftarrow{\Q}^N_0(\d\bx)$, it follows by \eqref{eq.rQN.charac} and \eqref{eq.FBpotentials}, that $\mathds{Q}^N = \overleftarrow { \Q}^N$ and that $\hat \varphi^{N}(\bx)=\sum_{i=1}^N%Q_c
\varphi^{i,N}(x^{i}), \Q^{N}_T\ae\, \bx $ as desired.
This proves the claim.

	\medskip
	
It follows from \eqref{eq:rvar.probl.particle}, and Girsanov's theorem, that $\overleftarrow{\Q}^N_{\bx}$ is optimal for 
\[
\sup_{\P^{N,\boldsymbol{\alpha}}  \in \Gamma^N(\delta_{\bx})}
		\E^{\P^{N,\balpha}}\bigg[
		%\exp\bigg( 
		- 
		\psi^N( \X_T) 
		+ 
		\eps
		%\bigg) 
		- 
		\frac12\int_0^T\sum_{i=1}^N\|\sigma^\dagger \alpha^i_t\|^2\diff t
		\bigg] 
		=
		\log\Big( 
		\E^{\P_{\boldsymbol{x}}^N}\big[
		\exp\big( 
		-  
		\psi^N( \X_T)
		+ 
		\eps
		\big)\big] \Big)
		=
		\log\big(
		\overrightarrow u (0,\bx)\big)
\]
where the equality for $\overrightarrow u $ in the statement follows from \cite[Section 5]{leonard2011stochastic}.
	The assertion on $\hat\balpha$ follows from classical verification arguments since $\overrightarrow u $ is a classical solution to \eqref{pde.forwardun}.
\end{proof}

In light of the previous result, the convergence of the sequence of potentials will follow from gradient estimates for solutions of backward Kolmogorov equations.
First observe that, without loss of generality, we can (and will) choose the functions $\varphi^{i,N}$ such that $a^{i,N}=\langle %Q_c
\varphi^{i,N}, \mufin\rangle =0$ for all $i=1,\dots,N$.
In fact, if the integrals are not zero, then %(because $Q_c$ is translation invariant in the sense that $Q_c(\phi + a) = Q_c\phi + a$ for $a\in \R$)
 the shifted functions $\hat\varphi^{i,N} := \varphi^{i,N} - a^{i,N}$ and $\hat\psi^N: = \psi^N + \sum_{i=1}^Na^{i,N}$ satisfy \eqref{eq:maximizing.seq.particle}--\eqref{eq:conv.3}.\medskip

%{\color{cyan}In particular we have $Q_c\varphi(\theta,\cdot) = \varphi^{i,N}(x)$ for $(i-1)/N<\theta\le i/N$.}
	%This induces a sequence of Borel functions $\varphi^N$ on $[0,1]\times \R^\xdim$.

%Let $E$ be an Euclidean space, $U\subseteq E$ and $\nu$ be a Borel measure on $U$.
For $p\geq 1$, we recall the Sobolev space $W^{1,p}(\R^\xdim,\nu)$ of weakly differentiable functions $v:\R^\xdim\longrightarrow \R$ such that $v\in \L^p(\R^\xdim,\nu)$ and the weak derivative $\nabla v\in \L^p(\R^\xdim,\nu)$, equipped with the norm 
% In what follows, we will denote by $H^1(K,\mu)$ the Sobolev space of real-valued continuous functions on $K$ admitting one weak derivative ({\color{blue}check clear definition in a book. Can we tak $K$ closed?}).
% 	As usual, we equip $H^1(K, \mu)$ with the norm
\begin{equation*}
		\|v\|_{W^{1,p}(\R^\xdim,\nu)} \coloneqq \|v\|_{\L^p(\R^\xdim,\nu)} + \|\nabla v\|_{\L^p(\R^\xdim,\nu)}.
\end{equation*}
For $p=2$, this is a Hilbert space, and we use the notation $H^1(\R^\xdim,\nu)\coloneqq W^{1,2}(\R^\xdim,\nu)$.

\medskip

%{\color{blue}$W^{1,p}(\R^m,\mufin)$ was not defined.}\\

    Our goal now is to derive uniform $H^1$ bounds for the sequences of average functions induced by $(\varphi^{i,N})_{N\geq 1}$.
	To do so, following %\cite{Bakry-Emery85}, see also 
	\cite[Section 1.4.2]{bakry2014analysis}, we introduce the so-called carr\'e du champ operator as well as the carr\'e du champ it\'er\'e:\footnote{We always use the Hilbert-Schmidt norm on matrices.}
	\begin{align*}
		\Gamma(f,g) &\coloneqq \frac12\Big(L^N(fg) -fL^N(g) - gL^N(f) \Big ), \; \Gamma_2(f,g) \coloneqq \frac12\Big( L^N\Gamma(f,g) - \Gamma(f, L^Ng) - \Gamma(g, L^Nf)\Big),
	\end{align*}
	and $\Gamma(f) \coloneqq \Gamma(f,f)$, $\Gamma_2(f) \coloneqq \Gamma_2(f,f)$.
Notice that
\begin{align*}
		\Gamma(f,g) 
		= 
		\frac12\langle \Sigma^\top  \nabla f, \Sigma^\top\nabla g\rangle, \text{ and, }\;  
		2\Gamma_2(f)
		=
		\frac12  \|\Sigma^\top \nabla^2 f\Sigma  \|^2 {-\langle  \Sigma^\top \nabla f, \Sigma^\top \nabla \Bc \, \nabla f\rangle}.
	\end{align*}
	
%\[
%\boxed{
%2\Gamma_2(f)
%=
%\frac12{\rm Tr} \!\left(
%A\nabla^2f\,A\nabla^2f
%\right)
%-
%\left\langle
%A\nabla f,\,
%(\nabla B)\nabla f
%\right\rangle.
%}
%\]
%
%Equivalently, since $A=\Sigma\Sigma^\top$,
%\[
%\boxed{
%2\Gamma_2(f)
%=
%\frac12
%{\rm Tr} \!\left(
%\Sigma\Sigma^\top\nabla^2f\,
%\Sigma\Sigma^\top\nabla^2f
%\right)
%-
%\left\langle
%\Sigma^\top\nabla f,\,
%\Sigma^\top(\nabla B)\nabla f
%\right\rangle.
%}
%\]
%
%Here
%\[
%\nabla B\in\R^{\xdim N\times \xdim N},
%\]
%\camilo{I am getting some extra $\sigma$ factors in the expression for $\Gamma_2(f)$. Also is this the expression for $\Gamma_2(f)$ or $2\Gamma_2(f)$? I get }. 
%\camilo{ In addition, we will equip the set $[0,1]$ with its uniform measure denoted $\eta$.
%	We denote by $\lambda$ the Lebesgue measure on $\R^\xdim$.}
\begin{proposition}
	\label{pro:H1bound.psiN}
		Let $N\geq 1$, $\varphi^{1,N},\dots, \varphi^{N,N}$ be as in {\rm \Cref{lemma.approx.potential}}.
	We have that
		\begin{equation*}
			\sup_{N\ge1}\frac1N\sum_{i=1}^N\| %Q_c
			\varphi^{i,N}\|^2_{H^1(\R^\xdim,\mufin)}<+\infty.
		\end{equation*}
%		\begin{equation*}
%			\sup_{N\ge1}\frac1N\sum_{i=1}^N\| %Q_c
%			\varphi^{i,N}\|_{H^1(\R^\xdim,\mufin)}<+\infty,   \mbox{ and } \sup_{N\ge1, i\leq N}\|%Q_c
%			\varphi^{i,N}\|_{\L^2(\R^m,\mufin)}<+\infty. 
%			% : =\sup_N \Big(\|\varphi^N(\theta,\cdot)\|_{\L^2(\mufin)} + \|\nabla \varphi^N(\theta,\cdot)\|_{\L^2(\mufin)}\Big)<\infty.
%		\end{equation*}
	\end{proposition}

\begin{proof}
%Let us note that, for $\mu^{N,k}\coloneqq \P^N\circ(X^1,\dots,X^k)^{-1}$ the $k$-particle marginal, by \cite[Theorem 2.6]{chen2025uniform}, $\mu^{N,k}\longrightarrow \bar \mu^{\otimes k}$, weakly, and thus 
%\begin{align}\label{eq.bounded.initial}
%\frac{1}N \Hc(\muin^{\otimes N}|\mu^N) \leq C.
%\end{align}
Let us define
\(
\overrightarrow v^N(t,\bx)\coloneqq\log (\overrightarrow u^N(t,\bx)),
\)
and note that by \eqref{eq.FBpotentials}, \eqref{pde.forwardun} and the reversibility of $\P^N$, we have that
\begin{align}\label{eq.phi.equal.v}
    \sum_{i=1}^N %Q_c
    \varphi^{i,N}(x^{i}) =\overrightarrow v^N(0,\bx), \;    \Q^{N}_T\ae\, \bx.
\end{align}

Moreover, by \eqref{pde.forwardun} and the chain rule, the function $\overrightarrow v^N$ solves the equation
\begin{align}\label{pde.forwardvn}
	 \partial_t \overrightarrow v^N  +L^N \overrightarrow v^N + \Gamma(\overrightarrow v^N)=0,\;
	\overrightarrow v^N(T,\boldsymbol{x}) =-\psi^N(\bx) +\eps.  
\end{align}

Let us now consider the Fisher information
	\[
	I_t^N 
	\coloneqq 
	\E^{\overleftarrow \Q^N}\big[ \|\nabla \log \overrightarrow{u}^N(t, \X_t))\|^2 \big]
	= 
	\E^{\overleftarrow \Q^N}\big[ \|\nabla  \overrightarrow{v}^N(t, \X_t))\|^2 \big].
	\]
Now, we have that there is a constant $C>0$ such that
\begin{align}
	\notag
		C &\ge   \frac1N \big( J^N_e(\Q^N)+ \Hc(\muin^{\otimes N}|\mu^N)\big)= \frac1N \Hc( \Q^N|\P^N)  
		  = \frac1N  \Hc( \overleftarrow  \Q^N|\overleftarrow \P^N)  \notag \\
		 &= \frac1N \big( J^N_e(\overleftarrow \Q^N)+ \Hc(\Q^{N}_T|\mu^N)\big)\notag 
		  \geq  \frac1N \E^{\overleftarrow \Q^N}\bigg[\frac12\int_0^T\sum_{i=1}^N  { \| \sigma \sigma^\top \nabla_{x^i}\log \overrightarrow u^N(t,\X_t)\|^2} \diff t \bigg] \notag \\
		& \geq   \frac1N  \frac{\Lambda}{2}\int_0^T I^N_t \diff t  
		\label{eq:bound.I.vs.H}
\end{align}
where the first inequality follows from \eqref{eq:maximizing.seq.particle} and \eqref{assump.parti.conv.N}, the first equality follows, for instance, from \cite[Lemma 4.3]{hernandez2025marginal}, the second one holds since relative entropy remains constant under measurable maps, the third equality follows from \cite[Proposition 4.2]{cattiaux2023time}, the second inequality follows by the non-negativity of relative entropy and \Cref{prop:rQN=Palpha}, and the last one is due to \Cref{assump.Df}.\medskip

%{
%    \color{cyan}I don't get the justifications. I thought the second inequality follows by \eqref{eq:rep.boldalpha.N} and positivity of the relative entropy, whereas the last inequality follows by \Cref{assump.Df}.
%    Also, I don't really understand how the footnote is used.
%}
%\footnotetext{By \cite[Proposition 4.2]{cattiaux2023time}, letting $\overleftarrow u^N(t,\bx)\coloneqq \overrightarrow u^N(T-t,\bx)$, we have the equality of the forward and backwards evolving expressions, respectively,
%\begin{align*}
%\E^{\overleftarrow \Q}\bigg[
%	\int_0^T 
%	\| \sigma^\top\nabla_{x_i}\log \overrightarrow u^N(t,\X_{t})\|^2
%	\d t
%	\bigg] 
%	=
%	\E^{\Q}\bigg[
%	\int_0^T 
%	\| \sigma^\top\nabla_{x_i}\log \overleftarrow u^N(t,\X_{t})\| ^2
%	\d t
%	\bigg] .
%\end{align*}
%}

Let us set
\[
\overrightarrow{L}^{N,\overleftarrow \Q^N} f
	\coloneqq 
	L^Nf 
	+
	2 \Gamma(\overrightarrow v^N, f)
\] 
and define $\tilde I^N_t := \E^{\overleftarrow \Q^N}[\Gamma(\overrightarrow v^N(t, \X_t))] = \frac12\E^{\overleftarrow \Q^N}[\| \Sigma\Sigma^\top \nabla \overleftarrow v^N(t, \X_t)\|^2]$, which is the Fisher information up to the scaling $\Sigma\Sigma^\top$.
\medskip

By It\^o's formula applied to $F(t,\boldsymbol{x}):=\Gamma \overrightarrow v^N(t,\boldsymbol{x})$, we have
\begin{align}\label{eq.ito.scaledFisher}
		\tilde I^N_t- \tilde I^N_0 
		&= 
		\int_0^t 
		\E^{\overleftarrow \Q^N}\Big[
		\partial_t\Gamma(\overrightarrow v^N)(r,\X_r) 
		+ 
		\overrightarrow L^{N,\overleftarrow \Q^N}\Gamma(\overrightarrow v^N)(r,\X_r)
		 \Big]
		 \diff r\notag \\ 
		& = 
		\int_0^t  
		\E^{\overleftarrow \Q^N}
		\Big[
		\partial_t\Gamma(\overrightarrow v^N)(r,\X_r) 
		+ 
		L^N\Gamma(\overrightarrow v^N)(r,\X_r) 
		+ 
		2 \Gamma(\overrightarrow v^N, \Gamma(\overrightarrow v^N))(r,\X_r) 
		\Big] 
		\diff r.  
\end{align}

One can check that the carr\'e du champ satisfies $\partial_t\Gamma(f) = 2\Gamma(f, \partial_tf)$.
	Thus, substituting the PDE for $\overrightarrow v^N$ gives
\begin{align*}
\partial_t\Gamma(\overrightarrow v^N) 
 	=  
 	2\Gamma\big(\overrightarrow v^N, - L^N\overrightarrow v^N - \Gamma(\overrightarrow v^N) \big) 
	= -2\Gamma(\overrightarrow v^N, L^N\overrightarrow v^N) - 2\Gamma(\overrightarrow v^N, \Gamma(\overrightarrow v^N)).
\end{align*} 

Back in \eqref{eq.ito.scaledFisher}, it follows from the definition of the carr\'e du champ it\'er\'e that 
\begin{align*}
\tilde I^N_t- \tilde I^N_0  
	& =  
	\int_0^t  
	\E^{\overleftarrow \Q^N}\Big[
	L^N\Gamma(\overrightarrow v^N)(r,\X_r) 
	- 
	2\Gamma(\overrightarrow v^N, L^N\overrightarrow v^N)(r,\X_r) 
	- 2\Gamma(\overrightarrow v^N, \Gamma(\overrightarrow v^N))(r,\X_r) 
	+ 2 \Gamma(\overrightarrow v^N, \Gamma(\overrightarrow v^N))(r,\X_r)  
	\Big]
	\diff r\\ 
	& = 
	\int_0^t 
	\E^{\overleftarrow \Q^N}\Big[ 
	2\Gamma_2(\overrightarrow v^N)(r,\X_r)
	\Big]
	\diff r.
\end{align*}

Let us now recall that by \Cref{assump.Df}\ref{assump.Df.2}, we have %\camilo{check}
\[
2 \Gamma_2(\overrightarrow v^N) 
	= 
	\frac12 \|\Sigma^\top    \nabla^2\overrightarrow v^N\Sigma \|^2 
	-
	\langle  \Sigma^\t \nabla \overrightarrow v^N ,\Sigma^\top  \nabla \Bc\,  \nabla \overrightarrow v^N\rangle 
	\ge 
	\frac12 \|\Sigma^\top  \nabla^2\overrightarrow v^N \Sigma \|^2  
	+ 
	\rho \|\Sigma^\top \nabla \overrightarrow v^N \|^2\geq 0.
\]
	
Hence, the Fisher information is increasing along the solution $\overrightarrow{v}^N$ of the Kolmogorov equation \eqref{pde.forwardvn}.
	In particular, we have $\tilde I^N_0 \le\tilde I^N_t $ for all $t\in [0,T]$, which, due to \eqref{eq:bound.I.vs.H} and \Cref{assump.Df}\ref{assump.Df.1} yields
	\begin{equation*}
	 \frac1N I^N_0
	 \le  
	 \frac1{N\Lambda} \tilde I^N_0
	 \le   
	 \frac{1}{N \Lambda T}
	 \int_0^T 
	 \tilde I^N_t 
	 \diff t 
	 \le  
	 \frac{\|\Sigma\Sigma^\top\|_\infty^2}{N \Lambda T}
	 \int_0^T  
	 I^N_t 
	 \diff t  
	 \le 
	 \frac{2C\|\Sigma\Sigma^\top\|_\infty^2}{\Lambda^2 T}.%\Hc(\Q^N|\P) .
	\end{equation*}
Thus, $ \frac1N I^N_0$ is bounded. 
	Recalling {\color{black}$\Q^{N}\circ (X^{i,N}_T)^{-1} = \mufin$ for all $i=1,\dots,N$} and \eqref{eq.phi.equal.v},
	it follows that
	% \begin{equation*}
	% 	\sup_{N\ge1, i\leq N}\int_{\R^m}|\nabla \varphi^{i,N}(x)|^2\mufin(\d x) \le \sup_{N\ge 1}\E^{\overleftarrow \Q^N}[|\nabla \overrightarrow v^N(0,\X_0)|^2] \le \sup_{N\ge 1}  \frac{1}{N}I^N_0<+\infty.
	% \end{equation*}
    \begin{align*}
        \sup_{N\ge1}\frac1N\sum_{i=1}^N\int_{\R^m}\|\nabla %Q_c
        \varphi^{i,N}(x)\|^2\mufin(\d x) 
        &= 
        \sup_{N\ge 1}\frac1N\sum_{i=1}^N\E^{ \Q^N}[\|\nabla %Q_c
        \varphi^{i,N}(X^{i,N}_T)\|^2]\\
        &= \sup_{N\ge 1}\frac1N\E^{ \Q^N}[\|\nabla \overrightarrow v^N(0,\X_T)\|^2]\\
        &= \sup_{N\ge 1}\frac1N\E^{\overleftarrow \Q^N}[\|\nabla \overrightarrow v^N(0,\X_0)\|^2] \le \sup_{N\ge 1}  \frac{1}{N}I^N_0<+\infty.
    \end{align*}
%    So which should introduce the average function $\phi^N$ before/in the statement of this proposition and state that $\|\phi^N\|_{H^1(\R^\xdim,\mufin)}<\infty$.}
	{Moreover, by Poincar\'e inequality and the fact that $\varphi^{i,N}$ has $\mufin$-mean zero, there is a constant $C>0$, independent of $N$, such that}
	\begin{equation*}
		\|%Q_c
		 \varphi^{i,N} \|_{\L^2(\R^\xdim,\mufin)} \le C \|\nabla %Q_c
		 \varphi^{i,N} \|_{\L^2(\R^\xdim,\mufin)}.
	\end{equation*}
%	\begin{equation*}
%		\|%Q_c
%		\varphi^{i,N}\|_{\L^2(\R^\xdim,\mufin)} \le C \|\nabla %Q_c
%		\varphi^{i,N}\|_{\L^2(\R^\xdim,\mufin)}.
%	\end{equation*}
%    {
%    \color{blue}why is there $\rho$ above? did we assume $\mufin$ satisfies Poincar\'e inequality with constant $1/\rho$?
%    }
	This shows that $\sup_{N\ge 1}\frac1N\sum_{i=1}^N\|\varphi^{i,N}\|^2_{H^1(\R^\xdim,\mufin)}<+\infty$ as desired.
\end{proof}

We will now use the estimate derived above to obtain converging subsequences of appropriate transformations of the potentials.
In the statement below, $B_R\subseteq \R^\xdim$ denotes the ball centered around the origin of radius $R>0$.
\begin{proposition}\label{prop.conv.phi^N}
%	{\color{yellow}There is $1\le p<2$ and} 
Let $N\geq 1$, $\varphi^{1,N},\dots, \varphi^{N,N}$ be as in {\rm \Cref{lemma.approx.potential}}, $I_i\coloneqq [(i-1)/N, i/N),\, i=1,\dots N$, and let $\eta$ denote the Lebesgue probability measure on $[0,1]$.
		Let
	\begin{equation}\label{eq.phiN.lift}
		\bar\phi^N(x)\coloneqq\frac{1}{N}\sum_{i=1}^N \varphi^{i,N}(x), \; x\in \R^\xdim,
		\mbox{ and, }
		\varphi^N:[0,1]\times \R^\xdim\longrightarrow \R, \, \varphi^N(\theta, x) \coloneqq 
		\begin{cases} 
		\varphi^{i,N}(x) &  (\theta,x) \in I_i \times \R^\xdim\\
		0 &\mbox{otherwise}.
		\end{cases}
	\end{equation}
%	where $I_i\coloneqq [(i-1)/N, i/N),\, i=1,\dots N.$\medskip
	Then we have:
	\begin{enumerate}[label=$(\roman*)$, ref=.$(\roman*)$,wide,  labelindent=0pt]
	\item The sequence $\varphi^N$ is bounded in $\L^2([0,1]; H^1(\R^\xdim, \mufin))$ and along a subsequence, $\varphi^N\longrightarrow\varphi$ weakly in $\L^2([0,1]; H^1(\R^\xdim, \mufin))$;
	\item $\bar\varphi^N\longrightarrow \bar\varphi$ with $\bar\varphi(x) \coloneqq \int_{0}^1\varphi(\theta,x)\eta(\d \theta) $ weakly in $H^1(\R^\xdim,\mufin)$ ;
	\item $\bar\varphi^N\longrightarrow \bar\varphi$ strongly in $\L^2_{\mathrm{loc}}(\R^\xdim, \lambda)$;
	\item $\bar \varphi^N \longrightarrow \bar\varphi$ strongly in $\L^p(\R^\xdim, \mufin)$ for every $1\le p<2$;
	\item $\bar \varphi \in H^1(B_R,\lambda)$ for every $R>0$.
\end{enumerate}
    % There is a subsequence $( \varphi^{N_k})_{k\geq 1}$ of $( \varphi^N)_{N\ge1}$ converging weakly in $\L^2(\eta\otimes\mufin, [0,1]\times \R^m)$ to a function $\varphi$ and $\bar\varphi^{N_k}$ converges weakly in $ \L^2(\R^\xdim,\mufin)$ to $\bar\varphi(x)\coloneqq \int_0^1\varphi(\theta,x)\eta(\d\theta)\in \L^2(\R^\xdim,\mufin)$ and also in $\L^p(\R^\xdim,\mufin)$ for every $1\le p<2$.
    % Moreover $\bar\varphi \in H^1(B_R,\lambda)$ for all $R>0$.
\end{proposition}
\begin{proof}
	$(i)$: Let us first note that
	\[
		\sup_{N\ge1} \|  \bar\varphi^{N}\|_{H^1(\R^\xdim,\mufin)}^2
        \leq 
       \sup_{N\ge1}\frac1N\sum_{i=1}^N\| %Q_c
			\varphi^{i,N}\|^2_{H^1(\R^\xdim,\mufin)},
        \mbox{ and, }
       \sup_{ N\geq 1}\|\varphi^{N}\|_{\L^2([0,1]\times \R^m,\eta\otimes\mufin)}^2
       \le  
        \sup_{N\ge1}\frac1N\sum_{i=1}^N\| %Q_c
			\varphi^{i,N}\|^2_{\L^2(\R^\xdim,\mufin)}.
	\]
	By \Cref{pro:H1bound.psiN}, 
    %\sout{for every $i\in\N$ we have $\sup_{N\in \N}\|\varphi^{i,N}\|_{H^1(\R^\xdim,\mufin )}\le C$ for a constant $C$ that does not depend on $i$.Thus,} 
    it holds $\sup_{N\in \N}\|\bar \varphi^{N}\|^2_{H^1(\R^\xdim,\mufin)}<\infty$.
    In particular, this shows that $\varphi^N$ is bounded in $\L^2([0,1]; H^1(\R^\xdim, \mufin))$ and thus that the sequence admits a converging subsequence in the weak topology.\medskip

    \item
	$(ii)$ and $(iii)$: We claim that $\bar\varphi^N, \nabla\bar\varphi^N  \in \L^2_{\rm loc}(\R^\xdim,\lambda)$.
	In fact, under \Cref{assump.Df}\ref{assump.Df.5} there is a constant $C_R>0$ which does not depend on $N$ such that
	\[
	\int_{B_R}{\color{black}( \|\nabla\bar\varphi^N\|^2 +}  | \bar\varphi^N(x)|^2) \lambda(\d x)=\int_{B_R} {\color{black}(\|\nabla\bar\varphi^N\|^2 +}| \bar\varphi^N(x)|^2) \frac{\d\lambda}{\d \mufin}(x) \mufin(\d x)\leq C_{R} \int_{B_R} ( {\color{black}\|\nabla\bar\varphi^N\|^2 +} | \bar\varphi^N(x)|^2 )  \mufin(\d x)<+\infty.
	\]

	Now, recall that by Rellich-Kondrachov Theorem, see, e.g., \citet[Theorem 5.7.1]{evans1994partial} and the ensuing remark, the space $H^1(B_R, \lambda)$ is compactly embedded into $\L^2(B_R,\lambda)$.
	%, where $B_R$ is the ball centered around the origin and of radius $R$.
	Thus, there is a subsequence $N_k$, depending on $R$, such that the restriction to $B_R$ of $\bar\varphi^{N_k}$ converges in $\L^2(B_R,\lambda)$ to a function $\hat \varphi$.
	Using a standard diagonalization argument yields that we can find a subsequence again denoted $N_k$ such that $\bar\varphi^{N_k}$ converges to $\hat \varphi$ in $\L^2_{\mathrm{loc}}(\R^\xdim,\lambda)$.
	Moreover, since $H^1(B_R,\lambda)$ is a Hilbert space, it is reflexive, and thus boundedness in $H^1(B_R,\lambda)$ of $\bar\varphi^N$ implies that  it admits a subsequence converging in the sense of distributions.
	But since strong convergence implies convergence in the sense of distributions, it follows by uniqueness of the limit that $\hat \varphi \in H^1(B_R,\lambda)$ for all $R>0$.

	% another diagonal argument allows to construct a further subsequence that does not depend on $i$.
	%In fact, for $i=1$, there is a subsequence $N_k^1$ such that $\varphi^{1,N_k^1} \longrightarrow \phi^1$ in $\L^2_{\mathrm{loc}}(\R^\xdim,\mufin)$.
	%From that subsequence, we extract a subsequence $N_k^2$ such that $\varphi^{2,N_k^2} \longrightarrow \phi^2$ in $\L^2_{\mathrm{loc}}(\R^\xdim,\mufin)$.
	%Inductively, we construct a subsequence $ N_k^k$ such that for all $i$, $\varphi^{i,N_k^k} \longrightarrow \varphi^i$ in $\L^2_{\mathrm{loc}}(\R^\xdim,\mufin)$.
	%For simplicity we will put again $N_k = N^k_k$.
	
	\medskip

	On the other hand,
	again thanks to \Cref{pro:H1bound.psiN},
    %\sout{$\sup_{N\geq1, i\leq N}\|\varphi^{i,N}\|_{\L^2(\mufin)}<+\infty$,} 
    we have that $\sup_{ N\geq 1}\|\varphi^{N}\|_{\L^2(\eta\otimes\mufin, [0,1]\times \R^m)}<+\infty$.
	This implies that up to a subsequence, $(\varphi^N)_{N\geq 1}$ converges to a weak limit $%{\color{yellow}\varphi}U
	\varphi \in \L^2([0,1]\times \R^\xdim,\eta\otimes \mufin)$.
	Recall that the integral operator $\Ic:\L^2([0,1]\times \R^\xdim)\longrightarrow \L^2(\R^\xdim)$, $(\Ic f)(x) \coloneqq \int_0^1f(\theta,x) \eta(\d\theta)$ for $f\in \L^2([0,1]\times \R^\xdim)$, is continuous and thus weakly continuous, and satisfies, $\Ic \varphi^N =\bar\varphi^N$.
	Hence, it follows that $\bar\varphi^N$ converges weakly to $\bar\varphi=\int_0^1%{\color{yellow}\varphi}U
	\varphi(\theta,\cdot) \eta(\d\theta)=\Ic \varphi $.
	By uniqueness of the weak limit and Jensen's inequality, $\hat\varphi=\bar\varphi\in \L^2(\R^\xdim,\mufin)$.
	%By uniqueness of the limit (and using that by definition $\phi^N(x) = \int_0^1\varphi^N(\theta,x)  \eta(\d\theta)$) we therefore have $\phi = \int_0^1\varphi(\theta,\cdot)  \eta(\d\theta)$.
	% and thus by Fatou's lemma that $\|\phi^i\|_{\L^2(\R^\xdim, \mufin)}$.
\medskip

$(iv)$:    Let  $p\in [1,2)$.
    We now claim $\bar\varphi^{N_k}$ converges in $\L^p(\R^\xdim,\mufin)$.
	Indeed, note that for every $R>0$, we have
	% there is $q>1$, depending only on $p$, such that
	\begin{align}\label{eq.conv.phiN}
		\|(\bar\varphi^{N} - \bar\varphi)\1{B_R^c}\|^p_{\L^p(\R^\xdim,\mufin)} 
		%\le \sup_{N\geq1, i\leq N}\|(\varphi^{i,N} - \phi)\1{B_R^c}\|^p_{\L^p(\R^\xdim,\mufin)} 
		\le \sup_{N\geq1} \|(\bar\varphi^N - \bar\varphi)\|_{\L^2(\R^\xdim,\mufin)} \mufin(B_R^c)^{(p-1)/p}
		\le \frac{C+\|\bar\varphi\|_{\L^2(\R^\xdim,\mufin)}}{R^{(p-1)/p}} \bigg(\int_{\R^\xdim}  |x| \mufin(\d x)\bigg)^{(p-1)/p} %\camilo{\frac{C}{R^{1/q}} \mufin(\R^m)^{1/q} }
	\end{align}
	where the first inequality holds by the definition of $\bar \varphi^N$ and Hölder's inequality, and the last inequality follows from Markov's inequality.
	Thus, using the triangle inequality, we have that
	\begin{align*}
		\|\bar\varphi^{N_k} - \bar\varphi\|_{\L^p(\R^\xdim,\mufin)}^p &\le \|\bar\varphi^{N_k} - \bar\varphi\|^p_{\L^p(B_R,\mufin)} + \|(\bar\varphi^{N_k} -\bar \varphi)\1{B_R^c}\|^p_{\L^p(\R^\xdim,\mufin)}.
		%&\le  \|\phi^{N_k} - \phi\|^p_{\L^p(B_R,\mufin)} + \frac{C +\|\phi\|_{\L^2(\R^\xdim,\mufin)} }{R^{1/q}}\bigg(\int_{\R^\xdim}  |x| \mufin(\d x)\bigg)^{1/q}  .
	\end{align*}
	The claim follows from \eqref{eq.conv.phiN} and taking the limits subsequently in $N_k$ and $R$.\medskip

	$(v)$: This follows from $(ii)$ and (the proof of) $(iii)$.
%\sout{	The same argument but with $\nabla \phi^N$ instead of $\phi^N$ shows that the convergence holds in $W^{1,p}(\R^\xdim,\mufin)$.}
%
%    \medskip
%
%    {
%    \color{cyan}We now show that for $\eta\as$ $\theta$, the function $u(\theta,\cdot)$ is the weak $c$-transform in the $x$-variable of a Borel measurable function $\varphi:[0,1]\times\R^\xdim\to \R$.
%    Define the $\theta$-wise range
%    \[
%    \mathfrak R_c
%    :=
%    \left(
%    U\in L^p([0,1]\times \R^\xdim):
%    U(\theta,\cdot)\in\mathcal R_c
%    \text{ for }\eta\text{-a.e. }\theta
%    \right).
%    \]
%    We know that $\phi^N\to U$ weakly in $L^2([0,1]\times\R^\xdim,\eta\otimes \mufin)$.
%    Thus, by Mazur's lemma, there is a subsequence in the asymptotic convex hull of $\phi^N$ (again denoted $\phi^N$) which converges strongly in $L^2(\eta\otimes \mufin)$.
%    Since $\Rc_c$ is convex {\color{blue}assumption}, it follows that $\mathfrak R_c$ is convex as well, and thus this subsequence is again in $\mathfrak{R}_c$.
%    Moreover, taking a further subsequence, it follows that $\phi^N$ converges almost surely to $u$.
%    In particular, by Fubini's theorem, for $\eta\ae$ $\theta$, $\varphi^N(\theta,\cdot)$ converges to $U(\theta,\cdot)$ $\mufin\as$.
%    Since {\color{blue}we assumed} $\Rc_c$ closed for $\mufin\as$ convergence, it follows that $U(\theta,\cdot)\in \Rc_c$.
%    Thus, there is $\varphi(\theta,\cdot)$ such that $U(\theta,x) = Q_c\varphi(\theta,x)$ $\mufin\as$
%    }
    % {
    % \color{blue} also argue that $\varphi^N$ converges weakly and (assume $Q_c$ weakly continuous, so) that $Q_c\varphi = U$
    % }
\end{proof}

With this, we are in a position to establish the announced upper bound.
The compactness result above is deliberately stated for the empirical average of the particle potentials.
This is sufficient because the finite-particle dual functional is invariant under permutations and concave in the vector of potentials.
Hence, symmetrization cannot decrease its value.
The limisup proof below will replace the particle-dependent potentials by their average and only then pass to the mean-field limit.
\begin{proposition}
\label{prop:upper.bound.particle}
	Under the conditions of {\rm \Cref{thm:parti.conv.N}}, it holds
	\begin{equation*}
		\limsup_{N\longrightarrow  \infty}\frac1NV^N_e(\muin,\mufin) 
		\le 
		V_e(\muin,\mufin).
	\end{equation*}
%	Moreover, there is $\phi\in H^1(\R^\xdim,\mufin)$ such that 
%	$V_e(\muin,\mufin) = \langle \Phi^\phi, \muin \rangle - \langle \phi,\mufin\rangle $.
\end{proposition}
\begin{proof}
	Starting from \eqref{eq:maximizing.seq.particle} we can write
	\begin{align}
	\notag
		\frac{1}{N}V^{N}_e(\muin,\mufin) 
		&\le  
		\frac1N\Psi^{N,\boldsymbol{\varphi}^N}( \muin^{\otimes N}) 
		- 
		\frac1N\sum_{i=1}^N\langle \varphi^{i,N}, \mufin\rangle 
		+ 
		\frac{\varepsilon}{N}\\
		& \label{eq:maximizing.seq.particle.NEW} = 
		\inf_{\boldsymbol{\alpha} \in \Gamma^N(\muin^{\otimes N})}
		 \frac1N\sum_{i=1}^N\E^{\P^{N,\boldsymbol{\alpha}}}\bigg[ \int_0^T\frac12\|\sigma^\dagger \alpha^i_t\|^2\diff t 
		 +
		 %Q_c
		 \varphi^{i,N}(X_T^{i,N}) \bigg] 
		 - 
		 \frac1N\sum_{i=1}^N\langle \varphi^{i,N}, \mufin\rangle 
		 + \frac{\varepsilon}{N}.
	\end{align}
	The goal of the next steps is to take the limit on the right hand side above.

	\medskip
	
	{\bf Step 1}. \emph{Reduction to bounded, Lipschitz controls}.
	Without loss of generality, we can assume that the infimum in the first term on the right-hand side above is taken over bounded processes $\alpha^i$.
	This can be achieved by approximating square-integrable controls with bounded ones.
	We can further reduce the infimum over Lipschitz, Markovian controls (i.e. $\alpha^{i}_t$ is a Lipschitz function of $\X_t$).
	In fact, let $\Gamma^N_{\mathrm{Lip}}(\muin^{\otimes N})$ be the set of elements of $\Gamma^N(\muin^{\otimes N})$ of the form $\balpha = (\alpha^1,\dots,\alpha^N)$ with $\alpha^i_t = a^i_t(\X_t), t\in [0,T]$, for a Borel-measurable function $a^i:[0,T]\times \R^N\longrightarrow \R$ that is Lipschitz continuous in its second argument, uniformly in $t$.
	Let $\varepsilon>0$ and $\balpha \in \Gamma^N(\muin^{\otimes N})$ bounded be such that
	\begin{align*}
		\frac1N\Psi^{N,\boldsymbol{\varphi}^N}( \muin^{\otimes N}) 
		\ge 
		\frac1N\sum_{i=1}^N\E^{\P^{N,\boldsymbol{\alpha}}}\bigg[ 
		\int_0^T\frac12\|\sigma^\dagger \alpha^i_t\|^2\diff t 
		+ 
		%Q_c
		\varphi^{i,N}(X_T^{i,N}) 
		\bigg] 
		+ \varepsilon.
	\end{align*}
	By the mimicking theorem, \cite[Theorem 3.6]{brunick2013mimicking}, there is a filtered probability space $(\mybar{\Omega}, \mybar{\Fc},\mybar \F ,\mybar{\P})$ that supports an $N$-dimensional Brownian motion $(\mybar{W}^1,\dots, \mybar{W}^N)$ and an adapted process $\mybar{\X} = (\mybar{X}^1,\dots, \mybar{X}^N)$ such that, with $a^i_t(\bx) \coloneqq \E^{\P^{N,\balpha}}[\alpha^i_t|\X_t = \bx], t\in [0,T]$, $\mybar{\P}\circ \mybar{\X}^{-1}_t = \P^{N,\balpha}\circ \X^{-1}_t$ for all $t \in [0,T]$, and
	\begin{equation}
		\label{eq:mimicking.sde}
		\d\mybar{X}_t^i 
		= 
		\big( a^i_t(\mybar{\X}_t)  + b(\mybar{X}_t^{i,N}, m^N_{\overline{\X}_t}) \big) \diff t 
		+ 
		\sigma \diff\mybar{W}_t^i,\; \mybar\P\as 
	\end{equation}
	Note that $a^i$ is a bounded measurable function.
	Let $\hat{\alpha}_t^i\coloneqq a^i_t(\X_t), t\in [0,T]$ and $\hat\balpha := (\hat\alpha^1,\dots,\hat\alpha^N)$.
	Then, because under \Cref{assump.Df}\ref{assump.Df.2} the SDE \eqref{eq:mimicking.sde} has a unique weak solution, we have $\P^{N,\hat{\balpha}}\circ \X^{-1}_t = \overline{\P}\circ \mybar{\X}^{-1}_t, t\in [0,T]$.
	Hence, by Jensen's inequality, we have
	\begin{align}\label{eq.mimick.estimate}
		 \frac1N\Psi^{N,\boldsymbol{\varphi}^N}( \muin^{\otimes N}) 
		&\ge 
		\frac1N\sum_{i=1}^N\E^{\P^{N,\boldsymbol{\alpha}}}\bigg[ 
		\int_0^T\frac12\|\sigma^\dagger a^i_t(\X_t)\|^2\diff t 
		+ 
		%Q_c
		\varphi^{i,N}(X_T^{i,N}) 
		\bigg] 
		+ \varepsilon\notag \\
		& = 
		\frac1N\sum_{i=1}^N\E^{\bar{\P}}\bigg[
		\int_0^T\frac12\|\sigma^\dagger a^i_t(\mybar{\X}_t)\|^2\diff t 
		+ 
		%Q_c
		\varphi^{i,N}(\mybar X^{i,N}_T) 
		 \bigg ]  
		 + \varepsilon\notag \\
		& =
		 \frac1N\sum_{i=1}^N\E^{\P^{N,\hat{\balpha}}}\bigg[
		 \int_0^T\frac12\|\sigma^\dagger \hat{\alpha}^i_t\|^2 \diff t 
		 + 
		 %Q_c
		 \varphi^{i,N}(X^{i,N}_T) 
		  \bigg ]  
		  + \varepsilon.
	\end{align}
	Let us now approximate the function $a^i$ by standard mollifications $a^{i,\delta}$ as $\delta\downarrow 0$ and put
	\begin{equation*}
		\dfrac{\d\P^{N,\delta}}{\d\P^N} \coloneqq  \Ec\bigg(\sum_{i=1}^N\int_0^\cdot \hat \alpha^{i,\delta}_t\diff W^i_t \bigg )_T\mbox{ with } \alpha^{i,\delta}_t \coloneqq a_t^{i,\delta}(\X_t).
	\end{equation*}
	
	Notice that, since $\hat\alpha, \hat\alpha^\delta$ are bounded, as $\delta\downarrow 0$ we have that $\sum_{i=1}^N\int_0^T\hat{\alpha}^{i,\delta}_t\diff W_t^i$ converges to $\sum_{i=1}^N\int_0^T\hat{\alpha}^{i}_t\diff W^i_t$ in any $\L^p(\Omega^N,\P^N)$, and then, $\frac{\d\P^{N,\delta}}{\d\P^N}$ converges to $\frac{\d\P^{N,\hat{\balpha}}}{\d\P^N}$ in $\L^2(\Omega^N,\P^N)$.
	Thus, for any $i=1,\dots,N$, as $\delta \downarrow 0$ we have that $\E^{\P^{N,\delta}}[%Q_c
	\varphi^{i,N}(X^{i,N}_T)] \longrightarrow \E^{\P^{N,\hat{\balpha}}}[%Q_c
	\varphi^{i,N}( X^{i,N}_T)]$, and, the estimate 
		\begin{align*}
		&\bigg|\E^{\P^{N, \delta}}\bigg[\int_0^T\|\sigma^\dagger \hat\alpha_t^{i,\delta}\|^2\diff t \bigg] 
		- 
		\E^{\P^{N, \hat{\balpha}}}\bigg[\int_0^T\|\sigma^\dagger \hat\alpha_t^{i,N}\|^2\diff t \bigg]  \bigg| \\
		&
		 \le 
		\bigg| \Big(\E^{\P^{N, \delta}} - \E^{\P^{N, \hat{\balpha}}}\Big)  
		\bigg[ \int_0^T\frac12\|\sigma^\dagger \hat{\alpha}^{i,N}_t\|^2\diff t \bigg]  \bigg|
		+ 
		\frac12\bigg|\E^{\P^{N,\delta}}\bigg[\int_0^T\|\sigma^\dagger \hat{\alpha}_t^{i,\delta}\|^2 - \|\sigma^\dagger \alpha^{i,N}_t\|^2\diff t \bigg ]\bigg|,
	\end{align*}
	together with the fact that $\frac{\d \P^{N,\delta}}{\d\P^{N}}$ is bounded in $\L^2(\Omega^N,\P^N)$, shows the convergence of costs.
	Consequently, back in \eqref{eq.mimick.estimate}, we have
	\begin{align*}
		\frac1N\Psi^{N,\boldsymbol{\varphi}^N}( \muin^{\otimes N}) 
		&\ge 
		\lim_{\delta\longrightarrow  0} \frac1N\sum_{i=1}^N
		\E^{\P^{N,\delta}}\bigg[ \int_0^T\frac12\|\sigma^\dagger \alpha^{i,\delta}_t\|^2\diff t 
		+ %Q_c
		\varphi^{i,N}(X_T^{i,N}) \bigg] + \varepsilon
	\end{align*}
	and therefore, since $\varepsilon>0$ was taken arbitrary,
	\begin{align*}
		\frac1N\Psi^{N,\boldsymbol{\varphi}^N}( \muin^{\otimes N}) 
		&\geq 
		\inf_{\balpha\in \Gamma^N_{\mathrm{Lip}}(\muin^{\otimes N})} 
		\frac1N\sum_{i=1}^N
		\E^{\P^{N,\balpha}}\bigg[ \int_0^T\frac12\|\sigma^\dagger \alpha^{i,N}_t\|^2\diff t 
		+ %Q_c
		\varphi^{i,N}(X_T^{i,N}) \bigg].
	\end{align*}
	And the above inequality is in fact an equality because the left hand side is clearly smaller than the right hand side.

	\medskip

	{\bf Step 2}. \emph{Terminal cost error.}
	%Let us now extend the functions $\varphi^{i,N}$ to functions on $[0,1]\times \R^\xdim$ as follows:
	%Equip $[0,1]$ with the uniform measure $\eta$ and the Borel-sigma algebra and let $\Theta\coloneqq (\theta^i)_{i\le N}$ be independent and uniformly distributed random variables on $[0,1]$.
	%	{
	%	\color{black}
	Consider the final dimensional dual functional
	\begin{equation*}
		D_N(\varphi^1,\dots, \varphi^N) := \frac1N\Psi^{N, \boldsymbol{\varphi}}(\nu^{\otimes N}_0) - \frac1N\sum_{i=1}^N\langle \varphi^i, \mufin\rangle 
	\end{equation*}
	with $\boldsymbol{\varphi} = (\varphi^1,\dots, \varphi^N)$ and where $\Psi^{N, \boldsymbol{\varphi}}(\nu^{\otimes N}_0)$ defined in \eqref{eq:psi-phi.N} is the value of the unconstrained $N$-particle control problem with terminal cost $\sum_{i=1}^N\varphi^i(X^{i,N}_T)$.
	Thus \eqref{eq:maximizing.seq.particle} can be re-written as
	\begin{equation*}
		\frac1NV^N_e(\muin,\mufin) \le D_N(\varphi^{1,N}, \dots, \varphi^{N,N}) + \frac{\varepsilon}{N}.
	\end{equation*}
	Recall the averaged function 
	$$
		\bar\varphi^N(x) := \frac1N\sum_{i=1}^N\varphi^{i,N}(x).
	$$
	We claim that
	\begin{equation*}
		D_N(\varphi^{1,N}, \dots, \varphi^{N,N}) \le D_N(\bar\varphi^N,\dots, \bar\varphi^N).
	\end{equation*}
	In fact, for a permutation $\pi \in S_N$, the set of permutations of $\{1,\dots,N\}$, define $\varphi^N_\pi := (\varphi^{\pi(1), N},\dots, \varphi^{\pi(N),N})$.
	Because the initial law is $\muin^{\otimes N}$ and interaction is through the empirical measures, the control problem \eqref{eq:psi-phi.N} is invariant under relabeling of particles.
	Thus, $D_N(\varphi^N_\pi) = D_N(\varphi^{1,N}, \dots, \varphi^{N,N})$ for every permutation $\pi$.
	In addition, $D_N$ is concave.
	Indeed, because $\Psi^{N, \boldsymbol{\varphi}}(\muin^{\otimes N})$ is an infimum, it is easily seen to be concave (as a function) of $\boldsymbol{\varphi} = (\varphi^1,\dots, \varphi^N)$ and the term $-\frac1N\sum_{i=1}^N\langle \varphi^i,\mufin\rangle$ is linear.
	Thus, by Jensen's inequality, 
	\begin{equation*}
		D_N\Big(\frac{1}{N!}\sum_{\pi \in S_N}\varphi^N_\pi \Big) \ge \frac{1}{N!}\sum_{\pi \in S_N}D_N(\varphi^N_\pi) = D_N(\varphi^{1,N}, \dots, \varphi^{N,N}).
	\end{equation*}
	But 
	\begin{equation*}
		\frac{1}{N!}\sum_{\pi \in S_N}\varphi^N_\pi = (\bar\varphi^N, \dots, \bar\varphi^N).
	\end{equation*}
	Thus, $D_N(\bar\varphi^N,\dots, \bar\varphi^N) \ge D_N(\varphi^{1,N},\dots, \varphi^{N,N})$ as announced.
	Therefore, % (putting $\boldsymbol{\bar\varphi} = (\bar\varphi^N,\dots, \bar\varphi^N)$)
	\begin{align*}
		\frac1NV^N_e(\muin,\mufin) &\le D_N(\bar\varphi^N, \dots, \bar\varphi^N) + \frac{\varepsilon}{N}\\
				%&= \frac{1}{N}\psi^{N, \bar{\boldsymbol{\varphi}}}(\muin^{\otimes}) - \langle \bar\varphi^N, \mufin\rangle + \frac{\varepsilon}{N}\\
				&\le \inf_{\boldsymbol{\alpha} \in \Gamma^N_{\mathrm{Lip}}(\muin^{\otimes N})}\frac{1}{N}\sum_{i=1}^N\E^{\P^{N,\boldsymbol{\alpha}}}\bigg[\int_0^T\frac12\|\sigma^\dagger \alpha^i_t\|^2dt + \bar\varphi^N(X^{i,N}_T) \bigg ]	- \langle \bar\varphi^N, \mufin\rangle + \frac{\varepsilon}{N}\\
				&\le \inf_{\balpha\in \Gamma^N_{\mathrm{Lip}}(\muin^{\otimes N})} \frac1N\sum_{i=1}^N\E^{\P^{N,\boldsymbol{\alpha}}}\bigg[ \int_0^T\frac12\|\sigma^\dagger \alpha^i_t\|^2\diff t + \bar\varphi(X^{i,N}_T)\bigg] - \langle \bar\varphi^N, \mufin\rangle + \frac{\varepsilon}{N} + e^N.
	\end{align*}
	%	}
	% Recall the functions $(\varphi^N,\varphi):[0,1]\times \R^\xdim \longrightarrow \R\times \R$, $\phi:\R^\xdim\longrightarrow \R$ and the partition $(I_i)_{i\leq I}\subseteq [0,1]$ defined in \Cref{prop.conv.phi^N}.
	% Let $(\theta^i)_{i\leq N}$, $\theta^i\in I_i$ to be specified.
	% % as $\varphi^N(\theta, x)\coloneqq \varphi^{i\camilo{,N}}(x)$ for $\theta \in (\theta^{i-1}, \theta^{i}]$.
	% %The function $\varphi^N$ is jointly measurable and 
	% For each $\balpha \in \Gamma^N_{\mathrm{Lip}}(\muin^{\otimes N})$ we have
	% \begin{align*}
	% \notag
	% 	\frac1N\Psi^{N,\boldsymbol{\varphi}^N}( \muin^{\otimes N}) 
	% 	&= 
	% 	 \inf_{\balpha\in \Gamma^N_{\mathrm{Lip}}(\muin^{\otimes N})} 
	% 	  \frac1N\sum_{i=1}^N\E^{\P^{N,\boldsymbol{\alpha}}}\bigg[
	% 	 \int_0^T\frac12\|\alpha^i_t\|^2\diff t 
	% 	 + %Q_c
	% 	 \varphi^N(\theta^i,X^{i,N}_T)\bigg]\\
	% 	&\leq  
	% 	\inf_{\balpha\in \Gamma^N_{\mathrm{Lip}}(\muin^{\otimes N})}   
	% 	\frac1N\sum_{i=1}^N\E^{\P^{N,\boldsymbol{\alpha}}}\bigg[ 
	% 	\int_0^T\frac12\|\alpha^i_t\|^2\diff t 
	% 	+ %Q_c
	% 	\varphi(\theta^i,X^{i,N}_T)\bigg] + e^N
	% \end{align*}
	Where, for some $\balpha_o\in \Gamma^N_{\mathrm{Lip}}(\muin^{\otimes N})$ we fixed, 
	\[
		e^N \coloneqq  \frac1N\sum_{i=1}^N\E^{\P^{N,\balpha_o}}\big[%Q_c
		\bar h^N(X^{i,N}_T)\big],\;  %Q_c\varphi(\theta^i,X^{i,N}_T)] %= \frac1N\sum_{i=1}^N\E^{\P^{N,\balpha_o}}[Q_c\varphi^N(\theta^i,X^{i,N}_T) - U(\theta^i,X^{i,N}_T)]
 \bar h^{N}(x)\coloneqq  \bar\varphi^{N}(x) - \bar\varphi(x).
\]

Now, since $\balpha_o\in \Gamma^N_{\mathrm{Lip}}(\muin^{\otimes N})$, on %a sufficiently rich filtered 
	the probability space $(\Omega, \Fc,\F, \P)$ supporting Brownian motions $(B^i)_{i\geq 1}$,\footnote{We may take $(\Omega^\infty, \Fc^\infty,\F^\infty, \P^\infty)$ the infinite product space induced by $(\Omega^N, \Fc^N,\F^N, \P^N)$.} the SDE
\begin{equation*}
	\d\Xc^{i,N}_t 
	= 
	\big( b(\Xc^{i,N}_t, m^N_{\Xc_t}) + \alpha^{i,N}_o(t,\Xc^{1,N}_t,\dots,\Xc^{N,N}_t)\big) \diff t 
	+ 
	\sigma  \diff B^i_t, 
	\mbox{ with, } 
	m^N_{\Xc_s}\coloneqq \frac1N\sum_{i=1}^N\delta_{\Xc^{i,N}_s},
\end{equation*}
admits a unique solution.
	Thus, $\P^{N,\balpha_o}\circ(X^{1,N},\dots, X^{N,N})^{-1}=\P\circ(\Xc^{1,N}, \dots, \Xc^{N,N})^{-1}$. \medskip

Let $k>0$. We may then estimate the error term as
\begin{align*}\label{eq.errors}
		e^N 
		=  
		\frac1N\sum_{i=1}^N\E^{\P}\Big[\bar h^N(\Xc^{i,N}_T)\Big]
		= 
		\frac1N\sum_{i=1}^N\Big(
		\E^{\P}\Big[ \bar h^N(\Xc^{i,N}_T)\1{|\Xc^{i,N}_T|> k}
		\Big] 
		+ 
		\E^{\P}\Big[
		\bar h^N(\Xc^{i,N}_T)\1{|\Xc^{i,N}_T|\le k}
		\Big] \Big)
		\eqqcolon 
		e^N_1 + e^N_2.
	\end{align*}
%	The term $e^N_2$ is handled below using Krylov's theorem as we did before (here we do not need to stop).
%	We study each error term separately. 

{\bf Step 3}. \emph{Estimating $e_1^N$.}
	By Cauchy--Schwarz and Markov's inequalities, we have 
	\begin{equation}\label{eq.e1}
		| e_1^N| 
		\leq 
		\frac{1}{Nk }
		\sum_{i=1}^N\E\big[
		|\bar h^{N}( \Xc^{i,N}_T)|^2
		\big]^{1/2}
		\E\big[\|\Xc^{i,N}_T\|^2\big]^{1/2}.
	\end{equation}

We claim that there is a constant $C>0$, independent of $N$ and $k$, such that
	\begin{equation}\label{eq.growth.hN}
   		\|\bar h^N\|_{\L^\infty( B_x, \lambda)}  
		\leq 
		C \e^{\kappa/2 \|x\|^2} 
		,\;  
		x\in  \R^\xdim,
	\end{equation}
	where $\kappa$ is the constant in \Cref{assump.Df}\ref{assump.Df.5}.    
	Indeed, let $x\in\R^\xdim$.
	It follows from \Cref{assump.Df}\ref{assump.Df.5} that
	\begin{align*}
    		\int_{ B_x }\|\nabla  \bar h^N(y)\|^2  \lambda (\d y) 
    		&\le 
    		C_\nu^{-1}\e^{\kappa(1 +\|x\|^2)}
    		\int_{ B_x }\|\nabla  \bar h^N( y)\|^2 \mufin (\d y)\\
    		&\le
    		C_\nu^{-1}\e^{\kappa(1 +\|x\|^2)}\bigg(
    		\frac1N \sum_{i=1}^N
    		\int_{B_x} \|\nabla \varphi^{i,N}(y)\|^2 \, \mufin (\d y) 
		+ 
		\int_{B_x} \int_0^1\|\nabla_{\! x} \varphi(\theta,y)\|^2\, \eta(\d\theta)\mufin(\d y)
		\bigg).
	\end{align*}
	Thus, thanks to Propositions \ref{pro:H1bound.psiN} and \ref{prop.conv.phi^N}, there is $C>0$ independent of $N$ such that $ \|\nabla  \bar h^N \|_{\L^2(B_x, \lambda) } \le C \e^{\kappa(1 +\|x\|^2)/2}$.
    	Arguing similarly for $\bar h^N$ we obtain the bound
	\begin{equation*}
		\| \bar h^N\|_{H^{1}( B_x, \lambda)} \le C \e^{\kappa(1 +\|x\|^2)/2}.
	\end{equation*}
	Consequently, by the Sobolev embedding theorem \cite[Theorem 5.6.3.6]{evans1994partial}, there is $C>0$ independent of $x$ and $N$ such that
	\begin{align}\label{eq.growth.sumhN}
		\|\bar h^N\|_{\L^\infty( B_x, \lambda)} \le C\| \bar h^N\|_{W^{1,2}( B_x,\lambda )} \leq   C\e^{\kappa/2\|x\|^2},
	\end{align}
	where $C$ may have changed from one inequality to another. 
	This proves the claim \eqref{eq.growth.hN}.\medskip

We now show that $(\Xc^{i,N})_{i\leq N, N\geq 1}$ has uniform square-exponential moments of order $\kappa>0$ as in \Cref{assump.Df}\ref{assump.Df.5}. 
	To this end, observe that because the control $\balpha_o$ is bounded, say by $C_A$, it follows from \Cref{assump.Df}\ref{assump.Df.2} that
%	\begin{equation*}
%		|\Xc_t^{i,N}| \le |\Xc^{i,N}_0| + \int_0^tk\Big(1 + |\Xc_{s}^{i,N}| + \Big(\frac{1}{N}\sum_{j=1}^N|\Xc^{j,N}_s|^2\Big)^{1/2}\Big) ds + At + |\sigma||W_t^{i}|, \; t\in [0,T].
%	\end{equation*}
%	Thus taking the square on both sides gives
	\begin{equation*}
		\|\Xc_t^{i,N}\|^2 
		\le 
		4\|\Xc^{i,N}_0\|^2 
		+ 4(\ell_b+C_A)^2T^2 
		+  8T\int_0^t \ell_b\Big( \|\Xc_{s}^{i,N}\|^2 + \frac{1}{N}\sum_{j=1}^N\|\Xc^{j,N}_s\|^2\Big) \d s 
		+ 4\|\sigma B^i\|_T^2, \; t\in [0,T],
	\end{equation*}
	where $\|B^i\|_T\coloneqq \sup_{t\in [0,T]} \|B_t^i\|$.
	Next, taking the average and using Gronwall's inequality, we find that
	\begin{align}\label{eq.estimate.average}
		\frac1N\sum_{i=1}^N\|\Xc_t^{i,N}\|^2 
		\le 
		4\e^{16\ell_b  T^2}
		\Big( \frac1N\sum_{i=1}^N \|\Xc^{i,N}_0\|^2 
		+ 
		(\ell_b+C_A)^2T 
		+ 
		 \frac1N\sum_{i=1}^N \|\sigma B^i\|_T^2 
		 \Big).
	\end{align}
%We now plug this above to get
%the constant in front of the average of BMs term should be (mod $|\sigma|^2$) the same as for the average $X_0$'s
%	\begin{align*}
%		|\Xc_t^{i,N}|^2 & \le 4\Big(|\Xc^{i,N}_0|^2  + 8kT^2 e^{16kT^2} \frac1N\sum_{j=1}^N |\Xc^{j,N}_0|^2  \Big)  + 4(k+A)^2T^2 + 32k e^{16kT^2}  (k+A)^2T^3 +8Tk\int_0^t |\Xc_{s}^{i,N}|^2 ds\\
%		& \quad  + 4|\sigma||^2 (\sup_{t\in [0,T]}  |W_t^{i}|^2 + 8k T^2e^{16kT^2}\frac1N\sum_{i=1}^N \sup_{t\in [0,T]} |W^i_t|^2  )
%	\end{align*}	
	
Replacing this in the previous estimate and applying Gronwall's inequality again, letting $c\coloneqq 4(\ell_b+C_A)^2T^2  \e^{8\ell_b T^2}(1 + 8 \ell_b T \e^{16\ell_b T^2})$, we find that
%	\begin{align*}
%		|\Xc_t^{i,N}|^2 &\le 4e^{8kT^2}\Big(1 + 8kT^2 e^{16kT^2}\Big) |\Xc^{i,N}_0|^2 + 4e^{8kT^2}(k+A)^2T^2 + 32k e^{16kT^2 + 8kT^2}  (k+A)^2T^3  + 4e^{8kT^2}(1 + 2T^2k)|\sigma|^2|W_T^{i}|^2\\
%		& = c_1 |\Xc^{i,N}_0| + c + c_2|W_T^i|^2.
%	\end{align*}
%\camilo{
	\begin{align*}
		\|\Xc_t^{i,N}\|^2 
		\leq 
		c+  
		4\e^{8\ell_b T^2}
		\Big(\|\Xc^{i,N}_0\|^2  + 8\ell_b T^2 \e^{16\ell_b T^2} \frac1N\sum_{j=1}^N \|\Xc^{j,N}_0\|^2  \Big) 
		+ 
		4\e^{8\ell_b T^2}  
		\Big(\|\sigma B^i\|_T^2  + 8\ell_b  T^2\e^{16\ell_b T^2}\frac1N\sum_{j=1}^N \|\sigma B^j\|_T^2 \Big).
	\end{align*}
Since $\X^{i,N}_0$ and $B^i$ are independent, we deduce that for $c_1$ as in \Cref{assump.Df}\ref{assump.Df.5} %$c_1\coloneqq 4e^{8kT^2}(1 + 8kT^2 e^{16kT^2})$
	\begin{align}\label{eq.e1.exp}
		\sup_{N\geq 1, i\leq N}\E\big[\e^{\kappa/2 \|\Xc^{i,N}_T\|^2}\big] 
		\le 
		\e^{\kappa c/2}
		\E\big[\e^{\kappa c_1\|\Xc_0\|^2}\big]
		\E\big[\e^{\kappa  c_1 \|\sigma B^1\|_T^2}\big]
		<+\infty,
	\end{align}	
	where $\Xc_0$ has law $\muin$. 
	Indeed, under \Cref{assump.Df}\ref{assump.Df.4}--\ref{assump.Df}\ref{assump.Df.5} both terms are finite since $\kappa c_1<c_0$ and $ \kappa \lambda_{\rm max}(\sigma\sigma^\top) c_1<1/2T$, see \cite[Corollary 3.2]{ledoux1991probability}.
	%{\color{red}Did you cite it or did I cite it? there is not corollary 3.2 in this book. \camilo{That was me. I just double check and it is correct. Now, I see in the literature folder another PDF with the same name as the intended reference but actually for a different book. The correct PDF is already there}	}
	\medskip
	
	From \eqref{eq.e1}, \eqref{eq.growth.sumhN} and \eqref{eq.e1.exp} we deduce that for any $k>0$ there is $C>0$, independent of $N$ and $k$, such that
	\begin{align}\label{eq.e1N.final}
	| e_1^N| \leq \frac{C}k.
	\end{align}

{\bf Step 4}. \emph{Estimating $e_2^N$.} 
	Let $R>k$ be fixed and 
	\begin{align*}
	\tau^{i,N}_R 
	\coloneqq 
	\inf\big\{t>0: 
	\|\Xc^{i,N}_t\| \vee M^N_t \ge R
	\big \}
	\mbox{ with } 
	M^N_t 
	\coloneqq 
	\Big(\frac1N\sum_{j=1}^N\|\Xc^{j,N}_t\|^2\Big)^{1/2}.
	\end{align*}
	Moreover, to alleviate the notation, let us write the dynamics of $\Xc^{i,N}$ as
	\begin{equation*}
		\diff \Xc^{i,N}_t = \beta^{i,N}_t\diff  t + \sigma \diff  B^i_t, \mbox{ where, } \beta^{i,N}_t \coloneqq b(\Xc^{i,N}_t, m^N_{\Xc_t}) + \alpha^{i,N}_o(t,\Xc^{1,N}_t,\dots,\Xc^{N,N}_t),
	\end{equation*}
	and note that %because $b$ is Lipschitz and of linear growth, 
	$\beta^{i,N}_t$ is bounded uniformly in $i,N$ on $\{t\le \tau^{i,N}_R\}$. 
	Indeed, since $\| \alpha^{i,N}_0\|\leq C_A$, thanks to \Cref{assump.Df}\ref{assump.Df.2} we have 
	\[
	\|\beta^{i,N}_t\| \le \|b(0,\delta_0)\| + 2\ell_bR +A \eqqcolon  K_R,\; \P\as ,\mbox{ on }\{t\le \tau^{i,N}_R\}.
	\] 

Let us now further decompose the error term $e^N_2$ as
	\begin{equation*}
		e^N_2 
		=
		\frac1N\sum_{i=1}^N 
		\E\Big[\bar h^N(\Xc^{i,N}_T)\1{\|\Xc^{i,N}_T\|\le k}\1{\tau^{i,N}_R\le T}\Big] 
		+ 
		\frac1N\sum_{i=1}^N\E\Big[\bar h^N(\Xc^{i,N}_T)\1{\|\Xc^{i,N}_T\|\le k}\1{\tau^{i,N}_R> T}\Big]
		\eqqcolon
		e^N_{2,1} + e^N_{2,2}.
	\end{equation*}
	
Let us estimate $e^N_{2,1}$.
	Since $\{\tau^{i,N}_R\le T\} = \{\sup_{t\le T}\|\Xc^{i,N}_t\|\ge R\}\cup \{\sup_{t\le T}M^N_t\ge R\}$, subsequently using Cauchy-Schwarz inequality, \eqref{eq.growth.sumhN} and \eqref{eq.e1.exp}, Markov's inequality, and \eqref{eq.estimate.average}, it follows that there is a constant $C>0$, independent of $N,k, R$, which may change from line to line, such that
	\begin{align*}
		| e^N_{2,1} |
		&\le 
		\Big(
		\frac1N\sum_{i=1}^N
		\E\big[
		|\bar h^N(\Xc^{i,N}_T)|^2
		\big]\Big)^{1/2}
		\Big(\frac1N\sum_{i=1}^N
		\P(\tau^{N,i}_R\le T)
		\Big)^{1/2}\\
		&\le 
		C\Big(
		\frac1N\sum_{i=1}^N
		\P(\tau^{N,i}_R\le T)
		\Big)^{1/2}\\
		& \le 
		\frac{C}{R}
		\Big(\frac1N\sum_{i=1}^N
		\E\Big[\sup_{t\in [0, T]}\|\Xc^{i,N}_t\|^2\Big]\Big)^{1/2}
		\le \frac{C}{R}.
	\end{align*}
%	where the last inequality follows from \eqref{eq.estimate.average}.
%	%Since $b$ is Lipschitz and $\alpha^{i,N}_0$ bounded, standard SDE estimates show that 
%	%$\frac1N\sum_{i=1}^N\E\big[\sup_{t\le T}|\Xc^{i,N}_t|^2\big]\le C$.
%	Hence, we have
%	\begin{equation*}
%		e^N_{2,1} \le \frac{C}{R^{1/2}}.
%	\end{equation*}
%	
%	\medskip
%	
%{
%	\color{blue}	
	We now estimate $e^N_{2,2}$.
	% is where the localization step is very useful.
	%Let %$a = \sigma\sigma^\top$ and put 
	Let 
	\begin{equation*}
	 	\frac{\d\Q^{\smallerfont{i}\!,\smallerfont{N}\!,\smallerfont{R}}}{\d \P} 
		\coloneqq
		\Ec\bigg(-\int_0^\cdot \vartheta^{\smallerfont{i}\!,\smallerfont{N}\!,\smallerfont{R}}_s\d B^i_s \bigg )_T, 
		%Z^{i,N,R}_T
		\mbox{ where, }
		\vartheta^{\smallerfont{i}\!,\smallerfont{N}\!,\smallerfont{R}}_t
		\coloneqq 
		\sigma^\dagger
		%\sigma^\top a^{-1}
		\beta^{i,N}_t\1{t\le \tau^{i,N}_R}.
	 \end{equation*}
	Since under %%uniform ellipticity of $\sigma$
	\Cref{assump.Df}\ref{assump.Df.1} we have $\|\vartheta^{\smallerfont{i}\!,\smallerfont{N}\!,\smallerfont{R}}_t\|^2\le K^2_R/\Lambda$,
	it follows that $\Q^{\smallerfont{i}\!,\smallerfont{N}\!,\smallerfont{R}}$ is a probability measure equivalent to $\P$ and by Girsanov's theorem 
	$\mybar{B}^i_t
	\coloneqq  
	B^i_t 
	+ 
	\int_0^t\vartheta^{\smallerfont{i}\!,\smallerfont{N}\!,\smallerfont{R}}_s\diff s
	$ 
	is a $\Q^{i,N,R}$-Brownian motion.
	 Moreover, since $\sigma\vartheta^{\smallerfont{i}\!,\smallerfont{N}\!,\smallerfont{R}}_t = \beta^{i,N}_t\1{t\le \tau^{i,N}_R}$, we have that
	 \begin{equation*}
	 	\d\Xc^{i,N}_t = \beta^{i,N}_t\1{t >\tau^{i,N}_R}\diff t + \sigma \diff \mybar{B}^i_t, \; \Q^{\smallerfont{i}\!,\smallerfont{N}\!,\smallerfont{R}}\as
	 \end{equation*}
	 Consequently, 
	 \[
	 X^{i,N}_T = X^i_0 + \sigma \mybar{B}^i_T, \mbox{ on } \{\tau^{i,N}_R>T\}, \mbox{ and, } X^i_0 + \sigma \mybar{B}^i_T \mbox{ has density }\muin \ast \gamma_{\sigma,T},
	 \]
	 where %so that under $\Q^{\smallerfont{i}\!,\smallerfont{N}\!,\smallerfont{R}}$, 
	 %$X^i_0 + \sigma \mybar{B}^i_T$ has density , where 
	 $\gamma_{\sigma,T}$ denotes the density function of a Gaussian random variable with parameters $0$ and $\sigma\sigma^\t T$.
	 Indeed, $X^i_0$ has law $\muin$ and $\sigma \mybar{B}^i_T$ is an independent Gaussian with density $\gamma_{\sigma,T}$.
	 Notice that %Thus, , %, where $\gamma_{\sigma,T}$ is the density of $\sigma \widetilde{B}^i_T$ under $\Q^{i,N}$ (i.e. the Normal distribution with mean zero and variance $aT$).
	%and 
	$\|\muin\ast\gamma_{\sigma,T}\|_\infty\le C_{\sigma,T} \coloneqq \frac{1}{\sqrt{2\pi \Lambda T}}$.\medskip

	 On the other hand, for every $q\ge1$, we have that
	 \begin{align*}
	 	\E\Big[\Big(
		\frac{\d\P}{\d\Q^{\smallerfont{i}\!,\smallerfont{N}\!,\smallerfont{R}}}
		\Big)^q\Big] 
		&= 
		\E\bigg[\Ec\Big(
		q\int_0^\cdot \vartheta^{\smallerfont{i}\!,\smallerfont{N}\!,\smallerfont{R}}_t\diff \mybar{B}^i_t
		\Big)_T
		\exp\Big(
		\frac{q(q-1)}{2}
		\int_0^T\|\vartheta^{\smallerfont{i}\!,\smallerfont{N}\!,\smallerfont{R}}_t\|^2 \diff t
		\Big)\bigg]\\
	 	&\le 
		\E^{\Q^{\smallerfont{i}\!,\smallerfont{N}\!,\smallerfont{R}}}\bigg[
		\exp\Big(
		\frac{q(q-1)}{2}
		\int_0^T\|\vartheta^{\smallerfont{i}\!,\smallerfont{N}\!,\smallerfont{R}}_t\|^2 \diff t
		\Big)\bigg]\\
		&\leq 
		\exp\Big(\frac{q(q-1)}{2\Lambda}K^2_RT\Big) 
		\eqqcolon C_{q,R,T}.
	\end{align*}
	
With this, we may use H\"older's inequality, first with $p = \frac{q}{q-1}$ and then with conjugates $p^\prime,q^\prime$ such that $pp^\prime<2$, to estimate that %$e^N_{2,2}$ can be estimated as
	\begin{align*}
		|e^N_{2,2}| 
		&= 
		\bigg | \frac1N\sum_{i=1}^N
		\E^{\Q^{\smallerfont{i}\!,\smallerfont{N}\!,\smallerfont{R}}}\bigg[
		\frac{\d \P}{\d \Q^{\smallerfont{i}\!,\smallerfont{N}\!,\smallerfont{R}}}
		\bar h^N(\Xc^{i,N}_T) \1{\|\Xc^{i,N}_\smallfont{T} \|\le k}\1{\tau^{i,N}_R>T} 
		\bigg ]\bigg| \\
		&\le 
		C_{q,R,T}^{1/q}\frac1N\sum_{i=1}^N
		\E^{\Q^{\smallerfont{i}\!,\smallerfont{N}\!,\smallerfont{R}}}
		\Big[
		|\bar h^N(X^i_0 + \sigma \mybar{B}^i_T)|^p
		 \1{\|X^i_\smallfont{0} + \sigma \bar{B}^i_\smallfont{T}\|\le k}
		 \Big]^{1/p}\\
		&\le 
		C_{q,R,T}^{1/q}
		\bigg(\int_{\R^\xdim}
		|\bar h^N(x)|^p \1{\|x\|\le k} \muin\ast\gamma_{\sigma,T}(x)\diff x
		\bigg)^{1/p}\\
	 	&\le 
		C_{q,R,T}^{1/q}
		C_{\sigma,T}^{1/p}
		\bigg(\int_{\R^\xdim} 
		|\bar h^N(x)|^p\1{\|x\|\le k} w_T^{-1}(x) 
		\mufin(\d x)
		\bigg)^{1/p}\\
	 	&\le 
		C_{q,R,T}^{1/q}C_{\sigma,T}^{1/p}
		\bigg(\int_{B_k}
		%\1{\|x\|\le k}
		w_T^{-q^\prime}(x)\mufin(\d x)\bigg)^{\frac{1}{pq^\prime}}
		\|\bar h^N\|_{\L^{pp^\prime}(\R,\mufin)}.
	 	%&\le C_{q,R,T}^{1/q}C_{\sigma,T}^{1/p}
	 \end{align*} 
	%with the latter inequality following by H\"older inequality.
	%Letting $p$ and $p^\prime$ be , %then by the envelop condition on $w$ we still have 
	
Since, under \Cref{assump.Df}\ref{assump.Df.5} and $B_k$ is compact we have $ \int_{B_k} w_T^{-q^\prime}(x)\mufin(\d x)<\infty$.
This allows to deduce that there is a constant $C_{R,k}$ independent of $N$, such that 
	\begin{equation}
	\label{eq.e2}
	e^{N}_{2,2}\le C_{R,k}\|\bar h^N\|_{\L^{pp^\prime}(\R,\mufin)}.
	\end{equation}

\medskip

{\bf Step 5.} \emph{Finishing the proof.}
By \eqref{eq:maximizing.seq.particle.NEW} and Steps 3 and 4 of this proof, see \eqref{eq.e1N.final} and \eqref{eq.e2}, we have that for any $k>0$ there are constants $C$, independent of $N, R$ and $k$, and $C_{R,k}$ independent of $N$, such that
\begin{align*}
 	\frac{1}{N}V^{N}_e(\muin,\mufin)
		& \le 
		\inf_{\boldsymbol{\alpha} \in \Gamma^N (\muin^{\otimes N})} 		\frac1N\sum_{i=1}^N
		\E^{\P^{N,\boldsymbol{\alpha}}}\bigg[ 
		\int_0^T\frac12|\alpha^i_t|^2\diff t + \bar\varphi(X^{i,N}_T)
		\bigg] 
		- 
		\langle 
		\bar \varphi^{N}, \mufin \rangle
		+ 
		\frac{\varepsilon}{N} 
		+
		 e^N_{1} 
		 + 
		 e^N_{2,1} 
		 + 
		 e^N_{2,2} \\
		& \le 
		\inf_{\boldsymbol{\alpha} \in \Gamma^N (\muin^{\otimes N})} 		\frac1N\sum_{i=1}^N
		\E^{\P^{N,\boldsymbol{\alpha}}}\bigg[ 
		\int_0^T\frac12|\alpha^i_t|^2\diff t + \bar\varphi(X^{i,N}_T)
		\bigg] 
		- 
		\langle 
		\bar \varphi^{N}, \mufin \rangle
		+ \frac{\varepsilon}{N} + \frac{C}{k} + \frac{C}{R^{1/2}} + C_{R,k}\|\bar h^N\|_{\L^{pp'}(\R,\mufin)}.
%		& \le \inf_{\boldsymbol{\alpha} \in \Gamma^N_{\mathrm{Lip}}(\muin^{\otimes N})} \frac1N\sum_{j=1}^N\E^{\P^{N,\boldsymbol{\alpha}}}\bigg[ \int_0^T\frac12|\alpha^j_t|^2dt + \phi(X^{j,N}_T)  \bigg] + C \|\phi^{N} - \phi\|_{\L^2(\mufin)} + C\|\nabla H_N  -b^\delta \|^{1/8}_{\L^2}  - \langle \phi^{N}, \mufin\rangle + \frac{\varepsilon}{N} 
 \end{align*} 
 Thus, thanks to \Cref{prop.conv.phi^N}, taking the $\limsup$ above as $N$ goes to infinity for fixed $R,k$, and then as $R$ and $k$ go to infinity we find that
 %{\color{red}Issue here: We can show that $\varphi^N$ converges weakly to some function $\psi$ we also know that $Q_c\varphi^N$ converges weakly to some $Q_c\varphi$. But I dont know that $\varphi= \psi$.}
 %{\color{blue}The notation $H_N$ can be misleading here; $H_N$ does not really depend on $N$. I will adjust when we agree on the drift.}
\begin{align*}
	\limsup_{N\longrightarrow \infty}\frac1NV^N_e(\muin,\mufin) &\le \inf_{\alpha\in \Gamma(\muin)}\E^{\P^\alpha }\bigg[\int_0^T\frac12\|\alpha_t\|^2\diff t + %Q_c
	\bar\varphi( X_T)\bigg]  -\langle \bar\varphi,\mufin\rangle 
	% \\
	%& = \inf_{\alpha\in \Gamma(\muin)}\E^{\P^\alpha }\bigg[\int_0^T\frac12|\alpha_t|^2dt - \phi( X_T)\bigg] - \langle \phi,\mufin\rangle 
	%\le V(\muin,\mufin).
\end{align*}
where the convergence of the first term is due to \cite[Theorems 1 and 2 and Proposition 1]{djete2020mckean2} and that of the second and last terms are due to \Cref{prop.conv.phi^N}.
%Since $k$ was arbitrary, letting $k\longrightarrow \infty$ yields
%with $\phi(x) = \int_0^1\varphi(\vartheta, x)\eta(\vartheta)$.
%Thus, we have
That is, we have
\begin{equation}
	\label{eq:prof.conv.last.bound}
	\limsup_{N\longrightarrow \infty}\frac1NV^N_e(\muin,\mufin) \le \Psi^{\bar\varphi}(\muin) - \langle \bar\varphi, \mufin \rangle 
\end{equation}
%with $\phi(x) = \int_0^1\varphi(\theta,x)\eta(\d\theta)$ and
where the function $\Psi^\phi$ was introduced in the statement of Theorem \ref{thm:duality.gen}.
This yields the desired bound and concludes the proof.
\end{proof}

\begin{proof}[Proof of {\rm \Cref{cor:dual.exists.char}}]
	From \eqref{eq:prof.conv.last.bound} in the proof of Theorem \ref{thm:parti.conv.N}, we have
	\begin{align*}
		V_e(\muin,\mufin)\le \liminf_{N\longrightarrow\infty}\frac1NV^N_e(\muin,\mufin)\le \limsup_{N\longrightarrow\infty}\frac1NV^N_e(\muin,\mufin)\le \Psi^{\bar\varphi}(\muin)- \langle \bar\varphi,\mufin\rangle \le V_e(\muin,\mufin).
	\end{align*}
	This shows that $\bar\varphi\in H^1(\R^\xdim, \mufin)$ is optimal for $V_e$.
\end{proof}

\section{Proof of the dual characterization}
We now present here the proofs of \Cref{theorem.verification} and {\rm \Cref{cor.HJB.linearization}}.

\begin{proof}[Proof of {\rm \Cref{theorem.verification}}]
	Let us argue the equivalence first. 
	The proof relies on a known verification argument for McKean-Vlasov optimal control problems.
	By \cite[Proposition 3.7]{djete2019mckean}, a sufficiently regular solution $v$ of the HJB equation satisfifes $v(0, \muin)\le J^\phi(\gamma)$ for every admissible rule $\gamma$, with equality if and only if the control minimizes the integrated Hamiltonian $\d t\otimes \d\P^\gamma$-\as 
	%recall the verification theorem for the McKean-Vlasov control problem with value $\Psi^\varphi(\muin)$, see e.g. \cite[Proposition 3.7]{djete2019mckean}.
	Therein, the result is stated in terms of the Hamiltonian operator $\hat H$, defined in terms of the integrated map $\hat h$ evaluated at a sufficiently smooth function $v(t,\mu)$, as opposed to the integral of the pointwise Hamiltonian, which is the setting in \Cref{theorem.verification}. 
	That is, for $v \in \Cc^{1,2}_{\rm flat}$, \cite[Theorem 3.7]{djete2019mckean} considers
	\[
	\hat H[v] (t,\mu)\coloneqq  \inf_{u\in \L^2(\mu)} \hat h[v] (t,\mu,u), \; \mbox{ where } \hat h[v] (t,\mu,u)\coloneqq  \int_{\R^\xdim} h(t,x,\mu,\nabla \delta v(t,\mu)(x), \nabla^2 \delta v(t,\mu)(x),  u(x) ) \mu( \d x),
	\]
and where the infimum is over $\L^2(\mu)$, the collection of $\mu$-square integrable maps $u:(\R^\xdim,\mathcal{B}(\R^\xdim),\mu) \longrightarrow \Uc$.%{\color{red}so here what is $\L^2(\mu)$? shouldn't that be $\Uc$-valued $\mu$-square-integrable maps? Camilo: Yes}
\medskip

Therefore, let us argue that
\begin{align}\label{eq.ver.int.ham}
\hat H[v] (t,\mu)= \int_{\R^\xdim}  \bigg(  \inf_{{\bf u} \in \Uc }  h(t,x,\mu,\nabla \delta v(t,\mu)(x), \nabla^2 \delta v(t,\mu)(x),  {\bf u} )\bigg)  \mu(\d x).
\end{align}
	But this is precisely the content of \cite[Theorem 14.60]{rockafellar2009variational} which allows to interchange infimum and integration a decomposable function class and a normal integrand.
	Note that $ \L^2(\mu)$ is decomposable relative to $\mu$ (see the discussion preceding \cite[Theorem 14.60]{rockafellar2009variational}) and since $v \in \Cc^{1,2}_{\rm flat}$ and by the growth condition on $f$, it holds $\hat h[v](t,\mu,u)<+\infty$ on $\L^2(\mu)$.\medskip
	
%	{\color{red}(because/if we use $\L^2(\Uc,\mu)$ then we need to justify decomposability right? I think it works because $\Uc$ is closed.)}\\
%	\camilo{I agree with you. But this is not the way the result we are citing is stated. That issue remains even in \cite[Example 14.62]{rockafellar2009variational}, the last display in Detail therein minimizes over $u\in \R^d$.\medskip
	
%	 So we pick the data of \cite[Theorem 14.60]{rockafellar2009variational}, which matches our goal. 
%	 We are working with $\L^2(\Uc,\mu)$, as a space of measurable functions onto $\R^\xdim\times\S^\xdim$, on the left side of \eqref{eq.ver.int.ham}. But we need to match the $\inf$ on the left side.}\\
	
	Let $(t,\mu)\in [0,T]\times \Pc_2(\R^\xdim)$ be fixed, $\R^\xdim\times( \R^\xdim\times\S^\xdim) \ni  (x,{\bf u})\longmapsto h_{t,\mu}(x,\boldsymbol u )\coloneqq h(t,x,\mu,\nabla \delta v(t,\mu)(x), \nabla^2 \delta v(t,\mu)(x), {\bf u} )$ and $\chi_{\Uc}({\bf u})\coloneqq  0 \1{\Uc}({\bf u})+\infty \1{\Uc^c}({\bf u})$ denote the indicator of the set $\Uc$.
	Notice that, $x\longmapsto h_{t,\mu}(x,\boldsymbol u)$ is measurable, and, $\boldsymbol{u}\longmapsto h_{t,\mu}(x,\boldsymbol{u})$ is continuous since $f$ and the other terms defining $h$ are continuous in $\boldsymbol{u}$.
	Thus, thanks to \cite[Example 14.32]{rockafellar2009variational}, the function $\R^\xdim\times( \R^\xdim\times\S^\xdim) \ni  (x,{\bf u})\longmapsto  \tilde h_{t,\mu}(x, {\bf u} )\coloneqq h_{t,\mu}(x,\boldsymbol u ) +\chi_{\Uc}({\bf u})$ is a normal integral,
%	
%	
%{\color{yellow}	We claim that, for $(t,\mu)\in [0,T]\times \Pc_2(\R^\xdim)$ fixed, $\R^\xdim\times\R^\xdim\times \S^\xdim_+\ni  (x,{\bf u})\longmapsto  \widetilde h(t,x,\mu, {\bf u} )\coloneqq h(t,x,\mu,\nabla \delta v(t,\mu)(x), \nabla^2 \delta v(t,\mu)(x), {\bf u} )+\chi_{\Uc}({\bf u})$ is a normal integrand.
% 	Indeed, by \cite[Example 14.29]{rockafellar2009variational} it suffices to note that, on its effective domain, $\bar h(t,x,\cdot)$ is continuous for each $(t,x)\in [0,T]\times \R^\xdim$, and, since $v \in \Cc^{1,2}_{\rm flat}$, $\bar h (\cdot,\cdot,u)$ is clearly measurable for each $u\in \Uc$.
% 	Note also that, thanks to the coercivity assumption on $f$, the right-hand side of the above displayed equation is finite and,} 
%	
	and thanks to \cite[Theorem 14.60]{rockafellar2009variational},
	\begin{align*}
	\hat H[v] (t,\mu)  =  \inf_{u\in \L^2(\mu)}  \int_{\R^\xdim} \tilde h_{t,\mu}(x, u(x) ) \mu( \d x) 
	& =  \int_{\R^\xdim}  \bigg(  \inf_{{\bf u} \in \R^\xdim\times \S^\xdim }  \tilde h_{t,\mu}(x, {\bf u} ) \bigg)  \mu(\d x)\\
	& =  \int_{\R^\xdim}  \bigg(  \inf_{{\bf u} \in \Uc }  h(t,x,\mu,\nabla \delta v(t,\mu)(x), \nabla^2 \delta v(t,\mu)(x),  {\bf u} )\bigg)  \mu(\d x),
	\end{align*}
where the first and last equalities follow by definition of $\tilde h_{t,\mu}$ and the second is \cite[Theorem 14.60]{rockafellar2009variational}.	
	That is, \eqref{eq.ver.int.ham} holds and, moreover, the minimizers are in correspondence,
	\begin{equation}
	\label{eq:meas.selec.optim}
	\boldsymbol{u}^\star(t,\mu) \in \argmin_{u\in \L^2(\mu)} \hat h[v](t,\mu,u) \Longleftrightarrow \boldsymbol{u}^\star(t,x ) \in \argmin_{{\bf u}\in \Uc }  h(t,x,\mu,\nabla \delta v(t,\mu)(x), \nabla^2 \delta v(t,\mu)(x),  {\bf u} ), \mu\ae
	\end{equation}
	{\color{black}In conclusion, if $(i)$ holds, then by \eqref{eq:meas.selec.optim} $\hat H[v](t,\mu^\star_t) = \hat h[v](t,\mu^\star_t,u^\star(t,\mu^\star_t))$ $\mu_t^\star$\as, thus by \cite[Proposition 3.7]{djete2019mckean} it holds $\boldsymbol{u}^\star(t, X^\star_t)$ is optimal.
	Conversely, if $(ii)$ holds, then by \cite[Proposition 3.7]{djete2019mckean}, $\hat H[v](t,\mu^\star_t) = \hat h[v](t,\mu^\star_t,u^\star(t,\mu^\star_t))$ and using \eqref{eq:meas.selec.optim} and the fact that $\mu^\star_t = \P^\star\circ (X_t^\star)^{-1}$, we have $(i)$.}

\medskip

We are left only with the last two statements.
	Let $c\equiv \int_{\R^\xdim} \1{x\neq y} \rho(\d y)$.
	First, let $\varphi^\star$ be a mean field potential, $v^{\varphi^\star\!\!}$~be a solution to \eqref{eq:HJB}, and $\gamma^\star\in \Gamma_{\rm det} (\muin,\mufin)$ be such that $(i)$ holds.
	It follows from the established equivalence that 
	\[ V_c(\muin,\mufin)= \Psi^{\varphi^\star}(\muin) -\langle \varphi^\star,\mufin\rangle=\E^{\P^\star}\bigg[\int_0^Tf(t,X^\star_t,u^\star_t, \mu^\star_t) \diff t\bigg] +\langle \varphi^\star,\P^\star_T\rangle -\langle \varphi^\star,\mufin\rangle=\E^{\P^\star}\bigg[\int_0^Tf(t,X^\star_t,u^\star_t, \mu^\star_t) \diff t\bigg] .	\]
	Conversely, let $\gamma^\star\in \Gamma_{\rm det} (\muin,\mufin)$ be primal optimal, $\varphi^\star$ be such that $(i)$ holds with $v^{\varphi^\star\!\!}$ solution to \eqref{eq:HJB}. Then
	\[ \Psi^{\varphi^\star}(\muin) -\langle \varphi^\star,\mufin\rangle= v^{\varphi^\star\!\!}(\muin) -\langle \varphi^\star,\mufin\rangle = \E^{\P^\star}\bigg[\int_0^Tf(t,X^\star_t,u^\star_t, \mu^\star_t) \diff t\bigg] = V_e(\muin,\mufin).\qedhere
	\]
\end{proof}

\begin{proof}[Proof of {\rm \Cref{cor.HJB.linearization}}]
Since $v\in \Cc^{1,2,{\rm full}}_{\rm flat}$, to ease the notation, we set 
\[
	u(t,\mu)(x)\coloneqq \frac{\delta v}{\delta \mu}(t,\mu)(x) \mbox{ and } w(t,\mu)(y;x)\coloneqq \frac{\delta^2 v}{\delta \mu^2}(t,\mu)(y;x).
\]
Thus, since $v$ is a solution to \eqref{eq:HJB} with datum $\varphi$, we have that 
\begin{equation}\label{eq:HJBBCGL}
\begin{cases}
	\partial_tv(t,\mu)
	+
	\int_{\R^\xdim} \Big(
	\tfrac12 \Delta u(t,\mu)(x)
	-
	\nabla V \ast \mu(x) \cdot \nabla u(t,\mu)(x)
	-
	\tfrac12 \big| \nabla u(t,\mu)(x) \big|^2 
	\Big)\mu(\d x)
	=0, 
	\; (t,\mu) \in [0,T]\times \Pc_2(\R^\xdim) \\
v(T,\mu)=\int\varphi\,\d\mu  \mbox{ for all } \mu \in \Pc_2(\R^\xdim).
\end{cases}
\end{equation}
Moreover, for any smooth function $g$, we introduce $(\mathcal A_tg)(y)\coloneqq \tfrac12\Delta g(y)-\big(\nabla V\ast \mu^\star_t(y)+\nabla  u(t, \mu^\star_t)( y) \big)\cdot\nabla g(y)$. 
	We argue in four steps.\medskip

{\bf Step 1.} We collect some results from the flat calculus. 
	Let $t\in [0,T]$.
	The flat derivative of the map $\mu\longmapsto\int G(t,y,\mu)\,\mu(\d y)$, with $G$ built from $u$,
	$\nabla_yu$, $\Delta_yu$ and $\nabla W\ast\mu$ as in \eqref{eq:HJBBCGL}, exists and, up to the normalizing constant, satisfies
	\begin{equation}\label{eq:productrule}
	\frac{\delta}{\delta\mu}\Big(\int G(t,y,\cdot)\, \mu(\d y)\Big)(x)
	=
	G(t,x,\mu)+\int\frac{\delta G}{\delta\mu}(t,y,\mu)(x)\,\mu(\d y), \; x\in \R^\xdim. 
	%\mbox{ and } \frac{\delta}{\delta\mu}\big(\nabla W*\mu(y)\big)(x)=\nabla W(y-x)-\nabla W*\mu(y).
	\end{equation}
Moreover, 
	\begin{equation}\label{eq:step3}
	\int\mathcal A_t\big(w(t, \mu^\star_t)(\cdot ;x)\big)(y)\,\mu^\star_t(\d y)
	=
	\int\mathcal A_t\big(w(t, \mu^\star_t)(x;\cdot )\big)(y) \,\mu^\star_t(\d y)
	+
	\int\mathcal A_t \big(  u(t, \mu^\star_t)\big) ( y)  \,\mu^\star_t(\d y), \; x\in \R^\xdim. 
	\end{equation}
Indeed, \eqref{eq:productrule} is the product rule for the flat derivative, which holds since $v\in \Cc^{1,2,{\rm full}}_{\rm flat}$.
	Let us now argue \eqref{eq:step3}.
	Differentiating the normalization $\int u(t,\mu)(z)\,\mu(\d z)=0$ in $\mu$ in the direction $x$, we find thanks to \eqref{eq:productrule} that 
	\(
	\int w(t,\mu)  (z;x)\,\mu(\d z)=-u(t,\mu) (x).
	\)
	Let $x,y\in\R^\xdim$, and note that, since $v\in \Cc^{1,2,{\rm full}}_{\rm flat}$, the map $[0,1]^2\ni (\lambda,\lambda^\prime)\overset{g}{\longmapsto} v\big(t,\mu+\lambda(\delta_x-\mu)+\lambda^\prime(\delta_y-\mu)\big)$,
is $C^2$ on $[0,1]^2\cap \{\lambda+\lambda^\prime \le1\}$, with continuous crossed partial derivatives, and thus they agree at $(0,0)$.
	Computing them from the definition of $u$ and $w$ and using the two normalizations, we find that
	\(
	(\partial_{\lambda \lambda^\prime}^2g)(0,0)
	=
	w(y;x)+u(x)
	,\, 
	(\partial_{ \lambda^\prime\lambda}^2g)(0,0)
	=
	w(x;y)+u(y).
	\)
	That is,
	\[
	w(t,\mu) (y;x)
	=
	w(t,\mu) (x;y)+u(t,\mu)  (y)-u(t,\mu)  (x) \mbox{ for all }x,y\in\R^\xdim.
	\]
	{\color{black}Notice that $w$ is not symmetric here because the flat derivative $u$ is defined by the normalization $\int u(t,\mu)(z)\mu(\d z) =0$.
	The preceding identity is the corresponding symmetry relation for the second derivative.}
	With this, \eqref{eq:step3} follows applying $\mathcal A_t$ in the variable $y$ above at $\mu=\mu^\star_t$ and integrating with respect to $\mu^\star_t$.
%: $\mathcal A_t$ has no zeroth-order term, so it annihilates $u(x)$, and the last term being independent of $x$.
\medskip

{\bf Step 2.} Let $U(t,x)\coloneqq u(t,\mu^\star_t)(x)$. We claim that, %letting $\mathcal A_t$ act on $y\longmapsto w(x;y)$,
\begin{equation}\label{eq:step2}
	\partial_t U(t,x)
	=
	\partial_t u(t,\mu^\star_t)(x)
	+
	\int\mathcal A_t\big(w(x;y)\big)\,\mu^\star_t(\d y),
	 \; (t,x)\in [0,T) \times \mathbb R^\xdim.
\end{equation}
Let  $x\in \mathbb R^\xdim$ be fixed. 
	By assumption, the map $(t,\mu)\longmapsto u(t,\mu)(x)$ belongs to $\mathcal C^{1,2}_{\rm flat}$, with flat derivative $w(t,\mu)(x;\cdot)$ in the space variable.
	Recall that by assumption under $\gamma^\star$ the drift is $-\nabla W\ast \mu^\star_t(X_t^{\star})-\nabla u(t, \mu^\star_t)( X_t^{\star}),\, \d t\otimes\d \P^{\gamma^\star}\ae$
	Then, \eqref{eq:step2} holds for ${\rm Leb}\ae\; t\in [0,T]$, thanks to \cite[Proposition~5.102]{carmona2018probabilisticI}, which holds since $v\in \Cc^{1,2,{\rm full}}_{\rm flat}$, along the flow of marginals $\mu^\star_\cdot$ and the constant process $x$. 
	That the identity holds for every $t\in[0,T)$ follows by continuity of both sides, thanks to the regularity of $v$.
	\medskip

{\bf Step 3.} We claim that there is $\varrho_o:[0,T)\longrightarrow \R$ such that, for any $(t,x)\in [0,T) \times \mathbb R^\xdim$,
\begin{equation}\label{eq:step1}
\begin{aligned}
	&\partial_tu(t,\mu^\star_t)(x) 
		+
		\tfrac12\Delta U(t,x)
		-
		\nabla V \ast\mu^\star_t(x)\cdot\nabla U(t,x)
		-
		\tfrac12|\nabla U(t,x)|^2
		+
		\int\mathcal A_t\big(w(y;x)\big)\,\mu^\star_t(\d y)\\
	&\; -\int\nabla V(y-x)\cdot\nabla U(t,y)\,\mu^\star_t(\d y)
		+
		\int \nabla V \ast \mu^\star_t(y)\cdot\nabla U(t,y)\,\mu^\star_t(\d y)
		=
		\varrho_o(t,\mu^\star_t).
\end{aligned}
\end{equation}
%where $\mathcal A_t$ acts on the \emph{first} slot $y$ of $w(y;x)$.
	Let $(t,x)\in [0,T) \times \mathbb R^\xdim$ be fixed.
	For any $(y,\mu) \in \R^\xdim\times \Pc(\R^\xdim)$, let  
	\[
	\Phi(t,\mu)\coloneqq \partial_tv(t,\mu)+\int G(t,y,\mu)\,\mu(\d y),
	\mbox{ and }
	G(t,y,\mu)\coloneqq \tfrac12\Delta_y u(t,\mu)(y) - \nabla V \ast \mu (y) \cdot \nabla_y u(t,\mu)(y)  -\tfrac12 |\nabla_yu(t,\mu)(y) |^2,
%	\mbox{ with } 
%	u(y)=u(t,\mu)(y).
	\] 
	By \eqref{eq:HJBBCGL}, $\Phi\equiv0$ on $[0,T)\times\Pc_2(\R^d)$. 
	Consider the map $[0,1]\ni \lambda\longmapsto\Phi\big(t,(1-\lambda)\mu+\lambda\delta_x\big)$ and note it vanishes on $[0,1]$, with right derivative at $\lambda=0$ equal to $\fder{\Phi}(t,\mu)(x)-\int\fder{\Phi}(t,\mu)(y) \,\mu(\d y)$. 
	Thus,
	\begin{equation}\label{eq:dPhi}
	\fder{\Phi}(t,\mu)(x)=\varrho_o(t,\mu)\; \mbox{for all }x\in\R^\xdim, \; \varrho_o(t,\mu)\coloneqq \int\fder{\Phi}(t,\mu)(y) \mu(\d y) .
	\end{equation}
	Let us find $\fder\Phi$. 
	By the commutation property of $v\in \Cc^{1,2,{\rm full}}_{\rm flat}$, $\fder{(\partial_tv)}=\partial_tu$, and, for $u=u(t,\mu)$, $w=w(t,\mu)$,
	\[
	\frac{\delta G}{\delta\mu}(t,y,\mu)(x)
		=
		\tfrac12\Delta_yw(y;x)
		-
		\big(
		\nabla V\ast \mu(y)
		+
		\nabla_yu(y)
		\big)
		\cdot\nabla_y w(y;x)
		-
		\big(
		\nabla V(y-x)-\nabla V\ast \mu(y)\big)
		\cdot\nabla_yu(y),
	\]
	
	The result follows differentiating in $\Phi\equiv0$, evaluating at $\mu=\mu^\star_t$, where $u(t,\mu^\ast_t)(\cdot)=U(t,\cdot)$, thanks to \eqref{eq:productrule}. 
	%and the second-order operator in $y$ is $\mathcal A_t$, \eqref{eq:dPhi} becomes
	\medskip

{\bf Step 4.} We prove \eqref{eq:BCGL17}.
	Replace \eqref{eq:step3} into \eqref{eq:step1}, and use \eqref{eq:step2} to replace $\partial_tu(t,\mu^\star_t)(x)+\int\mathcal A_t(w(x;y))\,\mu^\star_t(\d y)$ by $\partial_tU(t,x)$. 
	Notice the $w$-terms cancel and, since $W$ is symmetric, we have that $\nabla V(y-x)=-\nabla V(x-y)$.
	Thus, for all $(t,x)\in[0,T)\times\R^\xdim$,
\begin{equation}\label{eq:linHJB}
	\begin{cases}
	\partial_tU
	+
	\tfrac12\Delta U
	-
	\nabla V \ast \mu^\star_t\cdot\nabla U
	-
	\tfrac12|\nabla U|^2
	+
	\int_{\R^d}\nabla V (x-y)\cdot\nabla U(t,y)\,\mu^\star_t(\d y)
	=
	\varrho(t),\\
	U(T,\cdot)=\varphi-\int\varphi\,\d\mu^\star_T .
	\end{cases}
\end{equation}
	where
	\(
	\varrho(t)\coloneqq \varrho_o(t,\mu^\star_t)-\int\mathcal A_tU(t,y)\,\mu^\star_t(\d y) -\int\nabla W\ast \mu^\star_t(y) \cdot\nabla U(t,y)\, \mu^\star_t(\d y) 
	\).
	The terminal condition follows from $u(T,\mu)(x)=\fder{}\big(\int\varphi\,\d\mu\big)(x)=\varphi(x)-\int\varphi\,\d\mu$ due to the normalization. 
	For $\psi$ as in the statement, note that $\nabla\psi=-\nabla U$, $\partial_t\psi=-\partial_tU-\varrho(t)$, and $\psi(T,\cdot)=-\varphi$, which is 
	\eqref{eq:BCGL17}.	
%	 The Fokker--Planck
%equation is \eqref{eq:FP} with $-\nabla U=\nabla\psi$, and the form of
%the optimal control is.
%{\color{red}since we dont show smoothness of $C(T)$, should we pick $\psi^*$ in $C^{1,2}([0,T)\times \R^\xdim)\cap C([0,T)\times \R^\xdim)$? Camilo: Yes}
\end{proof}

\begin{appendix}
\section{Appendix}

The next result is shown in \cite[Lemma 2]{djete2020mckean2} in the context of approximating weak rules in their argument.
Since this is not argued for general $\P\in\Pc(\muin)$, to keep the manuscript self-contained, we provide a short proof.

\begin{lemma}\label{lem:flowdetermined}
For $t\in[0,T]$, let %define the (continuous) stopping map
$\theta_t:\bar \Omega\longrightarrow \bar\Omega$, 
	$\theta_t(\bar {\mathbf x},\bar {\mathbf u},\bar {\mathbf w})
	\coloneqq (\hat{\mathbf x}_{t\wedge\cdot},\hat{\mathbf u}_{t\wedge\cdot},\hat{\mathbf w})$,
	and 
	\(
	\iota:\Pc(\bar \Omega)\longrightarrow
	C([0,T],\Pc(\bar \Omega)),
	\iota(\bar \mu)\coloneqq (\bar\mu\circ\theta_t^{-1})_{t\in[0,T]}
	\).
	Then $\iota$ is well defined and continuous, 
	and for every $\P\in\Pc(\muin)$,
	\(
		\bar\mu=\iota(\bar\mu_T),\; \P\as
	\)
\begin{proof}
That $\iota$ is well-defined is clear.
	To prove the continuity, recall that $\mathcal P(\bar \Omega)$ and $C([0,T],\Pc(\bar \Omega))$ are endowed with the weak and uniform topologies, respectively.
	Notice first that for fixed $\bar \mu$, $t\longmapsto\bar \mu\circ\theta_t^{-1}$ is continuous.
	Let $\bar \mu^n\longrightarrow \bar \mu$ weakly.
	Then, $\{\bar \mu^n, n\geq 1\}\cup\{\bar \mu\}$ is tight.
	Thus, $\bar\mu^n\circ\theta_t^{-1}\longrightarrow \bar\mu\circ\theta_t^{-1}$, for $t\in [0,T]$ fixed.
	By the equicontinuity of the paths on a compact $K\subseteq\bar\Omega$, we then have that $\bar\mu^n\circ\theta_t^{-1}\longrightarrow \bar\mu\circ\theta_t^{-1}$, uniformly in $t\in [0,T]$.
	That is, $\iota(\bar \mu^n)\longrightarrow \iota(\bar \mu)$ in $C([0,T],\mathcal P(\bar\Omega))$.\medskip

	Let $t\in[0,T]$. By \Cref{def.controlrule}\ref{def.controlrule.3} at $T$ and $t$, both $\bar\mu_t$ and $\bar\mu_T\circ\theta_t^{-1}$ are versions of $\mathcal L^{\P} ((X_{t\wedge\cdot},U_{t\wedge\cdot},W) \, |\, \overline{\Gc}_T )$, hence they are equal $\P\as$
	By continuity of both sides in $t$, $\bar\mu$ and $\iota(\bar \mu_T)=(\bar\mu_T\circ\theta_t^{-1})_{t\in [0,T]}$ agree $\P\as$
\end{proof}
\end{lemma}
\end{appendix}

%\begin{enumerate}[label=$(\roman*)$, ref=.$(\roman*)$,wide,  labelindent=0pt]
%\item 
%\end{enumerate}

% Uncomment for bibliography

\bibliography{Bibliography}

\end{document}